\documentclass[10pt]{amsart}
\usepackage{amssymb} 
\usepackage{a4, amsthm} 

\theoremstyle{plain}
\newtheorem{theorem}{Theorem}[section]
\newtheorem{lem}[theorem]{Lemma}
\newtheorem{prop}[theorem]{Proposition}
\newtheorem{cor}[theorem]{Corollary}
\newtheorem{exa}[theorem]{Example}

\theoremstyle{definition}

\newtheorem*{rem} {Remark} 
\newtheorem{remark}[theorem]{Remark} 

\font\ninerm=cmr9  \font\eightrm=cmr8  \font\sixrm=cmr6
\font\ninei=cmmi9  \font\eighti=cmmi8  \font\sixi=cmmi6
\font\ninesy=cmsy9 \font\eightsy=cmsy8 \font\sixsy=cmsy6
\font\ninebf=cmbx9 \font\eightbf=cmbx8 \font\sixbf=cmbx6
\font\nineit=cmti9 \font\eightit=cmti8 
\font\ninett=cmtt9 \font\eighttt=cmtt8 
\font\ninesl=cmsl9 \font\eightsl=cmsl8

\font\twelverm=cmr12 at 15pt
\font\twelvei=cmmi12 at 15pt
\font\twelvesy=cmsy10 at 15pt
\font\twelvebf=cmbx12 at 15pt
\font\twelveit=cmti12 at 15pt
\font\twelvett=cmtt12 at 15pt
\font\twelvesl=cmsl12 at 15pt
\font\twelvegoth=eufm10 at 15pt

\font\tengoth=eufm10  \font\ninegoth=eufm9
\font\eightgoth=eufm8 \font\sevengoth=eufm7 
\font\sixgoth=eufm6   \font\fivegoth=eufm5
\newfam\gothfam \def\goth{\fam\gothfam\tengoth} 
\textfont\gothfam=\tengoth
\scriptfont\gothfam=\sevengoth 
\scriptscriptfont\gothfam=\fivegoth

\catcode`@=11
\newskip\ttglue

\def\tenpoint{\def\rm{\fam0\tenrm}
  \textfont0=\tenrm \scriptfont0=\sevenrm
  \scriptscriptfont0\fiverm
  \textfont1=\teni \scriptfont1=\seveni
  \scriptscriptfont1\fivei 
  \textfont2=\tensy \scriptfont2=\sevensy
  \scriptscriptfont2\fivesy 
  \textfont3=\tenex \scriptfont3=\tenex
  \scriptscriptfont3\tenex 
  \textfont\itfam=\tenit\def\it{\fam\itfam\tenit}%
  \textfont\slfam=\tensl\def\sl{\fam\slfam\tensl}%
  \textfont\ttfam=\tentt\def\tt{\fam\ttfam\tentt}%
  \textfont\gothfam=\tengoth\scriptfont\gothfam=\sevengoth 
  \scriptscriptfont\gothfam=\fivegoth
  \def\goth{\fam\gothfam\tengoth}
  \textfont\bffam=\tenbf\scriptfont\bffam=\sevenbf
  \scriptscriptfont\bffam=\fivebf
  \def\bf{\fam\bffam\tenbf}%
  \tt\ttglue=.5em plus.25em minus.15em
  \normalbaselineskip=12pt \setbox\strutbox\hbox{\vrule
  height8.5pt depth3.5pt width0pt}%
  \let\big=\tenbig\normalbaselines\rm}

\def\ninepoint{\def\rm{\fam0\ninerm}
  \textfont0=\ninerm \scriptfont0=\sixrm
  \scriptscriptfont0\fiverm
  \textfont1=\ninei \scriptfont1=\sixi
  \scriptscriptfont1\fivei 
  \textfont2=\ninesy \scriptfont2=\sixsy
  \scriptscriptfont2\fivesy 
  \textfont3=\tenex \scriptfont3=\tenex
  \scriptscriptfont3\tenex 
  \textfont\itfam=\nineit\def\it{\fam\itfam\nineit}%
  \textfont\slfam=\ninesl\def\sl{\fam\slfam\ninesl}%
  \textfont\ttfam=\ninett\def\tt{\fam\ttfam\ninett}%
  \textfont\gothfam=\ninegoth\scriptfont\gothfam=\sixgoth 
  \scriptscriptfont\gothfam=\fivegoth
  \def\goth{\fam\gothfam\tengoth}
  \textfont\bffam=\ninebf\scriptfont\bffam=\sixbf
  \scriptscriptfont\bffam=\fivebf
  \def\bf{\fam\bffam\ninebf}%
  \tt\ttglue=.5em plus.25em minus.15em
  \normalbaselineskip=11pt \setbox\strutbox\hbox{\vrule
  height8pt depth3pt width0pt}%
  \let\big=\ninebig\normalbaselines\rm}

\def\eightpoint{\def\rm{\fam0\eightrm}
  \textfont0=\eightrm \scriptfont0=\sixrm
  \scriptscriptfont0\fiverm
  \textfont1=\eighti \scriptfont1=\sixi
  \scriptscriptfont1\fivei 
  \textfont2=\eightsy \scriptfont2=\sixsy
  \scriptscriptfont2\fivesy 
  \textfont3=\tenex \scriptfont3=\tenex
  \scriptscriptfont3\tenex 
  \textfont\itfam=\eightit\def\it{\fam\itfam\eightit}%
  \textfont\slfam=\eightsl\def\sl{\fam\slfam\eightsl}%
  \textfont\ttfam=\eighttt\def\tt{\fam\ttfam\eighttt}%
  \textfont\gothfam=\eightgoth\scriptfont\gothfam=\sixgoth 
  \scriptscriptfont\gothfam=\fivegoth
  \def\goth{\fam\gothfam\tengoth}
  \textfont\bffam=\eightbf\scriptfont\bffam=\sixbf
  \scriptscriptfont\bffam=\fivebf
  \def\bf{\fam\bffam\eightbf}%
  \tt\ttglue=.5em plus.25em minus.15em
  \normalbaselineskip=9pt \setbox\strutbox\hbox{\vrule
  height7pt depth2pt width0pt}%
  \let\big=\eightbig\normalbaselines\rm}

\def\twelvepoint{\def\rm{\fam0\twelverm}
  \textfont0=\twelverm\scriptfont0=\tenrm
  \scriptscriptfont0\sevenrm
  \textfont1=\twelvei\scriptfont1=\teni
  \scriptscriptfont1\seveni 
  \textfont2=\twelvesy\scriptfont2=\tensy
  \scriptscriptfont2\sevensy 
   \textfont\itfam=\twelveit\def\it{\fam\itfam\twelveit}%
  \textfont\slfam=\twelvesl\def\sl{\fam\slfam\twelvesl}%
  \textfont\ttfam=\twelvett\def\tt{\fam\ttfam\twelvett}%
  \textfont\gothfam=\twelvegoth\scriptfont\gothfam=\ninegoth 
  \scriptscriptfont\gothfam=\sevengoth
  \def\goth{\fam\gothfam\twelvegoth}
  \textfont\bffam=\twelvebf\scriptfont\bffam=\ninebf
  \scriptscriptfont\bffam=\sevenbf
  \def\bf{\fam\bffam\twelvebf}%
  \tt\ttglue=.5em plus.25em minus.15em
  \normalbaselineskip=12pt \setbox\strutbox\hbox{\vrule
  height9pt depth4pt width0pt}%
  \let\big=\twelvebig\normalbaselines\rm}

\def\Ad{{\rm Ad}}
\def\a{{\mathfrak{a}}}\def\b{{\mathfrak{b}}}
\def\c{{\mathfrak{c}}}\def\f{{\mathfrak{f}}}
\def\s{{\mathfrak{s}}}
\def\o{{\mathfrak{o}}}

\def\u{{\mathfrak{u}}}
\def\p{{\mathfrak{p}}}    
\def\k{{\mathfrak{k}}}

\def\q{{\mathfrak{q}}}
\def\g{{\mathfrak{g}}}
\def\l{{\mathfrak{l}}}
\def\h{{\mathfrak{h}}}
\def\Z{{\mathbb{Z}}}
\def\B{{\mathbb{B}}}

\def\C{{\mathbb{C}}}
\def\R{{\mathbb{R}}}
\def\HH{{\mathbb{H}}}

\def\P{{\mathbb{P}}}
\def\O{{\mathcal O}}

\def\nsmallskip{\smallskip\noindent}
\def\bbigskip{\bigskip\bigskip}
\def\nbigskip{\bigskip\noindent}

\def\nmedskip{\medskip\noindent}
\def\buildunder#1#2{\mathrel{\mathop{\kern0pt #2}
\limits_{#1}}}
\def\buildover#1#2{\buildrel#1\over#2}

\def\qq{/\kern-.185em /}
\def\half{{1\over 2}}

\def\pn{\par\noindent}
\def\sn{\smallskip\noindent}
\def\mn{\medskip\noindent}
\def\bn{\bigskip\noindent}
\def\REM #1{}

\begin{document}

\title[On univalence of Riemann domains]{On univalence of equivariant Riemann domains over the 
complexification of a non-compact, Riemannian symmetric space}

\author[Geatti]{L. Geatti}
\author[Iannuzzi]{A. Iannuzzi}

\address{Laura Geatti and Andrea Iannuzzi: Dip. di Matematica,
II Universit\`a di Roma  ``Tor Vergata", Via della Ricerca Scientifica,
I-00133 Roma, Italy}
\email{geatti@mat.uniroma2.it,
       iannuzzi@mat.uniroma2.it}
  
\thanks {\ \ {\it Mathematics Subject Classification
(2000):} 32D26, 32Q28, 53C35, 32M05}

\thanks {\ \ {\it Key words}:
Riemann domain, semisimple Lie group, symmetric space}

\bbigskip

\begin{abstract}
Let  $\, G/K\,$ be a  non-compact,  rank-one, Riemannian symmetric space
and let $\,G^\C\,$ be the universal
complexification of $\,G$.
We prove  that a holomorphically separable,
$\,G$-equivariant, Riemann domain  over $\,G^\C/K^\C\,$
is necessarily univalent, provided
that $\,G\,$ is not a covering of $\,SL(2, \R)$.
 As a consequence of the above statement  one obtains a
univalence result for holomorphically separable,
$\,G \times K$-equivariant Riemann domains over $\,G^\C$.  Here  $\,G \times K\,$ acts 
on $\,G^\C\,$ by 
left and right translations. 
The proof of  such results involves a detailed study of the $\,G$-invariant complex geometry of  the quotient $\,G^\C/K^\C$,  including a complete  classification of all its Stein $\,G$-invariant 
 subdomains.\end{abstract}
 
\maketitle

\section{Introduction}

\bigskip
Let  $\,G\,$ be a connected Lie group and let $\,Y\,$ be a complex $\,G$-manifold,
i.e.  a complex manifold endowed with a real-analytic action of $\,G\,$ by holomorphic transformations. Consider the  action of $\,G\,$
on its universal complexification $\,G^\C\,$ by left translations.
A $\,G$-equivariant local biholomorphism 

\sn
$$p\,:\, Y \, \longrightarrow \, G^\C$$ 

\sn
is by definition a 
{\it $\,G$-equivariant Riemann domain over $\,G^\C$}.
It is of interest to determine  conditions under which $\,p\,$ is injective, i.e. under which 
the Riemann 
domain is {\it univalent}.  
  
 One  motivation comes  from the classical
problem of describing the envelope of holomorphy  of a  domain
in a Stein manifold.    If $\,G =\R^n$, as a consequence of 
Bochner's tube  theorem, the envelope of holomorphy of
a $\,G$-invariant domain in $\,G^\C\,$ is a univalent $\, G$-equivariant Riemann domain over
$\,G^\C$.
An analogous statement for 
$\,G\,$ a compact Lie group is due to O. S. Rothaus  (\cite {Rt}).  
The above 
results were later generalized to arbitrary 
  holomorphically separable,
 $\,G$-equivariant Riemann domains over $\,G^\C$. They were 
 also   extended to a larger class of Lie groups, including for example 
 the product  of a compact  and a simply connected  nilpotent Lie group
 (see  \cite{CL}, \cite{Ia}, \cite{CIT}). 
  Note that since $\,G^\C\,$ is Stein (see \cite{He2}),  
 holomorphic separability of $\,Y\,$ is a necessary condition
 for   $\,p\,$ to be injective.
  
 Another motivation comes from the problem of extending 
 to a global action
 the local $\,G^\C$-action induced by a
$\,G$-action on a reduced complex space.
Indeed the univalence of $\,G$-equivariant  Riemann domains over $\,G^\C\,$ 
turns out to be
a necessary condition for the existence of such an extension
(see \cite{Pa},$\,$\cite{HeIa},$\,$\cite{CIT}).

When $\,G\,$ is a non-compact, real semisimple Lie group,
univalence of holomorphically, separable $\,G$-equivariant 
Riemann domains over $\,G^\C\,$ does not hold in general. 
For $\,G=SL(2,\R)$, 
a Stein counter-example was pointed out to us by K. Oeljeklaus
(see Sect. 8).
 The image of this Riemann domain in  $\,G^\C\,$
 is also invariant under right  $\,K$-translations and its construction is based on the existence
of non-trivial coverings  of the $\,K$-orbits in $\,G^\C$. Here $\,K\,$ a
maximal compact subgroup in $\,G$.
Observe  also that $\,SL(2,\C)/SL(2,\R)\,$ is simply connected.  Thus this  
example gives a negative answer to the  question whether the
simple-connectivity of the quotient $\,G^\C/G\,$ 
is a sufficient condition for  univalence 
of $G$-equivariant Riemann domains over $G^\C$
(cf.~\cite{CL}).

Let $\,G \,$ be a connected, non-compact, real simple Lie group and  
let $\,K\,$ be a maximal compact subgroup of $\,G$.
The group $\,G\,$ is not necessarily embedded in $\,G^\C$, but it is 
 assumed to have finite center.  Consider the action of the
product group $\,G\times K\,$ on $\,G^\C$ by left and right translations. One of the 
results of  this paper is the following theorem (Thm.~\ref{UNIVALENCE3}).

\mn
{\bf Theorem.} {\sl Let $\,G/K\,$ be a non-compact, rank-one, Riemannian symmetric space.
A holomorphically separable,
$\,G\times K$-equivariant Riemann domain  $\,p \colon Y  \to G^\C\,$ is
univalent, provided that $\,G\,$ is not a covering of $\,SL_2(\R)$.}

\nmedskip

Note that since $\,Y\,$  embeds equivariantly  into its envelope of holomorphy 
(cf.~\cite{Ro} and Sect.~2), there is no loss of generality in assuming that
$\,Y\,$ is Stein. Then  a  result of P. Heinzner (\cite{He1}) implies that
the categorical quotient $\,Y\qq K\,$  is  also Stein.
By performing categorical $\,K$-reduction on both $\,Y\,$ and
$\,G^\C$, one can associate to  $\,p \colon Y  \to G^\C\,$ a Stein, 
$\,G$-equivariant Riemann domain $\,q \colon Y \qq K \to G^\C/K^\C$.
A suitable  characterization of the univalence of $\,q\,$ (see 
Prop.~\ref{UNIVALENCE1})
 implies that $\,p\,$ is univalent if  $\,q\,$ is univalent.
Then  the above  theorem 
is a consequence of the following  one, which is the main result  of the paper
(Thm.~\ref{UNIVALENCE2} and Rem.~\ref{FINITECENTER}).

\mn
{\bf Theorem.}  {\sl A holomorphically separable, $\,G$-equivariant
Riemann domain  $\,q \colon \Sigma \to G^\C/K^\C\,$ is univalent,  provided that $\,G\,$ is not a covering of $\,SL_2(\R)$.}

\mn

The  proof of  this theorem is carried out as follows.
First we show that, 
with few exceptions, the map $\,q\,$ 
is  injective on every $\,G$-orbit. For   principal $\,G$-orbits this is done by determining their topology. The result is then extended to the remaining  $\,G$-orbits by a general argument (Sect.~5). 
As a consequence there exists a $\,G$-invariant domain in
$\,\Sigma\,$ on which $\,q\,$ is injective.

Next we show that such domain can be enlarged
to the whole $\,\Sigma\,$.
This is done  by successively 
 lifting to $\Sigma$ 
local slices for principal $G$-orbits in $\,G^\C/K^\C$. Since such slices are one-dimensional and  $q$ is injective on $G$-orbits, each  lifting  determines a $G$-invariant domain   in $\Sigma$ on which $q$ is injective.  
The main difficulty 
 is to ensure monodromy around singular $G$-orbits.
For this we combine  a detailed
description
of the $\,G$-orbit   structure of $\,G^\C/K^\C$  with the  
 complex-geometric properties of certain
non-Stein, $\, G$-invariant  domains in $\,G^\C/K^\C$.

By the above univalence result, all Stein, $\,G$-equivariant
Riemann domains over $\,G^\C/K^\C\,$ can be regarded as Stein,
invariant  domains in $\,G^\C/K^\C$. We carry out their 
classification in  Theorem \ref{STEIN}.

For $\,G/K\,$ of arbitrary rank, 
recent investigations due to several authors have indicated
  an interplay between  the  complex geometry  of distinguished Stein, $\,G$-invariant domains in  $\,G^\C/K^\C\,$ (see \cite{KS}, \cite{FHW} and references therein) and the harmonic analysis on  $G$-symmetric spaces contained in $G^\C/K^\C$.
Envelopes of holomorphy of $\,G$-invariant  domains in $\,G^\C/K^\C\, $
might give new insights
in this context. We hope the present paper to be a first step for
 for further investigations  on symmetric spaces of  higher rank. 

\bn
The paper is organized as follows. 

\sn
In Section 2 we recall some basic notions and  results from geometric invariant
theory.

\sn
In Section 3, from a Stein 
$\,G\times K$-equivariant Riemann domain  $\,p\colon Y\to
G^\C\,$ we obtain a Stein, $\,G$-equivariant, 
Riemann domain $\,q\colon Y \qq K \to G^\C/K^\C$. 
We also show  that $\,p\,$ is univalent if  $\,q\,$ is univalent.

\sn
In Section 4 we give a detailed description of the $\,G$-orbit structure 
of $\,G^\C/K^\C\,$ when $\,G/K\,$ is a  non-compact,
rank-one, Riemannian symmetric space.
We also describe an explicit model for the space $\,G^\C/K^\C\,$ in the cases
  $\,G=SO_0(n,1)\,$ and $\,G=SU(n,1)$.

\sn
In Section 5 we show that, with few exceptions,
a $\,G$-equivariant Riemann domain $\,q\colon Y\to G^\C/K^\C\,$ is  univalent
on every $\,G$-orbit.

\sn
In Section 6
we carry out a complete classification of  Stein,
$\,G$-invariant domains in
$\,G^\C/K^\C$.  When $\,G=SU(n,1)\,$ some  of these domains appear to be new.

\sn
In Section 7 we prove the univalence  result  for holomorphically separable,
$\,G$-equivariant Riemann domains over $\,G^\C/K^\C$. 

\sn
In Section 8 we obtain the result for holomorphically separable,
$\,G\times
K$-equivariant Riemann domains over $\,G^\C$.
We also  discuss some examples. 

\sn
In the Appendix we compute the Levi form of 
all non-closed hypersurface $\,G$-orbits in
$\,G^\C/K^\C$. The results of this computation are used in Sections 6
and 7.

\nbigskip
{\it Acknowledgments.} We wish to thank Karl Oeljeklaus for pointing out to 
us  Example \ref{KLAUS}.
We also thank Stefano Trapani
for suggesting the argument in Lemma \ref{KILLtheCENTER}.

\hfill \vfill 

\eject


\section{Preliminaries}

\bigskip
Let $\,G\,$ be a  connected, real  Lie group. A complex Lie group
$\,G^\C\,$ together with a Lie group homomorphism $\,\iota:G \to G^\C\,$
is called a {\it universal complexification} of $\,G\,$ if for every
Lie group homomorphism $\,\psi\,$ from $\,G\,$ to a complex Lie group
$\,H\,$ there exists a holomorphic homomorphism $\,\psi^\C:
G^\C \to H\,$ such that $\,\psi = \psi^\C \circ \iota$. 
A universal complexification of $\,G\,$ always exists and
is unique up to biholomorphisms (see \cite{Ho}). 

Assume  that $\,G\,$ is a  connected, real semisimple Lie group. 
Denote by $\g$ the Lie algebra of $G$ and by $\g^\C:=\g \oplus i\g$ its complexification. Then  the universal complexification of $G$ is a complex semisimple Lie group $G^\C$  with Lie algebra $\g^\C$.
If $G$ is a  real form of a simply connected complex semisimple Lie group $G^\C$, then its  universal complexification is   $\,G^\C$.  
Furthermore, if  $\,\Gamma\,$ is a central subgroup of $\,G\,$,
then the universal complexification of the quotient  group $G/\Gamma$ 
is  given by $G^\C/\Gamma$.
Note that every automorphism of   $\,G\,$  uniquely  extends to a holomorphic automorphism of  its universal complexification $\,G^\C$.

Let $\,K\,$ be a compact Lie group and $\,X\,$ a Stein $K$-space, i.e.  a reduced
Stein space with a real-analytic action of $\,K\,$ by holomorphic transformations. The 
{\it categorical quotient} $\,X\qq K\,$ of $\,X\,$  is defined by the following
 equivalence relation:  $\,x\sim y\,$ if and only if $\,f(x)=f(y)\,$ for every $\,K$-invariant
 holomorphic function $\,f\,$ on $X$.
We recall some basic properties of the categorical quotient (see \cite{He1}).
\pn

\bigskip
\begin{theorem}
\label{CATEGO}
Let $K$ be a compact Lie group and $X$ a Stein $K$-space. Then

\begin{itemize}
\item[$(i)$] the categorical quotient $X\qq K$ equipped with the 
algebra $\, {\mathcal O}(X)^K\,$ of  holomorphic $\,K$-invariant functions on $X$ 
is a Stein space and the projection
$\,\pi \colon X\rightarrow X\qq K\,$ is holomorphic,

\pn
\item[$(ii)$] for every $\,K$-invariant
holomorphic map $\,\psi\,$ from $\,X\,$ to a complex space 
$\,Y\,$  there exists a unique holomorphic map
 $\widehat   \psi \colon X\qq K\rightarrow Y$ making   the diagram
 
$$\def\mapright#1{\smash{
\mathop{\longrightarrow}\limits^{#1}}}
\def\mapdown#1{
\rlap{$\vcenter{\hbox{$\scriptstyle#1$}}$}{\ \ \big\downarrow \ \, \ }}
\begin{matrix}X &\mapright{\pi}&X \qq K\cr
\cr 
\mapdown{\psi}&\swarrow\rlap
{$\vcenter{\hbox{$\scriptstyle\widehat \psi$}}$}\cr
\cr
Y\cr\end{matrix}$$
 
 commute.
 \end{itemize}
\end{theorem}

\bigskip
If  the $\,K$-action on $X$ is the restriction of a $K^\C$-action, then the algebras of $K$-invariant  and  of $\,K^\C$-invariant holomorphic functions on $X$  coincide.
In particular they induce
the same equivalence relation on $X$ and  $\,X \qq K \cong X \qq K^\C.$
In this case, if all $K^\C$-orbits are closed, then $\,X \qq K^\C\,$ coincides
with the
usual orbit  space  $\,X/K^\C\,$ (cf.~\cite{Sn},~Thm.~3.8).
A $K$-action
on a Stein  space  can always be extended to a $K^\C$-action, as shown  by the following
theorem due to Heinzner (\cite{He1}).

\bn
\begin{theorem}
\label{GLOBALI}
Let $\,K\,$ be a compact Lie group and $\,X\,$ a Stein $\,K$-space. Then there exist  a 
Stein $\,K^\C$-space $\,X^*\,$ and a $\,K$-equivariant holomorphic embedding
$\,\iota\colon X\hookrightarrow X^*$ with the following properties:

\smallskip
\item{(i)} the map $\iota$ is open and $\,K^\C\cdot \iota(X) = X^*$,

\smallskip
\item{(ii)} for every $K$-equivariant holomorphic map $\varphi$ from $X$ into a complex $K^\C$-space $Z$, there exists a unique $K^\C$-equivariant holomorphic map 
$\varphi^*\colon X^*\rightarrow Z$ making the  diagram
$$\def\mapright#1{\smash{
\mathop{\hookrightarrow}\limits^{#1}}}
\def\mapdown#1{
\rlap{$\vcenter{\hbox{$\scriptstyle#1$}}$}{\ \ \big\downarrow \ \, \ }}
\begin{matrix}X &\mapright{\iota}&X^*\cr 
\cr
\mapdown{\varphi}&\swarrow\rlap
{$\vcenter{\hbox{$\scriptstyle\varphi^*$}}$}\cr
\cr
Z\cr\end{matrix}$$
commute,

\smallskip
\item{(iii)} the inclusion $X\hookrightarrow X^*$ induces an isomorphism between the categorical quotients $X\qq K$ and $X^*\qq K^\C$.

\end{theorem}

\bigskip
Observe that, since $\,K^\C\cdot \iota(X) = X^*$, if $\,X\,$ is non-singular,
then $\,X^*\,$ is also non singular. 
Let $\,X\,$ be a complex manifold and let $\,G\,$ be a Lie group.
A  Riemann domain over $\,X\,$ is a complex manifold
$\,Y\,$ together with a locally biholomorphic map $\,p:Y \to X$.
If both $\,X\,$ and $\,Y\,$ are $\,G$-manifolds and the map $\,p\,$ is 
$\,G$-equivariant, then we refer to $\,p:Y \to X\,$ as a
{\it $\,G$-equivariant Riemann domain}. 
If $X$ is Stein and  $\,Y\,$ is  holomorphically separable, then $\,Y\,$
embeds as an open domain in its
envelope of holomorphy $\, \widehat Y \,$ and the map $\,p\,$ extends to
a local biholomorphism $\,\widehat p: \widehat Y \to X\,$ (see \cite{Ro}).
Moreover the $\,G$-action on $\,Y\,$ extends to
 $\,\widehat Y\,$ and the map $\,\widehat p\,$ is
$\,G$-equivariant, i.e.   $\,\widehat p\colon \widehat Y\to X\,$ is a
Stein, $G$-equivariant Riemann domain. 
 
A Riemann domain $\,p:Y \to X$  is called {\it univalent} if the map $p$ is injective. 
Assume that  $X$ is Stein and $Y$ is holomorphically separable.  If $\,\widehat p\,$
is univalent, then $\,p\,$ is also univalent.
 Aiming at univalence results for holomorphically separable Riemann domains
 over $\,G^\C$, it is therefore not restrictive to start with  Riemann domains which are Stein.

\bigskip


\section{From Riemann domains over $G^\C$ to Riemann domains over $G^\C/K^\C$}

\bigskip

Let $G$ be a connected,  non-compact, real semisimple Lie group,  $K \subset G$  a maximal compact subgroup  and  $G^\C$ the universal complexification of $G$. 
Let $G\times K$ act on $G^\C$ by left and right translations, i.e.
$$(g,k)\cdot z:= gzk^{-1},\quad {\rm for}\quad (g,k)\in G\times  K,~ z\in G^\C.$$
In this section, to every Stein, $G\times K$-equivariant Riemann domain $p\colon Y\to G^\C$ we associate 
  a Stein,  $G$-equivariant Riemann domain $q\colon \Sigma\to G^\C/K^\C$.  
  We also show
   that the univalence 
of $q$ implies that of $p$.

\bn

Let $\,X\,$ be a Stein $\,K^\C$-manifold and let
$\,p\colon Y\to X\,$ be a Stein, $\,K$-equivariant Riemann domain.
By Theorem \ref{GLOBALI} there exist a Stein $K^\C$-manifold $Y^*$, a $K$-equivariant holomorphic open embedding  $\iota \colon Y\hookrightarrow Y^*$ and a $K^\C$-equivariant holomorphic map $ p^*\colon Y^*\to X$ such that the  diagram 
 $$\def\normalbaselines{\baselineskip20pt\lineskip3pt
\lineskiplimit3pt} 
\def\mapright#1{\smash{
\mathop{\hookrightarrow}\limits^{#1}}}
\def\mapdown#1{
\rlap{$\vcenter{\hbox{$\scriptstyle#1$}}$}{\ \ \big\downarrow \ \, \ }}
\begin{matrix}Y &\mapright{\iota}&Y^*\cr 
\cr
\mapdown{p}&\swarrow\rlap
{$\vcenter{\hbox{$\scriptstyle p^*$}}$}\cr
\cr
X\cr\end{matrix} $$
commutes.
Since $p$ is locally biholomorphic, $p^*$ is $K^\C$-equivariant and $Y^*=K^\C\cdot Y$,
one has that
  $p^*$ is   locally biholomorphic as well, i.e.
it defines  a Stein  $K^\C$-equivariant Riemann domain.
By Theorem  \ref{GLOBALI}  the spaces $\,Y^*\qq K^\C\,$ and $\, Y \qq K\,$ 
 are biholomorphic. Therefore Theorem  \ref{CATEGO} implies 
there exists  a holomorphic map $q\colon Y \qq K\to X\qq K^\C$ making the   diagram 
\smallskip
$$
\begin{matrix}
 Y^*  &   \longrightarrow \ \ \ \,   &    Y^*\qq K^\C \cong Y \qq K  \cr 
    \cr
     \quad  \,\downarrow  \    p^*   &    &   \downarrow \, q   \  \  &   \cr
                              \cr 
                    X \   &  \longrightarrow   & X\qq K^\C\  \ &       \cr
\end{matrix}$$

\nsmallskip
 commute. 
Here the horizontal arrows
denote the categorical quotient maps.

Assume that all $\,K^\C$-orbits in $\,X\,$ are closed.
We claim that all $\,K^\C$-orbits in $\,Y^*\,$ are closed as well.
Suppose by contradiction that there exists a non closed orbit
$\,K^\C \cdot y\,$ in $\,Y^*$. Let $\,K^\C \cdot z\,$ be a lower dimensional
orbit in its closure (see \cite{Sn}, Prop.$\,$2.3).
Since $\,p^*\,$ is locally biholomorphic and $\,K^\C$-equivariant,  the
orbit $\,K^\C \cdot p^*(z)\,$ lies in the closure of $\,K^\C \cdot p^*(y)\,$
and has lower dimension. In particular such orbits are distinct. It follows that
 the orbit $\,K^\C \cdot p^*(y)\,$ is not closed,  contradicting the assumption.

By the above claim,  the categorical quotients $\,X\qq K^\C\,$ and $\,Y^*\qq K^\C\,$
coincide  with  the orbit spaces $\,X/K^\C\,$ and $\,Y^*/K^\C\,$, respectively. 
As a consequence, the map 
$$q \colon Y\qq K \to X/K^\C$$

\nsmallskip
is locally biholomorphic, i.e. it defines a Stein  Riemann domain.
We refer to it as the {\it Riemann domain induced by $\,p\colon Y\to X$}. 
Next we prove  a general  univalence result for Stein, $K$-equivariant Riemann domains.


\nbigskip
\begin{prop}
\label{UNIVALENCE1}
Let $X$  be a Stein $\,K^\C$-manifold, all of whose  $\,K^\C$-orbits  are closed and have connected isotropy subgroups. Let $\,p\colon Y\to X\,$ be a Stein, $\,K$-equivariant Riemann domain and  $\,p^*\colon Y^*\to X\,$ its extension to the $\,K^\C$-globalization $\,Y^*\,$ of $\,Y\,$.  Then 

\begin{enumerate}
\smallskip
\item[$(i)$] the induced Stein, Riemann domain  $\,q:   Y \qq K \to  X/K^\C\,$ is univalent  if and only if 
$\,p^* :  Y^*  \to  X\,$ is~univalent,
 
 \smallskip
\item[$(ii)$] if $\,q:   Y \qq K  \to  X/K^\C\,$ is univalent, then $p\colon Y\to X$ is univalent.
\end{enumerate}
\end{prop}

\smallskip
\begin{proof}  
 $\,(i)\,$ If $p^*$ is injective, then it maps distinct $\,K^\C$-orbits in $\,Y^*\,$ onto distinct $K^\C$-orbits in $X$. As we already noticed,  since all $\,K^\C$-orbits in $\,X\,$
 are closed,
 the categorical quotients $\,X\qq K^\C\,$ and $\,Y^*\qq K^\C\,$ coincide  with  the orbit spaces $\,X/K^\C\,$ and $\,Y^*/K^\C\,$, respectively.
It follows that the induced map  
$ \, Y^*/ K^\C\longrightarrow X/ K^\C \,$ is injective. 
Moreover, by Theorem \ref{GLOBALI}, the space
  $\,Y\qq K\,$ can  be identified with $\,Y^*\qq K^\C$. As a result the induced
  Riemann domain
$\,q\colon Y \qq K \longrightarrow X/ K^\C\,$ is univalent.

Conversely, assume that $\,q\colon Y \qq K \longrightarrow X/ K^\C\,$ is univalent, i.e. that the map $ \,Y^*/ K^\C\longrightarrow X/ K^\C\,$ is injective.
By assumption the $\,K^\C$-isotropy subgroups in $\,X\,$ are connected,
thus    $\,p^*\,$ is  injective on every $\,K^\C$-orbit in $\,Y^*$. It follows  
that $\,p^* \colon Y^*\longrightarrow X\,$ is globally injective.
This concludes the proof of $\,(i)$. 
Statement $\,(ii)\,$ is a direct consequence of $\,(i)$.
\end{proof}


\bn

\begin{remark}
$\,${\rm In general, under the assumptions of the above proposition,  the univalence of  $\,  p:  Y \to X\,$
does not  imply  that of $\,q\colon  Y \qq K \to  X/K^\C$.
For instance, let $\,\C^*\,$ act on $\,\C \times \C^*\,$ and 
on $\,X := \C^* \times \C^*\,$ by multiplication on the second factor.
Define $\,p^*:\C \times \C^* \to X\,$ by 
$\,(z,w) \to (e^{ z}, w)\,$ and consider
$$ Y:= \{\,(z,w) \in \C \times \C^*  \,:\, {\rm Im }\, z < |w| < 2 \pi + {\rm Im} \,z  \,\}.$$
\nsmallskip
Then $\,Y\,$ is a  Stein $\,S^1\,$-invariant  subdomain of $\,Y^*=\C \times \C^*\,$  and the map $\,p:=p^*|_Y\,$ is injective.
Nevertheless the induced map $\,q:Y \qq S^1 \cong \C \to X/\C^* \cong \C^*$,  given 
by $\, z \to e^{  z}$, is not injective.\qed
 }
\end{remark}

\bn

Consider now the case when $\,X \,$ is the group $\,G^\C\,$ endowed
with the $\,G\times K$-action by left and right translations.
Let  $\,p\colon Y\to G^\C\,$ be a Stein, $\,G\times K$-equivariant Riemann domain. Note that  the actions of $\,G\,$ and $\,K\,$  commute on $\,G^\C\,$.
Thus  they also
commute on $\,Y$, due to the fact that  $\,p\,$  is equivariant
and locally injective.
Since the $\,K$-action on $\,G^\C\,$ is the 
restriction of a $\,K^\C$-action all of whose  orbits  are closed,
the spaces  $\,G^\C\qq K\, $ and  $\,G^\C/K^\C\,$ are biholomorphic.

By the universality property of the categorical quotient (cf.~Theorem~\ref{CATEGO}),
the $\,G$-actions on $\,Y\,$ and on $\,G^\C\,$ induce
$\,G$-actions on $\,Y \qq K\,$ and on $\,G^\C/K^\C$, respectively.
Moreover the induced Stein, Riemann domain 
$$q \colon Y\qq K \to G^\C/K^\C$$
is $G$-equivariant. 
By applying Proposition \ref{UNIVALENCE1} to this situation,
one obtains the following result.


\bigskip
\begin{cor}
\label{INI}$\ $ Let $p\colon Y\to  G^\C$ be a Stein,  $\,G\times K$-equivariant Riemann domain over $G^\C$ and let $\,q:  Y \qq K \to  G^\C/K^\C\,$  be the induced  Stein,
$\,G$-equivariant  Riemann domain over $G^\C/K^\C$. If $q$ is univalent, then
$\,p\,$ is univalent.
\end{cor}

 
\bigskip
\section{$G $-orbit structure of $G^\C/K^\C$}

 \bigskip
 
Let $G$ be a connected,  non-compact, real simple Lie group, $K\subset G$ a maximal compact subgroup and $G^\C$ the universal complexification of $G$. Assume that $G$ is embedded in $G^\C$. The quotient  $G/K$ is a  Riemannian symmetric space of the
non-compact type.  In this section we obtain  a complete description of the $G$-orbit structure of $G^\C/K^\C$ in the case when $G/K$ has rank one.

\sn
We recall some basic facts which hold for $G/K$ of arbitrary rank.
Let 
$\, \sigma\,$ denote the anti-holomorphic involution of
$\,G^\C\,$ relative to
$\,G\,$ and $\,\tau: G^\C \to G^\C\,$  the holomorphic extension
of the Cartan involution $\,\theta\,$ of $\,G\,$ with respect to $\,K$.
Note that the fixed point set of $\,\tau\,$ in $\,G^\C\,$
contains the complexification $\,K^\C\,$  of $\,K$. The commuting
composition $\, \sigma \circ \tau = 
\tau \circ \sigma\,$ is a  Cartan involution 
of $\,G^\C$. Denote by $U$ the corresponding compact real form. The  $U$-orbit of the base point $eK^\C$ in $G^\C/K^\C$  is diffeomorphic to the compact dual symmetric space $U/K$, and is embedded  in $G^\C/K^\C$ transversally to $G/K$.

\bn
\begin{remark}\label{QUOTOP}
\item{(i)} 
For every triple $\, (G, K, G^\C)\,$ as above, 
the manifold $\,G^\C/K^\C\,$ is simply connected.  To see this, denote by
$\,\widetilde G^\C\,$  and $\,\widetilde U \subset \widetilde G^\C\,$ the universal coverings of $\,G^\C\,$ and $\,U\,$
respectively. Let  $\,\widehat G\,$ be the real form of $\,\widetilde G^\C\,$ relative
to the lifting  of $\,\sigma\,$ to $  \,\widetilde G^\C\,$. The group $\,\widehat G\,$ is
connected (cf. \cite{Sb})  and is a finite covering of $\,G$. Hence  $\,G=\widehat
G/\Gamma$, where  $\,\Gamma\,$ is a finite central subgroup of $\,\widehat G$.
Similarly
$\,K=\widehat K/\Gamma$, where $\,\widehat K\,$ is a maximal compact subgroup of
$\,\widehat G\,$. One has
$\,G^\C\cong \widetilde G^\C/\Gamma\,$ (cf. Sect.$\,$2) and consequently
$\,U =\widetilde U/\Gamma$. As a consequence
there are isomorphisms 
$$U/K \cong \widetilde U/\Gamma/ \widehat K/\Gamma\cong \widetilde U/\widehat
K.$$ Since $\,\widehat K\,$ is connected, the quotient $\,\widetilde U/\widehat K\,$ is simply connected.  Moreover $ \, U/  K\, $ is a
topological retract of $\, G^\C/K^\C$.  Hence the claim follows.

\sn
\item{(ii)} 
From different triples $\, (G, K, G^\C)\,$  as above associated with the same Riemannian symmetric space one obtains  the same complexification $\,G^\C/K^\C$. Indeed the map
 $\,\widetilde G^\C/ \widehat K^\C  \to G^\C/ K^\C$, given by $\,g\widehat K^\C
 \to g \Gamma K^\C$, defines a biholomorphism.
Moreover the center of $\,G\,$ acts trivially on $\,G^\C/K^\C$. As a consequence,  
different triples $\, (G, K, G^\C)\,$ yield
the same $\,G$-orbit structure of $\,G^\C/K^\C\,$ and  $\,G$-equivariantly diffeomorphic orbits.
    \qed
\end{remark}

\bn
 
Closed $G$-orbits of maximal dimension form an open dense subset of 
$\,G^\C/K^\C\,$ and come in a finite number of orbit types.  We refer to them as 
{\it principal $\,G$-orbits}.  They have real codimension   equal to the rank
of $G/K$. {\it Singular orbits} are closed $G$-orbits which are not
principal.

The $\,G$-orbit structure of $\,G^\C/K^\C\,$
is closely related to the $\,G\times K^\C$-orbit
structure of $\,G^\C$. Then, slices for the closed $\,G$-orbits in~$\,G^\C/K^\C$ can be obtained 
by applying Matsuki's  results on double coset decompositions 
of groups arising from two involutions  (\cite{Ma}, Sect.$\,$4).

\bn

Let $\,\k\oplus\p\,$ be the Cartan decomposition of $\,\g\,$ with respect to $\,K$  
  and let $\,\a \,$ be a maximal abelian subspace of $\, \p$. Following Matsuki,
 we denote by $$\,\mathcal A:=\exp i\a\, K^\C \,$$ the image of the compact torus $\exp i\a$ in $G^\C/K^\C$.
The set $\, \mathcal A\,$ is a slice for all  closed $\,G$-orbits intersecting the compact dual symmetric space $\,  U/K$ in $G^\C/K^\C$.  It is called the  {\it fundamental Cartan subset}.
 The remaining slices for closed $\,G$-orbits in $G^\C / K^\C$ 
are of the form
$$\,\mathcal C:=\exp i \c \cdot z,$$ where 
$\,\c\,$ is an abelian semisimple subalgebra of $\,\g\,$ of the same dimension
as $\,\a$ and  $\,z \in \mathcal A\ $ is a {\it base point }
sitting on a  singular closed $\,G$-orbit. 
 Such sets $\,\mathcal C$ are called {\it standard Cartan subsets}. 
 
 By \cite{Ge2},   every standard Cartan subset $\mathcal C\,$ admits a base point $z$ with the following properties: 

{\it \sn
$\bullet$ there exists a subgroup $G'\subseteq G$ (possibly $G'$ is equal to $G$)  such that the isotropy subgroup  of $z$ in $G'$ coincides with the isotropy subgroup $G_z$  of $z$  in $G$. 

\sn
$\bullet$ the quotient $G' / G_z$ is a  pseudo-Riemannian symmetric space  of the same rank as~$G/K$,

\sn
$\bullet$ the slice representation of $\,G_z\,$  at $z$  is equivalent to the
isotropy representation of~$G' / G_z$.} 

\bn 

 More precisely,  let $\g'=\g_z\oplus \q'$ be the decomposition of the Lie algebra of $G'$  corresponding to the symmetric space $\,G'/G_z\,$ (when $ G'=G$, one has $\g=\g_z\oplus \q$).
Denote by $T(G\cdot z)_z$ the tangent space  of the orbit $\,G\cdot z\,$ at $\,z\,$
and by $N_z$   a complementary subspace of $T(G\cdot z)_z$ in $T(G^\C/K^\C)_z$. 
 Then $N_z\cong \q'$ 
 and  the  slice representation at $z$  is equivalent to the Adjoint representation of $G_z$  on~$\q'$.  
Moreover, both  $\a$ and $\c$ are maximal abelian subalgebras in $\q'$. 

Consider the twisted bundle $\,G \times _{G_z} \q'\,$  defined as the
orbit space of
$\,G \times \q'\,$ under the action of  $\,G_z\,$ given  by 
$\, h \cdot(g,X):=(gh^{-1},  Ad_h X)$. The group
$G$ acts on $\,G \times _{G_z}
\q'\,$ by 
$\,\widehat g\cdot [g,X]:=[\widehat g g,X]$. By Luna's slice Theorem
(\cite{Lu},
Prop.$\,$1.2), there exists an open
$\,Ad_{G_z}$-invariant neighborhood $\,V\,$ of $\,0\,$
in $\,\q' \,$ such that the map
\begin{equation} \label{LUNA}
G \times _{G_z} V\to G^\C/K^\C, \quad \quad [g,X] \to g \exp iX \cdot z \, 
\end{equation} 
is a $\,G$-equivariant diffeomorphism onto an open 
$\,G$-invariant  neighborhood of $\,z$ in~$G^\C/K^\C$.  
Non-closed $G$-orbits in $G \times _{G_z} V$ correspond to
non-closed
$Ad_{G_z}$-orbits in $V$. The standard Cartan subset $\,\mathcal
C\,$  in~$G^\C/K^\C$ is the image of the set  $\,\{e\} \times \c\,$
 via the above map.
 
 \smallskip
 Let us now assume that $\,G/K\,$ has rank one. Then 
 the $\,G$-orbit space  of $\,G^\C / K^\C\,$ can be completely determined.
Let  $\, \Delta_\a\,$ be the restricted root system of $\g$ with respect to 
$\,\a\,$ and  let
$$\,\g = Z_\g(\a) \bigoplus_{\alpha \in \Delta_\a} \g^\alpha\,, \quad
{\rm with } \quad  Z_\g(\a)=Z_\k(\a)\oplus \a\, ,  $$
  be the corresponding root decomposition. Here $Z_\g(\a)$ and $Z_\k(\a)$  denote  the centralizers of $\a$ in $\g$ and $\k$, respectively.
Let $\,\Gamma\,$ be the lattice in $\,\a\,$  given by
the kernel of the map $\,\a \to U/K \,$ defined  by $\, X \to \exp (iX) K $.
Since the symmetric space $U/K$ is  simply connected (cf.
Remark$\,$\ref{QUOTOP}), the lattice $\Gamma$ is given by 
$$\Gamma =\bigoplus_{\alpha \in \Delta_\a} \Z i \pi h_\alpha,$$
where $h_\alpha\in\a$ is uniquely determined by $ \alpha(h_\alpha)=2$
(cf.$\,$\cite{Hl}, Thm.$\,$8.5, p.$\,$322).  Denote by $\,W_K(\a)\,$  the  Weyl
group    of
$\,\a\,$ and let  the semidirect  product $\,W_K(\a) \ltimes \Gamma$ act on
$\,\a\,$ by 
$$\,(k, \gamma) \cdot A := \Ad_k A + \gamma\,.$$
Denote by $\,\a_0\,$  a  fundamental domain  for this action  and 
define $\,\mathcal A_0:=
\exp i\a_0K^\C $. 
Then every closed $G$-orbit through the fundamental Cartan subset
$\,\mathcal A   $ intersects $\,\mathcal A_0$ in a single point
(cf.~\cite{Ma},~ Thm.~3).

Let $\,z \in \mathcal A_0\,$  be a  base point for a standard Cartan subset $\mathcal C$. By \cite{Ge2} and by the local linearization (1), the $G$-orbit structure of $G^\C/K^\C$ in a neighbourhood of $z$ is modelled on the orbit structure of the 
tangent space of a   rank-one, pseudo-Riemannian symmetric space under the isotropy representation. It  can be described as follows.

\bn
\begin{remark} 
\label{ADJOINT} Let $G/H$ a rank-one, pseudo-Riemannian symmetric space. Assume 
that the group $H$ is connected.   Let $\g=\h\oplus \q$ be the corresponding Lie algebra decomposition and $ \,\q\cap\k\, \oplus\, \q\cap \p\, $ the Cartan decomposition of $\q$.  The isotropy representation of $G/H$ is equivalent to the Adjoint representation of $H$ on $\q$.
Denote by $B$ both the Killing form of $\g$ and  its restriction to $\q\setminus\{0\}$.  The signature of
$\,B \,$ on $\q$ is given by $\,(s^+, s^-)$, with 
$$\,s^+ := \dim (\q\cap \p)\,\qquad \,s^- := \dim (\q\cap \k).$$ 
For $\,r \in \R$, denote by  
$\,B_r\,$ the level hypersurface $\,\{ B =r\}\,$ 
 in $\q\setminus\{0\}$.
 In diagonalized
form  one has
$$\, B_r=\{x_1^2 + \dots+ x_{s^+}^2 -y_1^2 - \dots - y_{s^-}^2=r\}\,.$$  
Since $\,G/K\,$ has rank one, every $\,Ad_H$-orbit in $\q\setminus\{0\}$ is a hypersurface. Thus, 
by the connectedness of $\,H\,$ and that $Ad_H$-invariance of $B$, it  coincides with a connected component of some $\,B_r$. 
We distinguish four cases.

\begin{itemize}

\smallskip
\item[$(a)$] Assume $\,s^+=s^-=1.$  For every $\,r\not=0\,$ the level set
$\,B_r\,$ consists of two connected  components. They  intersect either 
$\,\a = \q \cap \p$ or 
$\,\c = \q \cap\k$ in opposite points, depending on whether $\, r>0$ or $\,r<0$. The nilcone $\,B_0$ consists of four non-closed $\,Ad_H$-orbits.

\smallskip
\item[$(b)$] Assume $\,s^+>1$ and $\,s^-=1.$  For $\,r>0\,$ the level set
$\,B_r\,$ consists of a single component intersecting 
$\,\q \cap \p\,$ in a sphere. Thus for every non-zero vector $A\in \q\cap \p$,  and  every  $\,\,t>0\,$, the   points 
$\,tA $ and  $\,-tA \,$ belong to the same 
$\,Ad_H$-orbit. 
If $\,r<0\,$ the level set $\,B_r\,$ consists of two connected  components, which  intersect $\,\c = \q \cap\k$ in  opposite   points. 
  The  nilcone $\,B_0$ consists of two non-closed $\,Ad_H$-orbits. 
  
\smallskip
\item[$(b)$] Assume $\,s^+=1$ and $\,s^->1.$  If $\,r>0\,$, the level set $\,B_r\,$ consists of two connected  components, which  intersect  $\,\a = \q \cap\p$ in  opposite   points.   If $\,r<0\,$, the level set $\,B_r\,$ intersects 
$\,\q \cap \k\,$ in a sphere. Thus for every non-zero vector $C\in \q\cap \k$   and every  $\,\,s>0\,$, the    points 
$\,sC$ and  $\,-sC\,$ belong to the same 
$\,Ad_H$-orbit. 
 The  nilcone $\,B_0$ consists of two non-closed $\,Ad_H$-orbits. 
 
\smallskip
\item[$(d)$]
Assume $\,s^+>1$ and  $\,s^->1$.  For every $\,r\not=0\,$ the level set
$\,B_r\,$ consists of a single connected component. It intersects
either  $\,\q \cap \p$  or 
$\,\q \cap \k\,$ in a sphere, depending on whether $r>0$ or $\,r<0$. Thus for every non-zero vector $A\in \q\cap \p$  and  every  $\,\,t>0\,$, the    points 
$\,tA $ and  $\,-tA \,$ belong to the same 
$\,Ad_H$-orbit. A similar statement holds true for  points $\,sC$ and $\,-sC\,$, with $C$ a non-zero vector in $ \q\cap\k$ and $s>0$.
The nilcone $\,B_0$ consists of a unique non-closed $\,Ad_H$-orbit. 
\qed
\end{itemize}

\end{remark}

\bigskip
In order to give further details, we recall the classification of rank-one,
Riemannian symmetric spaces  
of the non-compact type. For each space $M$ we list its real dimension, its standard presentation $\,G/K$,  and the dimensions of the restricted roots spaces of $\,\g$
(cf. \cite{Wo}, p.$\,$294 and \cite{Hl}, p.$\,$532).   

\bn
\vbox{
\noindent\bf Table 4.0
\bn 
\rm
\smallskip

\eightrm{

\centerline{
\offinterlineskip\hbox{\vbox{
\hrule\halign{&\vrule#&\strut\quad#\hfil\quad\cr
height2pt&\omit&&\omit&&\omit&&\omit&&\omit&\cr
&$M$&& $\dim M$&&$G/K$&&$\dim  \g^\alpha$&&$\dim \g^{2\alpha}$&\cr
height2pt&\omit&&\omit&&\omit&&\omit&&\omit&\cr
\noalign{\hrule}
height3pt&\omit&&\omit&&\omit&&\omit&&\omit&\cr
&$H^n (\R)$&&$n$&&$SO_0(n,1)/SO(n) ,~n\ge 2$&&$n-1$&&$0$&\cr
height3pt&\omit&&\omit&&\omit&&\omit&&\omit&\cr
&$ H^n (\C)$&&$2n$&&$SU (n,1)/U(n),~n\ge 2$&&$2(n-1)$&&$1$&\cr
height3pt&\omit&&\omit&&\omit&&\omit&&\omit&\cr
&$ H^n (\HH)$&&$4n$&&$Sp(n,1)/Sp(n)\times Sp(1),~n\ge 2$&&$4(n-1)$&&$3$&\cr
height3pt&\omit&&\omit&&\omit&&\omit&&\omit&\cr
&$ H^2(\C ay)$&&$16$&&$F_4^*/Spin(9)$&&$8$&&$7$&\cr
height3pt&\omit&&\omit&&\omit&&\omit&&\omit&\cr}
\hrule}}}
}
}

\smallskip
 \nbigskip
{\bf Remark.} The two dimensional symmetric space  $\,S_0(2,1)/SO(2)$ can also be identified with $ SU(1,1)/U(1)$ or  $  SL(2,\R)/SO(2)$.  The symmetric space  $\,S_0(3,1)/S O(3) $ can be identified with $\,SL(2,\C)/SU(2) $.

\medskip
\nbigskip 
{\bf  4.1 The reduced case}
 
 \bigskip
Assume that the restricted root system of $\,\g\,$ is {\it reduced}, i.e. it consists of two roots~$\,\{\pm \alpha\}$. 
This is the case of  the spaces $H^n(\R)$ in Table 4.0.
A fundamental domain for the action of
$\,W_K(\a) \ltimes \Gamma \,$ on $\,\a\,$  is given by $\,\a_0=\{\, A\in \a~:~ 0\le \alpha(A)\le \pi \, \}\,$ and there are three singular  orbits intersecting ${\mathcal A_0:= \exp i\a_0K^\C }$.
Their base  points  are given by $\,z_j=g_jK^\C$, for $\,j=1,2,3$.
Here  $\,g_j = \exp iA_j\,$ and
the elements $\,A_j \in \a_0\,$  satisfy the conditions
\begin{equation}
\label{BASEPOINT1}
\alpha(A_1)=0,\qquad  \alpha(A_2)=\pi/2 \qquad \alpha(A_3)=\pi,
\end{equation} 
respectively. 
 The $\,G\,$-orbits through $\,z_1\,$ and $\,z_3\,$ are 
 diffeomorphic to  the symmetric space $\,G/K\,$ and are embedded in 
 $\,G^\C/K^\C\,$ as totally real submanifolds of maximal dimension.
Moreover,  the $\,G$-orbit through   $\,z_2 \,$ is  a 
 rank-one, pseudo-Riemannian  symmetric space $\, G/H\,$ with involution $\,\tau_{z_2}=Ad_{g^{\ }_2}\circ \tau \circ Ad_{g^{-1}_2}$.  The space $\,G/H\,$
 is embedded in $\,G^\C/K^\C\,$ as a closed, 
 totally real submanifold of maximal dimension
 (see Lemma 2.11 and Rem.~2.13 in \cite{Ge1}).  
A  standard Cartan subset starting at $\,z_2\,$ is given by
$\,{\mathcal C}=\,\exp i \c \cdot z_2,$   where  
$\,\c=\R(X+\theta  (X) )\,$ and 
 $\,X  \,$ is  a non-zero vector in $\,\g^\alpha$.    
In the next lemma
we determine the $G$-orbit structure of $G^\C/K^\C$ in a neighbourhood of $z_2$.
Fix   a generator
$\,C\,$ of  $\,\c\,$. 

\bn
 \begin{lem} 
\label{LOCAL1} Assume that the restricted root system of $\g$ is reduced. Let $z_2\in \mathcal A_0$ be the base point of the Cartan subset $\mathcal C$.

\sn
(i) If $\dim G/K>2$, then the orbit $\,G \cdot z_2\,$ is simply connected.
In particular, the  isotropy subgroup  $H$ of $z_2$ in $G$ is connected. 

 \sn
(ii) For every $\,s>0$, the  points 
 $\ \exp( is C) \cdot z_2\ $ and $\  \exp( -is C) \cdot z_2\ $ lie  on
 the same $\,G$-orbit in $\,G^\C/K^\C$ if and only if $\dim\g^\alpha>1$. 

\sn
(iii) If $\dim\g^\alpha>1$, there are two non-closed
$\,G$-orbits in  $\,G^\C/K^\C\,$ containing
$\,G\cdot z_2   \,$ in their closure.
 If $\dim\g^\alpha=1$, such orbits are four.
\end{lem}

\smallskip
\begin{proof}
(i)  Using  the hyperquadric model (cf.$\,$Example \ref{HYPER}), one can verify that 
the orbit of $z_2$ is diffeomorphic to $ SO_0(n,1)/SO_0(n-1,1)$.  In particular, it is  topologically equivalent  to a sphere of dimension $n-1$ and is simply connected for $n>2$. In that case, the isotropy subgroup $H$ is connected, since $G$ is connected by assumption. 
When $n=2$, the orbit $G/H$ is not simply connected. The isotropy subgroup of $z_2$ is either connected (when $G=SO_0(2,1)$) or its quotient  by the ineffectivity subgroup is connected (when $G$ is a non trivial covering of $SO_0(2,1)$).

As a consequence,  (ii) and (iii) of the lemma follow  from
Remark \ref{ADJOINT}, provided  that $\,\dim (\q\cap \p)=1\,$ and
$\,\dim (\q\cap \k)=\dim\g^\alpha$. In order to show this,
 define  $\,\g[\alpha]:=\g^\alpha \oplus  \g^{-\alpha}$. Then $\,\g[\alpha]\,$ is a $\,\theta$-stable subspace of $\,\g\,$ of dimension equal to $\,2\dim \g^\alpha$. Let  
$\,\g[\alpha]= \g[\alpha]_\k\oplus \g[\alpha]_\p\,$ be 
its Cartan decomposition. The components 
$\,\g[\alpha]_\k$ and $\,\g[\alpha]_\p\,$ are  generated  by  vectors of the form
$$X +\theta ( X)  \quad {\rm and } \quad X -\theta (X)$$
respectively, where $\,X \,$ ranges through the elements
of a basis of
$\,\g^\alpha$. In particular
$\dim \g[\alpha]_\k=\dim \g[\alpha]_\p=\dim \g^\alpha.$
Consider the decomposition $\,\g=Z_\k(\a)\oplus \a\oplus \g[\alpha]\,$ 
and note that 
$\,\tau_{z_2}=Ad_{g_2^{\ }}\circ \tau  \circ Ad_{g^{-1}_2}
=Ad_{g_2^2}\circ \theta$. Since $\,Ad_{\exp i A_2}= e^{ad(iA_2)}$,
one has  
$$\,\tau_{z_2} =Id\ \ \ \ {\rm on} \ \ \ \ Z_\k(\a)\,, 
\qquad \tau_{z_2}=-Id  \ \ \ \  {\rm on} \ \ \ \ \a. $$
Since  $\,\alpha(A_2)=\pi/2$, one has 
$\,\tau_{z_2} =-\theta\, $ on $\,\g[\alpha]$.
It follows that
$\,\q:=Fix(-\tau_{z_2},\g)= \a\oplus \g[\alpha]_\k\,$.
In particular,
$\,\dim (\q\cap \p)=\dim \a = 1$ and $\,\dim (\q\cap \k)=\dim \g[\alpha]_\k=\dim\g^\alpha$,
as wished.
\end{proof}

 \mn
 
From the above discussion  and  Table 4.0  it follows that  in the reduced case  
the $G$-orbit space of $G^\C/K^\C$ can be  described by the following  diagrams.



\nbigskip

\noindent
{ $\blacklozenge$ $\,\ G/K=SO_0(2,1)/SO(2).$}

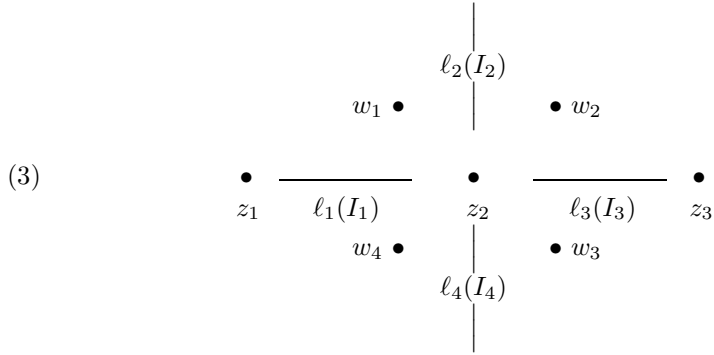
\nsmallskip
\begin{equation}
\label{DIAGRAM1}
\begin{matrix}  
&  &\Big|& &   \cr 
&\quad \quad  \quad \quad \quad \quad \quad \quad &   \ell_2( I_2) & &  \cr
&\ \ \ \  \quad \quad \quad \quad w_1 \ \bullet & \Big|&\bullet \  w_2 \quad \quad \quad \quad  \quad \ \ &  \cr 
&  && &   \cr 
 &\ \ \ \   \bullet \quad \overline{ \quad \quad \quad \quad \quad  }&\bullet&\overline{ \quad \quad \quad \quad \quad   }\quad \bullet  \ \ \ \ \ \ \ &\cr  
  &\  z_1 \quad \ \ \ \ell_1(I_1)  \      &\  z_2 & \ell_3(I_3)\ \ \   \quad  z_3\ \ \cr
  & \ \ \ \  \quad \quad \quad  \quad w_4\ \bullet&\Big|&\bullet\ w_3\ \  \quad \quad \quad \quad  \quad&  \cr 
  & &  \ell_4(I_4)& &  \cr 
  &  &\Big|& &   \cr 
  \end{matrix}
  \end{equation}


\nbigskip 
 $\blacklozenge$    $\,\ G/K=SO_0(n,1)/SO(n)$, $\,n>2$.

\bigskip
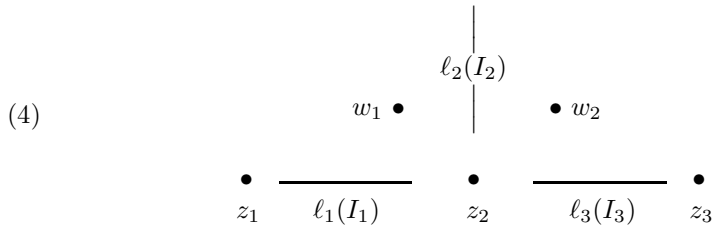
\begin{equation}
\label{DIAGRAM2}
\begin{matrix}  
&  &\Big|& &   \cr 
&\quad \quad  \quad \quad \quad \quad \quad \quad &   \ell_2( I_2) & &  \cr
&\ \ \ \  \quad \quad \quad \quad w_1 \ \bullet & \Big|&\bullet \  w_2 \quad \quad \quad \quad  \quad \ \ &  \cr 
&  && &   \cr 
  &\ \ \ \   \bullet \quad \overline{ \quad \quad \quad \quad \quad  }&\bullet&\overline{ \quad \quad \quad \quad \quad   }\quad \bullet  \ \ \ \ \ \ \ &\cr  
  &\  z_1 \quad \ \ \ \ell_1(I_1)  \      &\  z_2 & \ell_3(I_3)\ \ \   \quad  z_3\ \ \cr
 \end{matrix}
 \end{equation}

\bn 
\sn
Set $\,I_1=I_3=(0,1)\,$. For $\,j=1,3$, the maps $\,\ell_j: I_j \to G^\C/K^\C$, defined by 
\begin{equation}
\label{SLICE13}
\ell_1(t):= \exp(-itA_2) \cdot z_2 , \quad \ell_3(t):= \exp (itA_2) \cdot z_2,
\end{equation}

\sn 
parametrize the principal $G$-orbits through $\mathcal A_0$. One has
$$\,\mathcal A_0\ =\ z_1 \cup \ell_1(I_1)\cup z_2\cup \ell_3(I_3)\cup z_3.\,$$ 

\sn
Set
$\,I_2=I_4=(0,\infty)\,$. For $j=2,4$,  the maps $\,\ell_j: I_j \to G^\C/K^\C$,  defined by
\begin{equation}
\label{SLICE24}
\ \,  \ell_2(s):= \exp (isC) \cdot z_2 , \quad \quad \ \ell_4(s):= \exp (-isC) \cdot z_2,\, \end{equation}

\sn
parametrize the principal closed $G$-orbits through the standard Cartan subset $\,\mathcal C\,$ and

$$\mathcal C= \ell_2(I_2)\cup  z_2 \cup  \ell_4(I_4) .$$
The points 
$$w_1,w_2, w_3,w_4$$

\sn
represent the non-closed $G$-orbits containing the
singular orbit $G\cdot
z_2$ in their closure.

 \bn 
\begin{exa} 
\label{HYPER} {\it The complex hyperquadric.}
{\rm Let  $G=SO_0(n,1) $, with
$n\ge 2$, and  let $\,G^\C=SO(n,1,\C)\,$ be its universal  complexification. 
By definition $G^\C$ is the subgroup of $\,SL(n+1, \C)\,$ leaving invariant the
quadratic form of signature $(n,1)$. The space $G^\C/K^\C$  can be
identified with the
$G^\C$-orbit through $ (0, \dots , 1)$ which coincides
with the $n$-dimensional
complex hyperquadric  
$$M^\C = \{ (\xi_1, \dots , \xi_{n+1})\in \C^{n+1} 
\ : \   \xi_1^2+\ldots +\xi_n^2 - \xi_{n+1}^2=-1 \}.$$ 
Fix the elements 
$$A_2=\left ( \begin{matrix}  0  &  \dots & 0& 0  \cr
                                            \vdots  &     & \vdots& \vdots  \cr
                                           0  &  \dots &  0 & \frac{\pi}{2} \cr 
                                           0  &  \dots & \frac{\pi}{2} & 0 \cr 
\end{matrix} \right ) \quad {\rm and} \quad 
C =\left( \begin{matrix} 
                              0 & \dots  & \ 0  &  0 & 0 \cr 
                              \vdots &   & \vdots   &  \vdots & \vdots  \cr
                              0 & \dots  & 0  &  -2\ & 0  \cr
                              0 & \dots  & 2   &  \ 0 & 0 \cr 
                              0 & \dots  & \ 0  &  \ 0 & 0 \cr 
\end{matrix} \right )  \ $$

\nmedskip
 in $\g$ as  generators of $\a$ and $\c$,  respectively. Then points on 
 the singular orbits in $M^\C$ satisfying conditions
(\ref{BASEPOINT1})  are given by
$$ z_1=(0,\dots,0,1)\,,\quad  \quad z_2=(0,\dots,0,i,0)\,,\quad  \quad z_3=(0,\dots,0,-1).\, $$
The $G$-orbit of $z_2$ is diffeomorphic to
the pseudo-Riemannian symmetric
space $\,G/H\cong SO_0(n,1)/SO_0(n-1,1)$.   
The   slices  $\ell_1$ and $\ell_3$ are given by
$$\ell_1(t)=(\,0, \dots,0, i\sin \frac{\pi}{2}(1 -t),\cos \frac{\pi}{2}(1 -t)\,),\quad t\in (0,1),$$
$$\ell_3(t)= (\,0, \dots,0 , i\sin \frac{\pi}{2}(1 +t),\cos \frac{\pi}{2}(1 +t)\,),  \quad t\in (0,1). $$
The slices $\ell_2$ and $\ell_4$ are given by
$$ \ell_2(s)=(\,0, \dots , 0, \sinh 2s, i\cosh 2s, 0\,),\quad s>0, $$
$$ \ell_4(s)=(\,0, \dots , 0, -\sinh 2s, i\cosh 2s, 0\,),\quad s>0 .$$
The slice representation at $z_2$  is equivalent to  the linear action of $SO_0(n-1, 1)$  on 
$\R^{n}$.  
When $n=2$,  
we can choose representatives 
of the four non-closed
hypersurface $G$-orbits containing 
$G \cdot z_2$ in their closure to be
$$\,w_1=(-1,i ,-1),\ \    w_2=(1, i, -1),\ \  w_3=(1,i,1) \ \ w_4=(-1, i,  1)\,. $$
When $n>2$, the slice representation identifies
$ \ell_2 $ and $ \ell_4 $ and representatives 
of  the two non-closed hypersurface $G$-orbits containing 
$G \cdot z_2$ in their closure are for example
$$\,w_1=(-1,0,\dots,0, i,-1)\quad {\rm and} \quad   w_2=(1,0,\dots,0, i, -1).$$ \qed }
\end{exa}


\nbigskip
{\bf 4.2 The non-reduced case.}

 \bigskip
Assume that the restricted root system of $\,\g\,$  is non-reduced, i.e.  it
consists of four roots  $\, \{\pm \alpha,~\pm 2\alpha \}.$ 
This is the case of $H^n(\C)$, $H^n(\HH)$ and  $H^2(\C ay)$ in Table 4.0.
A fundamental domain 
for the action of $\,W_K(\a) \ltimes \Gamma \,$ in $\,\a\,$
is given by $\,\a_0=\{ \,A\in \a~:~ 0\le \alpha(A)\le \pi/2\,\},\ $ and there  are three singular  orbits intersecting $\,{\mathcal A_0}$.
 Their base  points  are given by $\,z_j=g_jK^\C$, for $\,j=1,2,3$.
Here  $\,g_j = \exp iA_j\,$ and
the elements $\,A_j \in \a_0\,$  satisfy the conditions
 \begin{equation}
 \label{BASEPOINT2}
  \alpha(A_1)=0\,, \qquad \alpha(A_2)=\pi/4 \,, \qquad \alpha(A_3)=\pi/2\,,
  \end{equation}
respectively.
The $\,G\,$-orbits through $\,z_1\,$ is diffeomorphic to  the symmetric space $\,G/K$,
the one through $\,z_3\,$ is diffeomorphic to a  rank-one, pseudo-Riemannian 
symmetric space $\, G/H$. Both orbits are embedded in 
 $\,G^\C/K^\C\,$ as totally real submanifolds of maximal dimension.
 (see Lemma 2.11 and Rem.~2.13 in \cite{Ge1}). 
The orbit 
of $\,z_2\,$ is a homogeneous space $\,G/H'$, with $\,H':=G_{z_2}$, and 
$\,\dim G/H' > \dim G/K$ (see Lemma 2.14 and  Rem.~2.15 in \cite{Ge1}). 
Set $\,G':=Z_G(g_2^4)$, where $Z_G(g_2^4)$
denotes the centralizer of $g_2^4$ in $G$.
 Then  $\,H'\,$ is contained in $\,G'\,$ and $\,G'/H'\,$ is a
 rank-one, pseudo-Riemannian symmetric space with involution 
 $\,\tau_{z_2}=Ad_{g^{\ }_2} \circ \tau \circ Ad_{g_2^{-1}}$. 
 Moreover,  the slice representation at $z_2$ is equivalent to the isotropy representation of $\,G'/H'\,$ (cf. \cite{Ge2}).
The  standard Cartan subset starting at $\,z_2\,$ is  given by
$\,{\mathcal C'}=\,\exp i \c' \cdot z_2$, where  
$\,\c'=\R(X +\theta (X ))\,$ and 
 $\,X  \,$ is a non-zero vector  in $\,\g^{2\alpha}$. 
If  $Z_\k(\a)\oplus \a\oplus \g^{\pm  \alpha}\oplus \g^{\pm 2\alpha} $ is the restricted root decomposition of $\g$, then the Lie algebra of $G'$ is given by 
\begin{equation}
\label{GIPRIMO}
 \g' =Z_\k(\a)\oplus \a\oplus \g^{\pm 2\alpha} . 
\end{equation}
Moreover, if $ \h'\oplus \q'$ is the 
 $\tau_{z_2}$-decomposition
  of $\g'$, then $\c'$ is a
maximal abelian subalgebra in $\q'$. Fix a generator $C'$ of  $ \c'$.

  \bn
\begin{lem} 
\label{LOCAL2}
\label{sing1} 
Assume that the restricted root system of $\g$ is non-reduced. 
Let $z_2\in \mathcal A_0$ be the base point of the Cartan subset $\mathcal C'$.

\sn
\item{(i)} The isotropy subgroup $H'$ of $z_2$ in $G$ is connected.

\sn
\item{(ii)} For every $t>0$, the points $ \exp( it C')\cdot  z_2$ and $  \exp( -it C')\cdot  z_2$ sit  on the same $G$-orbit if and only if  $\dim \g^{2\alpha}>1$.

\sn
\item{(iii)}  If  $\dim \g^{2\alpha}>1$, there  are  two non-closed $G$-orbits in $G^\C/K^\C$ 
containing $G\cdot z_2$ in their closure.  If  $\dim \g^{2\alpha}=1$, such orbits are four. 
 
\end{lem}

\mn\begin{proof}
 (i) The group $H'$ is connected if and only if $H'\cap K$ is connected.  Note that  $G'=Z_G(g_2^4)$ is $\theta $-stable, since so is $G$ and $\theta(g_2^4)=g_2^{-4}$.
Therefore   $H'\cap K$ is the common fixed point subgroup of the two involutions 
$\tau_{z_2 }$ and $\theta$ of  $G'$.  As a result,  $H'\cap K =Z_K(g^2_2)$. 
Now regard  $z_2$ as a point
on the compact dual symmetric space $U/K$  endowed with the $K$-action by left translations. 
Denote by $K_{z_2}$ the isotropy subgroup of $z_2$ in $K$. On the one
hand, $K_{z_2}=Z_K(g_2^2)$. On the other hand, since the isotropy subalgebra
$\k_{z_2}$ is given by $\k\cap Ad_{z_2}(\k)$, one sees that $\k_{z_2}$ has minimal dimension and
coincides with $Z_\k(\a)$ if and only if $\alpha(A_2)\not= m\pi$, for $m\in\Z$.  By
(\ref{BASEPOINT2}), it follows that   
 $K_{z_2}$  is 
 principal and consequently is equal to $Z_K(\a)$. 
Finally $Z_K(\a)$ is connected for all rank-one, Riemannian
symmetric spaces of dimension greater than two  (see \cite{Kn}  or Lemma
\ref{TOPORB1} for a direct proof). In conclusion
$$H' \cap K= Z_K(g_2^2)=K_{z_2}=Z_K(\a)$$
implying $(i)$.
 
Parts $(ii)$ and $(iii)$ of the lemma follow  by applying
Remark \ref{ADJOINT} to the symmetric space $G'/H'$, provided
that 
$$\dim \q'\cap \p=1 ,\qquad \dim \q'\cap \k=\dim\g^{2\alpha}.$$
 In order to show this, 
 define  $\,\g[2\alpha]:=\g^{2\alpha} \oplus  \g^{-2\alpha}$. Then $\,\g[2\alpha]\,$ is $\,\theta$-stable subspace of $\,\g\,$ of dimension equal to $\,2\dim \g^{2\alpha}$. Let  
$\,\g[2\alpha]= \g[2\alpha]_\k\oplus \g[2\alpha]_\p\,$ be 
its Cartan decomposition. The components 
$\,\g[2\alpha]_\k$ and $\,\g[2\alpha]_\p\,$ are  generated  by  vectors of the form
$$X +\theta (X)  \quad {\rm and } \quad X -\theta (X) $$

\nmedskip
respectively, where $\,X \,$ ranges through the elements of a basis of
$\,\g^{2\alpha}$. In particular
$\dim \g[2\alpha]_\k=\dim \g[2\alpha]_\p=\dim \g^{2\alpha}.$
One sees that 
$$\,\tau_{z_2} =Id\ \ \  {\rm on}  \  \ \ Z_\k(\a)\,, 
\qquad \tau_{z_2}=-Id   \ \  \ {\rm on} \    \ \ \a,  \qquad \tau_{z_2} = -\theta 
\  \ \ {\rm on} \  \ \ \g[2\alpha]. $$
Consequently
$\,\q':=Fix(-\tau_{z_2},\g')= \a\oplus \g[2\alpha]_\k\,$
and
$\,\dim (\q'\cap \p)=\dim \a = 1\,$.
Similarly, $\,\dim (\q'\cap \k)=\dim \g[2\alpha]_\k=\dim\g^{2\alpha}$,
as wished.

\end{proof}

\bn 

By Lemma 2.11 and  Remark$\,$2.13 in \cite{Ge1},    the $\,G$-orbit of   $\,z_3 \,$ is  a rank-one, 
 pseudo-Riemannian  symmetric space $\, G/H\,$ with involution $\,\tau_{z_3}=Ad_{g^{\ }_3}\circ \tau \circ Ad_{g^{-1}_3}$.  The space $\,G/H\,$
 is embedded in $\,G^\C/K^\C\,$ as a closed, 
 totally real submanifold of maximal dimension.  
The  standard Cartan subset starting at $\,z_3\,$ is given by
$\,{\mathcal C}=\,\exp i \c \cdot z_3 $,  where  
$\,\c=\R(X+\theta  (X) )\,$ and 
 $\,X  \,$ is  a non-zero vector in $\,\g^\alpha$.    
If $\g=\h\oplus \q$ is 
 the $\tau_{z_3}$-decomposition of $\g$, then $\c$ is a maximal abelian subalgebra in $\q$. Fix a generator $C$ of $\c$.

\bn
\begin{lem} 
\label{LOCAL3} 
Assume that the restricted root system of $\g$ is non-reduced.
Let $z_3\in \mathcal A_0$ be the base point of the Cartan subset $\mathcal C$.

\sn
\item{(i)} The orbit $\,G \cdot z_3\,$ is simply connected. In
particular the isotropy subgroup $H$ of $z_3$ in $G$  is connected.

\sn
\item{(ii)} For every $t>0$, the   points  $ \exp( it C)\cdot  z_3$ and $ \exp( -it C)\cdot  z_3$ sit  on the same $G$-orbit in $G^\C/K^\C$. 

\sn
\item{(iii)}  There  is precisely one non-closed $G$-orbit in $G^\C/K^\C$ containing
$G\cdot z_3$ in its closure. 
\end{lem}

\mn 
\begin{proof}
(i) Since by assumption $G$ is connected, we prove that  $\,H\,$ is connected by showing that 
the orbit $\,G\cdot z_3\,$ is simply connected. 
In order to do this, by Remark~\ref{QUOTOP},
it is sufficient to choose $\,G\,$ as in the standard presentation in Table~4.0.
Let  $G =SU(n,1)$.  By direct computations (cf.~Example~\ref{NONREDUCED1})  one finds that $G\cdot z_3\cong SU(n,1)/ U(n-1, 1) $.  This quotient  is topologically equivalent to the complex projective space $\C\P^{n-1}$. In particular,  it  is simply connected. \pn
Consider then  $G=Sp(n,1) $ or $G =F_4^* $. In both cases  the group $G$ is simply connected. Since $H$ is the fixed point subgroup of an involution of $G$, it is connected (cf. \cite{Sb}).
It follows that the quotient  is simply connected.

Parts (ii) and (iii)  follow  from Remark \ref{ADJOINT}, provided that
$\,\dim (\q\cap \p)=1+ \g^{2\alpha}\,$ and
$\,\dim (\q\cap \k)=\dim\g^\alpha$.
 In order to show this, define  $\,\g[\alpha]:=\g^\alpha \oplus  \g^{-\alpha}\,$ and
$\,\g[2\alpha]:=\g^{2\alpha} \oplus  \g^{-2\alpha} $.
Then both $\,\g[\alpha]$ and $\,\g[2\alpha]\,$  are
$\,\theta$-stable subspaces of $\,\g\,$ of dimension equal to
$\,\dim \g^\alpha$ and $\,2\dim \g^{2\alpha}\,$ respectively.
Let $\,\g[\alpha]_\k$, $\,\g[\alpha]_\p$, $\,\g[2\alpha]_\k\,$ and $\,\g[2\alpha]_\p\,$
be the  components of the respective Cartan decompositions. The  same arguments
as in  the proof of Lemma \ref{LOCAL1} and Lemma \ref{LOCAL2} show that 
$$\dim \g[\alpha]_\k=\dim \g[\alpha]_\p=\dim \g^\alpha\ \ \ \ \  \ \ \ 
\dim \g[2\alpha]_\k=\dim \g[2\alpha]_\p=\dim \g^{2\alpha}. $$

\nsmallskip
Moreover, one sees that 
$$\,\tau_{z_3} =Id\ \, \ {\rm on} \ \, \ Z_\k(\a), 
 \,  \quad \tau_{z_3}=-Id  \ \, \  {\rm on} \ \, \ \a,\quad \,\tau_{z_3} =-\theta\ \, \ {\rm on} \,  \ \ \g[\alpha], 
\,  \quad \tau_{z_3} = \theta\, \ \ {\rm on} \ \, \ \g[2\alpha]. $$
Since 
$$\g=Z_\k(\a)\oplus \a \oplus \g[\alpha]\oplus \g[2\alpha], $$ 
it follows that 
$\,\q:=Fix(-\tau_{z_3},\g)= \a\oplus \g[\alpha]_\k\oplus \g[2\alpha]_\p  $.
In particular, 
$\,\dim (\q\cap \p)=1 + \dim \g^{2\alpha} $ and 
 $\,\dim (\q \cap \k)=\dim \g^\alpha $,
as claimed.
\end{proof}

\bn
As a consequence  of the above lemmas and Table 4.0, in the non-reduced case  the
$\,G$-orbit space of $\,G^\C/K^\C\,$ can be represented by the following diagrams. 

\nbigskip
$\blacklozenge$  $\,\ G/K=SU(n,1)/ U(n)$, $\,n \ge 2$.

\medskip

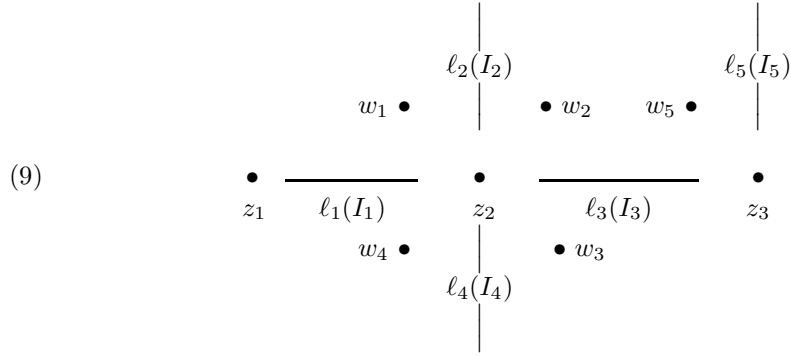
\begin{equation}
\label{DIAGRAM3}
\begin{matrix}  
&  &\Big|& & \Big|  \cr 
&\quad \quad  \quad \quad \quad \quad \quad \quad &   \ell_2( I_2) & &  \ell_5( I_5) \cr
&\ \ \ \  \quad \quad \quad \quad w_1 \ \bullet & \Big|&\bullet \  w_2  \quad \quad  w_5 \ \bullet & \Big|  \cr 
&  && &   \cr 
  &\ \ \ \   \bullet \quad \overline{ \quad \quad \quad \quad \quad  }&\bullet&\overline{ \quad \quad \quad \quad \quad \quad  }  & \bullet \cr  
  &\  z_1 \quad \ \ \ \ell_1(I_1)  \      &\  z_2 & \ell_3(I_3)  & z_3\cr
  & \ \ \ \  \quad \quad \quad  \quad w_4\ \bullet&\Big|&\bullet\ w_3\ \  \quad \quad    \   &   \cr 
  & &  \ell_4(I_4)& &  \cr 
  &  & \Big|& &     \cr 
  \end{matrix}
    \end{equation}
 
 \bn
 \nbigskip

\noindent
$\blacklozenge$  $\ \  G/K=Sp(n,1)/Sp(n)\times Sp(1)$, $\ n \ge 2$

\sn
$\ $  $\ \ G/K=F_4^*/Spin(9)$

\medskip
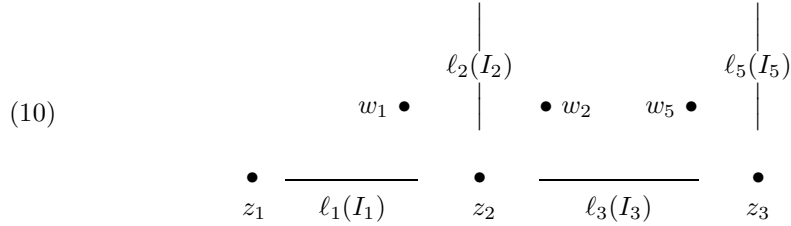
\begin{equation}
\label{DIAGRAM4}
\begin{matrix}  
&  &\Big|& & \Big|  \cr 
&\quad \quad  \quad \quad \quad \quad \quad \quad &   \ell_2( I_2) & &  \ell_5( I_5) \cr
&\ \ \ \  \quad \quad \quad \quad w_1 \ \bullet & \Big|&\bullet \  w_2  \quad \quad  w_5 \ \bullet & \Big|  \cr 
&  && &   \cr 
  &\ \ \ \   \bullet \quad \overline{ \quad \quad \quad \quad \quad  }&\bullet&\overline{ \quad \quad \quad \quad \quad \quad  }  & \bullet \cr  
  &\  z_1 \quad \ \ \ \ell_1(I_1)  \      &\  z_2 & \ell_3(I_3)  & z_3\cr
     \end{matrix} 
     \end{equation}

 \bn
Set $\,I_1=I_3=(0,1)$.  For $j=1,3$, define 
$\,\ell_j\colon I_j\longrightarrow G^\C/K^\C\,$ by
\begin{equation}
\label{SLICENR13}
\ell_1(t)=\exp(-i tA_2)\cdot z_2\,, \qquad  \ell_3(t)=\exp (itA_2)\cdot z_2\,.
\end{equation}

\noindent
The slices $\ell_1$ and $\ell_3$ parametrize the principal $G$-orbits through $\mathcal A_0$ and  
$$\,\mathcal A_0 \,=\,z_1 \cup \ell_1(I_1)\cup z_2\cup \ell_3(I_3)\cup z_3\,. $$
Set $I_2= I_4=(0,\infty)\,$.  For $j= 2,4 $, define 
$\,\ell_j\colon I_j\longrightarrow G^\C/K^\C\,$ by 
\begin{equation}
\label{SLICENR24}
\ell_2(s)=\exp( sC' )\cdot z_2\,, \qquad  \ell_4(s)=\exp (-sC')\cdot z_2\,. 
\end{equation}

\noindent
The slices $\ell_2$ and $\ell_4$ parametrize the principal $G$-orbits through the Cartan subset $\mathcal C'$ with base point $z_2$ and
$$\mathcal C'=  \ell_2(I_2)\cup z_2 \cup \ell_4(I_4).$$
Finally, set $I_5=(0,\infty)$ and define $\ell_5\colon I_5\to G^\C/K^\C$ by
\begin{equation}
\label{SLICENR5}
\ell_5(s)=\exp( s C)\cdot z_3\, .
\end{equation}

\noindent
The slice $\ell_5$ parametrizes the principal $G$-orbits through
the standard Cartan subset $\,\mathcal C\,$ with base point $\,z_3$. 
The points $\,w_1,\ldots,w_4\,$
represent  the   non-closed orbits containing $\,G\cdot z_2\,$ in their closure.
The point  $\,w_5\,$ represents  the non-closed orbit containing
$\,G\cdot z_3\,$ in its closure.

\bn
\begin{exa}
\label{NONREDUCED1}{\it A model in the non-reduced case}. {\rm 
Let $G=SU(n,1)$, with $\,n \ge 2$, be the subgroup of $SL(n+1,\C)$ leaving invariant
the hermitian form
$\langle z,w\rangle_{n,1}= z_1\bar w_1 +\ldots+z_n\bar w_n - z_{n+1} \bar w_{n+1}$ in
$\C^{n+1}$. Denote by $\sigma$ the conjugation of $G^\C=SL(n+1,\C)$ relative to $G$,
namely 
$\sigma(g)=I_{n,1}{}^t\bar g^{-1}I_{n,1}$. Denote by $\overline{\P}^n $
the complex
projective space endowed with the opposite complex structure,
i.e. the one for which the
map $\,\P^n\to \overline\P^n, ~[z]\mapsto [\bar z]\,$ is holomorphic. The group $G^\C$
acts holomorphically on $\P^n \times 
\overline{\P}^n $ by 
$$g\cdot([z],[ w]):=([g\cdot z],[\sigma(g)\cdot  w] ).$$

\nsmallskip
 Under this action $\,\P^n \times \overline {\P}^n\,$ consists of two orbits:
a closed
one given by 
  $\,\{\,([z],[ w]) \in \P^n \times \overline {\P}^n \ :\ 
  \langle z,  w\rangle _{n,1}=0\,\}\,$ and an open one given by its complement.
 The quotient $\,G^\C/K^\C\,$ can be  identified with the open orbit
$$ M^\C:= G^\C \cdot ([0: \dots :0:1],[0: \dots:0 :1])
=\P^n \times \overline {\P}^n \setminus \{\,\langle z,  w \rangle _{n,1} = 0\,\}.$$

\nsmallskip
Fix the elements
$$A_2=\left ( \begin{matrix}  0  &  \dots & 0& 0  \cr
                                            \vdots  &     & \vdots& \vdots  \cr
                                           0  &  \dots &  0 & \frac{\pi}{4} \cr 
                                           0  &  \dots & \frac{\pi}{4} & 0 \cr 
\end{matrix} \right ),~~  C'=\left ( \begin{matrix}  0  &  \dots & 0& 0  \cr
                                            \vdots  &     & \vdots& \vdots  \cr
                                           0  &  \dots &   i & 0 \cr 
                                           0  &  \dots & 0& -i \cr 
\end{matrix} \right ), ~~
C =\left( \begin{matrix} 
                              0 & \dots  & \ 0  &  0 & 0 \cr 
                              \vdots &   & \vdots   &  \vdots & \vdots  \cr
                              0 & \dots  & 0  &  -1\ & 0  \cr
                              0 & \dots  & 1   &  \ 0 & 0 \cr 
                              0 & \dots  & \ 0  &  \ 0 & 0 \cr 
\end{matrix} \right )  \ $$
in $\g$ as generators of $\a$, $\c'$ and $\c$, respectively.
Then  points on the singular orbits in $M^\C$ satisfying conditions
(\ref{BASEPOINT2}),
are given by
$$ z_1=([0:\ldots:0:1], [0:\ldots:0:1])\,,  \qquad z_2=([0:\ldots:0:i :1], [0:\ldots:0:-i :1])\,$$
$$ z_3=([0:\ldots:0:1:0], [0:\ldots:0:1:0]) .\, $$
 The $G$-orbit of $z_2$ is diffeomorphic to the homogeneous space $G/H'$,
where
$H'\cong  U(n-1)\times  SO(1,1)$. The group
$G'$ is isomorphic to $ U(n-1)\times SU(1,1) $ and  the quotient $G'/H'$ is
diffeomorphic to the two-dimensional rank-one, pseudo-Riemannian symmetric space
$SU(1,1)/SO(1,1)$.  The $G$-orbit of $z_3$ is diffeomorphic to the pseudo-Riemannian
symmetric space $SU(n,1)/ SU(n-1,1) $.  The   slices  $\ell_1$ and $\ell_3$ are given
by
$$\ell_1(t)=([0:\ldots:i\sin {\pi\over 4}(1-t):\cos {\pi\over 4}(1-t)], [0:\ldots:-i\sin {\pi\over 4}(1-t):\cos {\pi\over 4}(1-t)]) ,$$
$$\ell_3(t)= ([0:\ldots:i\sin {\pi\over 4}(1+t):\cos {\pi\over 4}(1+t)], [0:\ldots:-i\sin {\pi\over 4}(1+t):\cos {\pi\over 4}(1+t)]),$$
where $t\in(0,1)$.
The slices $\ell_2$ and $\ell_4$ are given by
$$ \ell_2(s)=([0:\ldots:ie^{-s} :e^{ s} ], [0:\ldots:-ie^{ s} :e^{ -s} ])\,, $$
$$ \ell_4(s)=([0:\ldots:ie^{s} :e^{ -s} ], [0:\ldots:-ie^{ -s} :e^{s} ])\,,$$
with $\,s>0$. Finally the slice $\ell_5$ is given by
$$ \ell_5(s)=([0:\ldots:\sinh s:i\cosh s:0], [0:\ldots:\sinh s:-i\cosh s:0])\, ,$$
with $\,s>0$.
The slice representation at $z_2$ is equivalent to the standard action of $SO(1,1)$ on $\R^2$. So there are four non-closed  $G$-orbits containing  $G\cdot z_2$ in their closure. We can choose representatives of 
 such orbits to be
$$\,w_1=([0:\ldots:0:1],[0:\ldots:-i:1])\,,\ \  w_2=([0:\ldots:i:1], [0:\ldots:1:0])\,
,$$
$$\ \  w_3=([0:\ldots:1:0],
[0:\ldots:-i:1])\,, \ \ w_4=([0:\ldots:i:1], [0:\ldots:0:1])\,.$$
 A representative for the unique  non-closed orbit containing  $G\cdot z_3$
 in its closure is given  by
$$w_5=([0:\ldots:1: -i : 1], [0:\ldots:1:i:1]).$$
\qed }
\end{exa}


\bn
\begin{remark}
\label{N=1}
 When $\,G=SU(1,1)$, the restricted root system 
of $\,\g\,$ is reduced. The quotient $\,G^\C/K^\C\,$ can be
identified with $\,\P^1 \times \overline{\P}^1 \setminus \{\langle z, w
\rangle_{1,1}=0 \}\,$ and the $\,G$-orbit space can be described as above,
except for the fact that the slice $\,\ell_5\,$ and the point $\,w_5\,$
must to be omitted. Moreover the $\,G$-orbit through $\,z_3\,$ is diffeomorphic
to the symmetric space $\,G/K$. Note that $\,SU(1,1)^\C/U(1)^\C \,$ is biholomorphic
to $\,SO_0(2,1)^\C / SO(2)^\C$. Thus it can also be identified with the
two-dimensional hyperquadric described in Example \ref{HYPER}.
\qed
\end{remark}

\bn

\section{Univalence on $G$-orbits in  $G^\C/K^\C$}

\bn

Let $\,G \,$ be a connected, non-compact,  real simple Lie group, $\,K \subset G\,$ a maximal compact subgroup  and $\,G^\C\,$ the universal complexification of $\,G$.  Assume that $\,G\,$ is embedded in $\,G^\C$.  
Consider a  $\,G$-equivariant Riemann domain 
$$q\colon \Sigma \to G^\C/K^\C.$$
The main goal of this section is to prove that
$\,q\,$ is injective on $\,G$-orbits, 
if $\,G/K\,$ is a  rank-one, Riemannian
symmetric space  of  dimension  greater than three.  We  first prove the result  for
principal $\,G$-orbits and  later we extend it to all $\,G$-orbits by a general
argument.  In most cases, the injectivity of $\,q\,$ on principal $\,G$-orbits  follows
from their simple connectedness.   The cases $\,\dim G/K =2,3\,$
are discussed separately.

\bn

Recall that by $(ii)$ of Remark \ref{QUOTOP}, different triples $\, (G, K, G^\C)\,$  associated with the same Riemannian symmetric space $G/K$ 
  yield   $G$-equivariantly diffeomorphic orbits in $G^\C/K^\C$.  
Let $\mathcal A_0$, $\mathcal C' $ and $\mathcal C $ be the standard Cartan subsets in $G^\C/K^\C$. Let  
$H$ be the isotropy subgroup  of the base point  of $\mathcal C $ and $H'$   the isotropy subgroup  of the base point  of $\mathcal C' $ (see  Lemmas \ref{LOCAL1}, \ref{LOCAL2} and  \ref{LOCAL3}). 
By Prop.$\,$3.4 and  Prop.$\,$3.15 in \cite{Ge1}, the principal orbits intersecting 
$\mathcal A_0$,  $\mathcal C$ and
$\mathcal C'$ have isotropy type 
$$Z_K(\a),\quad Z_H(\c)\quad \hbox{and}\quad  Z_{H'}(\c'),$$
respectively.

 \bn
\begin{lem}
\label{TOPORB1} Principal $G$-orbits of isotropy type $Z_K(\a)$ are simply connected
if and only if $\dim G/K>2$.
\end{lem} 

\mn
\begin{proof} 
An orbit $\,G/Z_K(\a)\,$ is topologically equivalent to $\,K/Z_K(\a)$. Consider the isotropy representation of $\,K\,$ on $\,\p$. The non-zero
$\,K$-orbits  in $\,\p$ are diffeomorphic to $\, K/Z_K(\a)\,$. Since $\,G/K\,$ has rank one,  they are also diffeomorphic to spheres  of dimension $\,\dim (G/K) -1$. Hence the statement follows.
\end{proof}


 \bn
 \begin{remark} \label{SO21A}  When $G=SO_0(2,1)$,  the isotropy subgroup  $\,Z_K(\a)$ is trivial. 
Therefore principal orbits of type $\,G/Z_K(\a)\,$ are diffeomorphic to
$\,SO_0(2,1)$ and  topologically equivalent to $\,SO(2)$.
In particular, they are not simply connected. 
\qed
\end{remark}

\bn
\begin{lem}
\label{TOPORB2} Principal $G$-orbits of isotropy type $Z_H(\c )$ are simply
connected,  except  when $\,G\,$ is one of the groups
$\, SO_0(2,1)$, $\,SO_0(3,1)\,$ or $\, SU(2,1)$. 
\end{lem} 

\sn
\begin{proof}
An orbit $\,G/Z_H(\c)\,$ is topologically equivalent to $\,K/Z_{K\cap H}(\c)$. 
We prove the lemma by discussing each case separately. 
Let  $G =SO_0(n,1)$.
Using the hyperquadric model given in  Example \ref{HYPER}, one checks that
 $$H\cong SO_0( n-1,1),\quad   Z_H(\c)\cong  SO_0( n-2,1),\quad  K/Z_{H\cap K}(\c)\cong SO(n)/SO(n-2).$$
In particular, $\,K/Z_{K\cap H}(\c)$ is diffeomorphic to a Stiefel manifold, which is simply connected   for $\,n>3$.

Consider next the  case  $\,G =SU(n,1) $, with $\,n \ge 3$.
Direct computations on the
model in Example \ref{NONREDUCED1} show that
\sn
$$H\cong U( n-1,1),\quad \quad  
Z_{K\cap H}(\c)\cong U(n-2)\times U(1),$$
$$ \quad K/Z_{K\cap H}(\c)\cong U(n)/(U(n-2)\times U(1)) .$$

\sn
Since, for $\,n \ge 3$, the embedding  $\,U(n-2) \to U(n)\,$ induces an epimorphism
 of fundamental groups, so does  the embedding $\,U(n-2) \times U(1) \to U(n)$.
  As a consequence, 
 $\,K/Z_{K\cap H}(\c)$  is simply connected.

 Finally, consider   $\,G =Sp(n,1) \,$ or $\,G =F_4^* \,$. 
Note that in both cases $\,K\,$  is simply connected. Therefore  $\,K/Z_{K\cap H}(\c)\,$
is simply connected provided that $\,Z_{K\cap H}(\c)\,$ is connected.
In order to show that this,
consider  the compact, rank-one, symmetric space $\,K/K\cap H\,$ and the  corresponding 
 isotropy representation of $\,K\cap H\,$ on $\,\k\cap \q$.
The non-zero $\,K\cap H$-orbits in  $\,\k\cap \q\,$ are of type   
$\,K\cap H/Z_{K\cap H}(\c)\,$ and are diffeomorphic to 
spheres  of dimension  $\, \dim (\k\cap \q) -1.$ 
Since $\,\dim (\k\cap \q)=\dim\g^\alpha  >2\,$  
 (cf. Table$\,$4.0), they are simply connected.
By  Lemma \ref{LOCAL1} or Lemma \ref{LOCAL3}, the group  $\,H\,$ and likewise its
maximal compact subgroup $\,K\cap H$ are connected.  Then
the  exact homotopy sequence of
the quotient $\,K\cap H/Z_{K\cap H}(\c)$, implies that the group
$\,Z_{K\cap H}(\c)\,$  is connected,
as wished. This completes the proof of the lemma.
\end{proof}

\bn
\begin{remark} 
\label{SU21}
When $G=SO_0(2,1)$,  direct computations using the model
described in Example \ref{HYPER} show that the isotropy subgroup 
$\,Z_H(\c)\,$ is trivial.  Therefore principal orbits of type $\,G/Z_H(\c)\,$ are
diffeomorphic to
$\,SO_0(2,1)$ and  topologically equivalent to $\,SO(2)$.
In particular, they are not simply connected. 

Similarly, when $\,G= SO_0(3,1)\,$ the isotropy subgroup  $\,Z_H(\c)\,$ is
isomorphic to $\,SO_0(1,1)$, which is connected. Therefore principal orbits
of type
$$\,G/Z_H(\c) \cong SO_0(3,1)/SO_0(1,1) \,$$ 

\sn
are topologically equivalent to $\,SO(3)\,$ and are not simply connected. 

When  $\,G=SU(2,1)$,  
direct computations using the model  described in Example \ref{NONREDUCED1}
show  that  the isotropy subgroup 
$\,Z_{K\cap H}(\c)\,$ is isomorphic to  $\,S(U(1) \times U(1))$,   
which is connected.  
Principal orbits of type $\,G/Z_H(\c)\,$ are topologically equivalent to 
$\, K/Z_{K \cap H}(\c) \cong U(2)/U(1)  \cong SO(3)$. Hence they are
not simply connected.

Note that in all the above cases, despite the fact that the orbits
are not simply connected, the corresponding
isotropy subgroups are connected.
\qed
\end{remark}

\bn
\begin{lem}
\label{TOPORB3} All principal $G$-orbits of type $Z_{H'}(\c')$ are simply
connected.
\end{lem} 

\sn
\begin{proof}
An orbit of type $\,G/Z_{H'}(\c')\,$ is topologically equivalent to $\,K/Z_{H'\cap
K}(\c')$. We prove that the latter quotient is simply
connected by discussing each case separately.

Consider first  $\,G =SU(n,1) $. Direct computations using the model constructed in
Example \ref{NONREDUCED1} show  that
 $\,Z_{H'\cap K}(\c')\cong U(n-1)$.
 Hence the quotient $\,K/Z_{H'\cap K}(\c')\cong  U(n)/U(n-1)$ is diffeomorphic to the sphere  $ S^{2n-1}\,$.
 In particular, it is  simply connected for all $n\ge 2$.

Next let $\,G = Sp(n,1) $ or $\,G = F_4^* $. Both $G$ and $K$ are simply connected. 
So the quotient 
$\,K/Z_{H'\cap K}(\c')\,$ is simply connected provided that
$\,Z_{H'\cap K}(\c')\,$ is connected.
In order to show this,  denote by $\,K'\,$
 the maximal compact subgroup  of
$\,G'$ (see Sect. 4.2). Since $\,H'\,$ is contained in$\, G'$, the groups  $\,H'\cap K$ and $ H'\cap K'$ coincide and
 are both connected by Lemma \ref{LOCAL2}. 
Consider  the compact, rank-one, symmetric space
$\,K'/(K'\cap H')\subset G'/H'\,$. 
The non-zero  orbits of the 
isotropy representation of $\,K'\cap H'\,$ on $\,\k'\cap \q'$ 
 are of type $\,K'\cap H'/Z_{K'\cap H'}(\c')\,$ and  are
diffeomorphic to spheres  of dimension equal to 
$\,\dim\g^{2\alpha} -1 $.
Since  
 $\,\dim \g^{2\alpha}>2\,$ (cf. Table$\,$4.0),  they are simply connected.    
  As $H'\cap K'$ is connected,  by the  exact  homotopy sequence of the quotient $\,K'\cap H'/Z_{K'\cap H'}(\c')$, the groups 
$\,Z_{K'\cap H'}(\c')$ and $ Z_{K\cap H'}(\c')$ are also connected.
 It follows that the quotients  $\,K/Z_{H'\cap K}(\c')\,$  and $\,G/Z_{H'\cap K}(\c')\,$ are simply connected, as desired. \end{proof}

\nbigskip
\begin{lem}
\label{INIGENERAL}
Let  $\,q:\Sigma \to Z\,$ be  a $\,G$-equivariant Riemann  domain.  
Assume that every $\,z\,$ in
$\,Z\,$ admits an arbitrary small neighbourhood   $\,V\,$ and a
sequence $\,\{z_n\}\,$ converging 
to $\,z\,$ with the property that   
both the isotropy subgroups $\,G_{z_n}\,$ and the intersections 
$\,G\cdot z_n \cap V \,$ are connected.
Then $\,q\,$ is injective on every $\,G$-orbit of $\,\Sigma$. 
\end{lem}

\smallskip
\begin{proof}
Assume by contradiction that  the map $q$ is not injective on the $G$-orbit through
some  $\,\zeta\,$  in $\,\Sigma$.   Then there exists
$\,h \in G \,$ with  $\,h \cdot \zeta \not= \zeta\,$ such that 
$\,q(h \cdot \zeta)=q(\zeta) $. Since $q$ is locally injective, one can choose an open
neighborhood  $\,V$ of $z:=q(\zeta)$ in $Z$ as in the assumption, and open neighbourhoods  $\,W_\zeta\,$
and $\,W_{h \cdot \zeta}\,$ of  $\,\zeta\,$  and $h \cdot \zeta\,$ in $\Sigma$, such that 
$\,W_\zeta \cap W_{h \cdot \zeta}= \emptyset$ and  the restrictions $\,q|_{W_\zeta}:W_\zeta\to V\,$
and $\,q|_{W_{h \cdot \zeta}}:W_{h \cdot \zeta} \to V\,$ 
are bijective.  Then there exists  a sequence $\{z_n\}$ in $Z$, 
 converging   to $z$, with the property that both the isotropy subgroups $\,G_{z_n}\,$ and the intersections $\,G\cdot z_n \cap V \,$ are connected.
 
Consider  the sequence $\,\{\zeta_n:= ( q|_{W_\zeta})^{-1}(z_n)\}$ in $W_\zeta$. Since $\{\zeta_n\} $ converges to $ \zeta$, for
$\,n\,$ large enough,   the points 
$\,h\cdot \zeta_n\,$ lie in $ W_{h\cdot \zeta}$. Consequently
  their images 
$\,q(h\cdot \zeta_n)=h\cdot q(\zeta_n)=h\cdot z_n\,$ lie in 
$\,V$. Since both $\,G_{z_n}\,$ and $\,G\cdot z_n \cap V \,$
are connected,   the set $\,\Omega_n:=\{\,g \in G \ : \ g\cdot z_n \in V\}\,$
is connected. Note that both $\,e\,$ and $\,h\,$ belong to
$\, \Omega_n$. Hence there exists a continuous path 
$\,\gamma:[0, 1] \to \Omega_n\,$ with $\,\gamma (0)=e\,$ and $\,\gamma (1)=h$.
By the $\,G$-equivariance of $\,q\,$ both 
paths 
$$t \mapsto 
(q|_{W_\zeta})^{-1}(\gamma(t) \cdot z_n)
\quad {\rm  and}\quad  t \mapsto \gamma(t) \cdot \zeta_n$$ 
in $\,\Sigma\,$ are liftings of $\,t \mapsto \gamma(t) \cdot z_n$, with initial point
$\,\zeta_n$. On the other hand $\, (q|_{W_\zeta})^{-1}(\gamma(1) \cdot z_n)\in W_\zeta\,$
while
$\,\gamma(1) \cdot \zeta_n \in W_{h \cdot \zeta}$, giving a contradiction.
\end{proof}

\nbigskip
As a consequence of the previous lemmas one obtains the main result of this section.


\bn
\begin{prop}
\label{ORBITINJ} Let $\,G \, $ be a connected, non-compact,  real simple Lie group such
that the Riemannian symmetric space  $\,G/K\,$ has rank one.  Assume that $\,G\,$
is embedded in
its  universal complexification~$G^\C$ and is different from the groups $\,SL(2,\R)\,$
and 
$Spin(3,1)$.  
Let $q\colon \Sigma\to G^\C/K^\C$ be a $G$-equivariant Riemann domain.
Then $q$ is injective on every $G$-orbit.
\end{prop}

\begin{proof}
\medskip 
We begin by proving the following claim.

\sn
{\it Claim.} The isotropy subgroups of all principal $\,G$-orbits
are connected.

\sn
{\it Proof of the claim.} Since $G$ is connected, the isotropy subgroups
of simply connected orbits are necessarily connected.
Hence by Lemmas \ref{TOPORB1} -- \ref{TOPORB3} we only need
to discuss the isotropy types $\,Z_K(\a)\,$ when $\,G\,$ has Lie algebra 
$\,\s \o_0(2,1)\, $ and 
the isotropy  types $\,Z_H(\c)\,$ when $G$ has Lie algebra 
 $\,\s \o_0(2,1)$, $\,\s \o_0(3,1)\,$
and $\,\mathfrak s\mathfrak u(2,1)$.

Let $\,\g=\s\o(2,1)$. When $\,G=SO_0(2,1)\,$ the isotropy subgroups of all principal
$\,G$-orbits
are connected, by Remarks~\ref{SO21A}~and~\ref{SU21}. Observe that $\,SO_0(2,1)\,$  is
centerless and that $\,SL(2,\R)\,$ is a double covering of 
$\,SO_0(2,1)$. Since  the universal complexification of $\,SL(2,\R)\,$ is $\,SL(2,\C)$,
which is simply connected,  no covering of $\,SO_0(2,1)\,$ other  than $\,SL(2,\R)\,$ admits
an embedding into its universal complexification. Hence the claims follows for every group
$\,G\not=SL(2,\R)\,$ which has Lie algebra $\,\s\o(2,1)\,$ and embeds in its universal
complexification.

Let $\,\g=\s\o(3,1)$. When $\,G=SO_0(3,1)\,$ the isotropy subgroup $\,Z_H(\c)\,$ is
connected, by
Remark \ref{SU21}.  Note that $\,SO_0(3,1)\,$ is centerless and $\,Spin(3,1)\,$ is
the only non-trivial covering of $\,SO_0(3,1)\,$ which embeds in its
universal complexification. Hence the claims follows for every group
$\,G\not=Spin(3,1),$ which has Lie algebra $\,\s\o(3,1)\,$ and embeds in its universal
complexification.

Finally, let $\,\g=\mathfrak s \mathfrak u (2,1)$. When
$\,G=SU(2,1)$, the isotropy
subgroup
$\, Z_H(\c)\,$ is connected, by Remark~\ref{SU21}. Thus the same holds true
for every connected  real Lie group covered by
 $\,SU(2,1)$.
Since no covering group of $\,SU(2,1)\,$ admits an embedding in its universal
complexification, the claim holds true for every 
$\,G\,$ which has Lie algebra $\,\s\u(2,1)\,$ and embeds in its universal
complexification.
This  concludes the proof of the claim.

In order to complete the proof of the proposition,
recall that the union of  principal $G$-orbits forms an open dense subset
of $G^\C/K^\C$.
Hence, by the above claim every point in $\,G^\C/K^\C\,$ can be approximated by
points with connected isotropy subgroups. Due to this fact and the description
of the slice representation
at closed $G$-orbits (cf. Lemma \ref{ADJOINT} and diagrams in Sect. 4),  
all assumptions  of Lemma \ref{INIGENERAL} are met and the statement follows.
 \end{proof}

 \bn
 \begin{remark} \label{2ORBINJ} When $G=SL(2,\R)$,  the isotropy subgroups
of all
principal $G$-orbits in $G^\C/K^\C$ consist of the central elements
$\{\pm I_2\}$.
As we shall see in Example \ref{NONINJ1}, in this case there exist Stein,
$G$-equivariant Riemann domains
which are not injective on $G$-orbits. Similarly, one can  construct
$\,G$-equivariant
Riemann domains which are not injective on $G$-orbits in the case $G=Spin(3,1)$.
However, by Theorem
\ref{UNIVALENCE2} such Riemann domains cannot be Stein.
\end{remark}



\bn

\section{$G$-invariant Stein domains in $\,G^\C/K^\C$}

\bn

Let $\,G/K\,$ be a non-compact, rank-one, Riemannian symmetric space. 
In this section we exhibit a complete classification of Stein  $\,G$-invariant domains in $\,G^\C/K^\C$.
The main ingredient is the computation of the Levi form of
hypersurface $\,G$-orbits in
 $\,G^\C/K^\C$, which is carried out in \cite{Ge1} and in
Appendix 9.
Most of the Stein domains in our list are known. However, 
for $\,G=SU(n,1)\,$ we present some examples which  appear
to be new.  By working out an explicit  model  of $\,G^\C/K^\C$,
we show that they are all biholomorphic to $\,\B^n \times \C^n$.

\medskip
The classification result is stated for 
the the standard presentations  of $\,G/K\,$ given in Table 4.0.
This is no loss of generality, since by  Remark \ref{QUOTOP} 
the $\,G$-orbit structure of  $\,G^\C/K^\C\,$
as well as the $\,CR$-structure and topology of $\,G$-orbits  do not depend on the
presentation of the symmetric space $\,G/K$.

Retain  the notation used in  diagrams (\ref{DIAGRAM1}),
(\ref{DIAGRAM2}), (\ref{DIAGRAM3}), (\ref{DIAGRAM4}) of  Section 4.
Consider the 
$\,G$-invariant domains in $G^\C/K^\C$ defined by

\begin{equation} \label{DOMAINS}
\begin{matrix} 
D_1(a)=G\cdot(z_1\cup \ell_1((a,1))),\qquad D_2(a)=G\cdot(z_3\cup \ell_3((a,1))),\quad
\,0\leq a <1\,, 
    \cr
\cr
    S_1(b)=G\cdot \ell_2((b,\infty)),\qquad S_2(b)=G\cdot \ell_4((b,\infty)),\qquad \,0 \le b
<\infty\,. \cr 
\end{matrix}  
\end{equation}

 \nbigskip
\begin{theorem}
\label{STEIN} Let $\,G/K\,$ be a  non-compact, rank-one, Riemannian symmetric
space.  
All Stein $G$-invariant domains in $G^\C/K^\C$ are given by the following table.  

\bn
\noindent 
{\bf Table 6.0.}
\mn
\eightrm{

\centerline{
\offinterlineskip\hbox{\vbox{
\hrule\halign{&\vrule#&\strut\quad#\hfil\quad\cr
height2pt&\omit&&\omit&&\omit&&\omit&&\omit&\cr
&\quad\quad\qquad\qquad \qquad\quad $G=$ && $SO_0(2,1)$&&$SO_0(n,1)$&&$SU (n,1)$&&$Sp(n,1), \, n \ge 2$&\cr
height2pt&\omit&&\omit&&\omit&&\omit&&\omit&\cr
&DOMAIN\qquad &&  &&$n \ge 3$&&$n \ge 2$ &&$F_4^* $&\cr
height2pt&\omit&&\omit&&\omit&&\omit&&\omit&\cr
\noalign{\hrule}
height3pt&\omit&&\omit&&\omit&&\omit&&\omit&\cr
&$D_1(a),~0\leq a <1$&&Stein&&Stein&&Stein&&Stein&\cr
height3pt&\omit&&\omit&&\omit&&\omit&&\omit&\cr
&$ D_2(a),~0\leq a <1$&&Stein&&Stein&&no  &&no &\cr
height3pt&\omit&&\omit&&\omit&&\omit&&\omit&\cr
&$ S_1(b),~0\leq b <\infty$&&Stein&&no  && no &&no  &\cr
height3pt&\omit&&\omit&&\omit&&\omit&&\omit&\cr
&$ S_2(b),~0\leq b <\infty$&&Stein&&no  &&no  && no &\cr
height3pt&\omit&&\omit&&\omit&&\omit&&\omit&\cr
&$ D_1(0)\cup G \cdot w_1\cup S_1(0)$&&Stein&&no  &&Stein&&no  &\cr
height3pt&\omit&&\omit&&\omit&&\omit&&\omit&\cr
&$ D_1(0)\cup G \cdot w_4\cup S_2(0)$&&Stein&&no  &&Stein&&no  &\cr
height3pt&\omit&&\omit&&\omit&&\omit&&\omit&\cr
& $ D_2(0)\cup G \cdot w_2\cup S_1(0)$&&Stein&&no  &&no  && no &\cr
height3pt&\omit&&\omit&&\omit&&\omit&&\omit&\cr
&$ D_2(0)\cup G \cdot w_3\cup S_2(0)$&&Stein&&no  &&no  && no &\cr
height3pt&\omit&&\omit&&\omit&&\omit&&\omit&\cr}
\hrule}}}
}
\end{theorem}

\bn
{\bf Remark.} The domains $D_1(0)$ and $D_2(0)$ are known as   Akhiezer-Gindinkin domains. They were introduced in \cite{AkGi}  for $G/K$ of arbitrary rank.
In the two-dimensional case, the domains $S_1(0)$ and $S_2(0)$  are related to the causal structure of the symmetric space $G/H= SO_0(2,1)/SO(1,1)$. Domains of this type were studied by Neeb  in  \cite{Ne}.

\bigskip

\noindent
{\it Proof of the theorem.} We first show that all the domains listed in the above table are Stein.
The Akhiezer-Gindikin domain $\,D_1(0)\,$ is  Stein by \cite{BHH}. 
For $\,0<a<1$, the domains 
$\,D_1(a)$ are  $\,G$-invariant subdomains
of  $\,D_1(0)\,$ containing the minimal orbit $\,G \cdot z_1 \cong G/K$.
Their Steiness follows for example from the 
non-linear convexity theorem in \cite{GiKr}.

When  $\,G= \,SO_0(n,1) $, with $n\ge 2$,  the domain $\, D_2(0)\,$ and its subdomains 
 $\, D_2(a)\,$, for  $\,0 < a < 1$, are Stein since  they are 
biholomorphic to $\,D_1(0)$ and $\, D_1(a) $, respectively.
One such biholomorphism is given for example 
by the map 
 
$$G^\C/K^\C \to G^\C/K^\C\,, \quad gK^\C \to g_3gK^\C\,,$$

\nmedskip
where $g_3= \exp iA_3$, with $\alpha(A_3)={\pi\over 2}$ (cf. (\ref{BASEPOINT1})).
Note that  $g_3\in  \,SO(n,1) \setminus \{SO_0(n,1)\}\,$, therefore
$\,g_3G=Gg_3$. As a consequence,
the above map  exchanges the singular orbits $G\cdot z_1$ and $G\cdot z_3$ and
maps $G\cdot \ell_1(a)$ onto $G\cdot \ell_3(a)$, for $\,0 < a < 1$.

When  $G=SO_0(2,1)$, the domains $S_1(0)$ and $S_2(0)$  and their subdomains
 $\,S_1(b)\,$ and 
$\,S_2(b)$, for $\,0<b<\infty$,
were shown to be
Stein  in \cite{Ne}.

The last four domains in the list contain  in their interior one of the
non-closed orbits $G\cdot w_i$, for some $i=1,\ldots,4$.
Their boundary  consists of  two
non closed $G$-orbits and the singular orbit in their closure.  All of them are Stein if
$G=SO_0(2,1)\cong SU(1,1)/\{\pm I_2\}$. Only  $D_1(0)\cup G \cdot w_1\cup S_1(0)$ and
$D_1(0)\cup G \cdot w_4\cup S_2(0)$ are Stein, when $G=SU(n,1)$, with $n>1$. These facts
are proved  in Example \ref{NONREDUCED2}  by constructing  explicit models of such
domains.

 \sn
 
In order to complete the classification, it remains to show that no $G$-invariant domains in $\,G^\C/K^\C\,$ are Stein other than the ones listed in  Table 6.0.
When $\,G=SO_0(2,1) \cong SU(1,1)/\{\pm Id\}$ and 
$\,G=SU(n,1)$, with $\,n \ge 2$, this is proved 
in Example \ref{NONREDUCED2}.

In all other cases, namely $\,SO_0(n,1)$, with $\,n > 1$, $\,Sp(n,1)$ and $F^*_4$, this 
follows  from the description of the $\,G$-orbit
space of $\,G^\C/K^\C\,$ given in diagrams (\ref{DIAGRAM2}),  (\ref{DIAGRAM3}), (\ref{DIAGRAM4}) 
and  from the computation of the Levi form of the  hypersurface $\,G$-orbits in $G^\C/K^\C$. Indeed,
by Prop.$\,$5.6 and Prop.$\,$5.21 in \cite{Ge1},  all principal orbits
have indefinite Levi form, except for the ones intersecting the slice $\,\ell_1\,$
(the domain $\,D_1(a)\,$ is Stein)
and,  only when the restricted root system of $\,\g\,$ is reduced,  the slice
$\ell_3$ (the domain $\,D_2(a)\,$ is Stein for $\,G=SO_0(n,1)\,$ ).  Moreover,  by Remarks
\ref{REMW23}  and \ref{REMW2} the Levi form of the non-closed hypersurface orbits  
$G\cdot w_2$ and $G\cdot w_5$ is indefinite. Since the boundary of a Stein
domain cannot
have indefinite Levi form, the theorem follows.

 \qed

\bigskip

Let us  illustrate the result 
 of Theorem \ref{STEIN} on the model of $G^\C/K^\C$ described in Example
 \ref{HYPER}. The Stein, $\,G$-invariant domains are studied by means of
 an appropriate
$G$-invariant function on $G^\C/K^\C$.


\bn
\begin{exa} \label{HYPER2}  \rm
  Let $G=SO_0(n,1)$. By  Example \ref{HYPER}, the quotient
   $G^\C/K^\C$ can be identified
with $M^\C:=\{\xi\in\C^{n+1}~:~\xi_1^2+\ldots+\xi_n^2-\xi_{n+1}^2=-1\} $. Assume $n > 2$.
Consider 
 the $\,G$-invariant
function $\,f: M^\C \to \R\,$ defined  
 by $$\,f(\xi_1, \dots ,\xi_{n+1}):=|\xi_1|^2+\dots +|\xi_n|^2-|\xi_{n+1}|^2 -1.$$
For every  $\,0<a<1$, the $\,G$-invariant domains
$\,D_1(a)\,$ and $\,D_2(a)$ coincide with the two connected components of 
the set  $\,\{\,\xi\in M^\C~:~  f(\xi)<r\,\}\,$, for some  $\,-2<r<0\,$.  Every  such
domain   is bounded by a single $\,G$-orbit on which the Levi form of $f$ is
positive definite.  Hence it is Stein. 

The  $\,G$-invariant domains
$\,D_1(0)\,$ and $\,D_2(0)$ coincide with the two connected components 
of the set  $\,\{\,\xi\in M^\C~:~  f(\xi)<0\,\}\,$.  They are bounded  by  the
non-smooth hypersurfaces 
$\, \partial D_1(0)=G\cdot (z_2\cup w_1)\,$ and $\, \partial D_2(0)=G\cdot (z_2\cup w_2)\,$, 
respectively. 
At all smooth points of $\,\partial D_1(0)\,$ and $\,\partial D_2(0)\,$
the Levi form of $\,f\,$ has  $\,n-2\,$ positive eigenvalues  and one zero eigenvalue. 
This is consistent with the fact that $\,D_1(0)\,$ and $\,D_2(0)$ are Stein.
The Levi form of $f$ is indefinite on all remaining hypersurface $G$-orbits.
Thus there are no other Stein $G$-invariant domains in $M^\C$.  
\qed
\end{exa}

\bn

Next we determine all Stein, $\,G$-invariant domains in $\,G^\C/K^\C \,$ in
the case $\, G=SU(n,1)\,$ by using the model of $\,G^\C/K^\C \,$ described in
Example \ref{NONREDUCED1} and Remark \ref{N=1}. This settles the missing cases in the proof of Theorem
\ref{STEIN}.


\bn
\begin{exa}
\label{NONREDUCED2}{\rm
Let $G=SU(n,1)$, with $\,n \ge 1$.
By   Example \ref{NONREDUCED1} the quotient $G^\C/K^\C$ can be identified with
$M^\C:=\P^n \times  \overline{\P}^n \setminus\{ \langle z,w\rangle_{n,1}=0\} $. 
Consider the $G$-invariant  function
 $f\colon M^\C \longrightarrow \R$ defined by
 $$f([z],[w])=-{{\langle z,z\rangle_{n,1}\langle w,w\rangle_{n,1}}\over { |\langle z,w\rangle_{n,1}|^2}}.$$

\noindent
$\bullet$ $\,G=SU(1,1)$

\nmedskip
 By computing the Levi form of $\,f\,$ on the 
 $\,G$-orbits in the level set $\,\{f=r\}$, with $\,r<0$, one shows that 
 the domains $\,D_1(a)$ and $\,D_2(a)\,$ are Stein for all $\,0<a<1$.
 Similarly one shows that $\,S_1(b)\,$ and $\,S_2(b)$ are Stein
 for every  $\,b>0$.
One can also verify  that the Levi form of $f$ on all non-closed hypersurface orbits $G\cdot w_1,\ldots, G\cdot w_4$ is identically
zero.
This is consistent with the fact that the domains $\,D_1(0)$, $\,D_2(0)$,
$\,S_1(0)$ and $\,S_2(0)\, $ are Stein.
 We claim that the domains 
 $$W_{1,1}:=D_1(0)\cup G \cdot w_1\cup S_1(0), \quad W_{1,2}:=D_1(0)\cup G \cdot w_4\cup S_2(0),$$
 $$W_{2,1}:=D_2(0)\cup G \cdot w_2\cup S_1(0), \quad W_{2,2}:=D_1(0)\cup G \cdot w_3\cup S_2(0),$$
 
 \nmedskip
 are Stein as well.
 By evaluating the hermitian forms
 $\,\langle z,  z \rangle _{n,1}\,$ and $\,\langle w,  w \rangle _{n,1}\,$ on
 the slices described in Example \ref{NONREDUCED1} and Remark \ref{N=1},
 one sees that such domains can be characterized as follows
$$W_{1,1} = \,\{\,
\langle z, w\rangle _{1,1} \not= 0 \ {\rm and} \  
  \langle z,  z\rangle _{1,1}<0\,\}\,,$$
$$\,W_{1,2} = \,\{\,
\langle z,  w\rangle _{1,1} \not= 0 \ {\rm and} \ 
  \langle w, w\rangle _{1,1}<0\,\}\,,$$
$$W_{2,1} = \,\{\,
\langle z, w \rangle _{1,1} \not= 0 \ {\rm and} \  
  \langle w,  w \rangle _{1,1}>0\,\}\,,$$
$$\,\, W_{2,2} = \,\{\,
\langle z,  w\rangle _{1,1} \not= 0 \ {\rm and} \ 
  \langle z, z\rangle _{1,1}>0\,\}\,.$$

\smallskip
  \noindent
As a consequence, the maps defined by
$$\Delta \times \C \to W_{1,1}, \quad (u,v) \to ([u:1],[\bar v :  1+ \bar u\bar v])\,,$$
 $$\C \times \Delta \to W_{1,2}, \quad (u,v) \to
 ([u:1+ u v],[ \bar v:1])\,,$$
 $$\Delta \times \C \to W_{2,1}, \quad (u,v) \to ([1+ uv : u],[1: \bar v])$$
 $$\C \times \Delta \to W_{1,2}, \quad (u,v) \to
 [1: u],[1+ \bar u \bar v: \bar v])\,, $$
 
 \nsmallskip
 are biholomorphisms. Here
 $\,\Delta\,$ denotes the unit disk in $\,\C$.
In particular the domains $\,W_{1,1} \ldots W_{2,2}\,$ are Stein,
 as claimed.

Other $G$-domains in $M^\C$ which are possibly Stein
can only be obtained as arbitrary unions of domains $\,W_{k,l}$,
for $\,k,l=1,2$.
We claim that such unions are not Stein. For instance, let us 
show that  $\,W_{1,1} \cup W_{2,1}\,$ is not Stein. Consider 
the Stein local chart
$$\phi: \C^2 \to \P^1 \times \overline \P^1 , \quad (u,v) \to ([u:1],[1:\bar v])\, .$$

\sn
Since the preimage
$$\,\phi^{-1}(W_{1,1} \cup W_{2,1})=\{\,(u,v) \in \C^2 \ : \ u\not=v \
{\rm and  \ either}\  |u|<1 \ {\rm or} \  |v|<1\,\}$$

\sn
is not Stein,
the domain  $\,W_{1,1} \cup W_{2,1}\,$ is not Stein either.
An analogous argument applies to the remaining cases.
  
\nbigskip
$\bullet$ $\,\ G=SU(n,1)$, with $\,n \ge 1$. 

\medskip
Using the $G$-invariant function $f$, one can  prove that the domains 
$D_1(a)$, for $a>0$, are Stein. One  can also verify that $D_1(0)$ coincides with a
connected component of the set $\{z\in M^\C~|~f(z)<0\}$ and that on the smooth part of
its boundary $\partial D_1(0)=G\cdot (w_1\cup z_2\cup w_4)$ the Levi form of $f$ has
non-negative eigenvalues. This is consistent with the fact that $D_1(0)$ is Stein. 
 
Moreover, the  Levi form of $f$ on the principal $G$-orbits  through the slices
 $\ell_2$, $\ell_3$, $\ell_4$, $\ell_5$ and on the non-closed hypersurface orbit $G\cdot w_5$  is indefinite. On the other hand, the  Levi form of $f$ is definite
 on the non-closed hypersurface orbits $G\cdot w_2$ and $G\cdot w_3$.
As a result, other $G$-invariant domains in $M^\C$ which are possibly Stein are only
  $$ W_{1,1}:=  D_1(0)\cup G \cdot w_1\cup S_1(0)\,,\quad  
 W_{1,2}:=D_1(0)\cup G \cdot w_4 \cup S_2(0)\,, \quad  W_{1,1} \cup W_{1,2} .$$
 First we show that $\,W_{1,1}\,$ and $\,W_{1,2}\,$ are indeed Stein.
By evaluating $\,<z,  z>_{n,1}\,$ and $\,<w,  w>_{n,1}\,$ on the slices described in Example  \ref{NONREDUCED1},  one sees that such domains can be characterized as follows
$$ W_{1,1} = \,\{\,([z],[  w]) \in \P^n \times \overline {\P}^n \ :\ 
<z, w>_{n,1} \not= 0 \ {\rm and} \  
  <z,  z>_{n,1}<0\,\}\,$$
$$ W_{1,2} = \,\{\,([z],[  w]) \in \P^n \times \overline {\P}^n \ :\ 
<z,  w>_{n,1} \not= 0 \ {\rm and} \ 
  <w, w>_{n,1}<0\,\}\,.$$
\smallskip
  \noindent

As a consequence the  maps
$$\B^n \times \C^n \to W_{1,1}, \quad (u,v) \to ([u:1],[v:1+\bar u_1\bar v_1 +
 \dots + \bar u_n\bar v_n])$$
 $$\C^n \times \B^n \to W_{1,2}, \quad (u,v) \to
 ([u:1+u_1 v_1 +
 \dots + u_n v_n ],[ \bar v:1])\,,$$

\sn
 are biholomorphisms. Here  $\, \B^n\,$ denotes the unit ball in  $\,\C^n$.
 In particular $\,W_{1,1}\,$ and $\,W_{1,2}\,$ are Stein, as claimed.

Next we show that the domain $\,\Omega:=W_{1,1} \cup W_{1,2}$, with boundary $\,\partial \Omega = G\cdot (w_2 \cup z_2 \cup w_3)$,
is not Stein. 
Assume by contradiction
that $\,\Omega\,$ is Stein.  Let $\,\c'\,$ be the abelian subalgebra generating
the Cartan subset $\,\mathcal C'\,$ (cf. Example
\ref{NONREDUCED1}). Let $\,T=\exp \c'\,$ be the corresponding compact torus in~$\,G\,$.  Consider the $T$-action on $\,\Omega\,$
and the induced local holomorphic $T^\C$-action.
By the globalization theorem in \cite{He1}, Section 6.6, the domain
$\,\Omega\,$ embeds
in its $T^\C$-globalization $\Omega^*$  as a $T$-invariant,  orbit-convex subset. 
By definition, this means that  the intersection of $\, \Omega \,$ with an $\exp i\c $-orbit in $\,\Omega^*\,$ is  connected.

Every $\,T^\C$-orbit through the slice $\,\ell_1\,$ is contained in
$\,\Omega$. 
Indeed in $\,M^\C\,$ one 
can verify that $\,\exp (isC' )\cdot \ell_1(t) =\,$
$$([0: \ldots:e^s i \sin \frac{\pi}{4}(1 - t):
e^{-s} \cos \frac{\pi}{4}(1 - t)], [0: \ldots  : -e^{-s} i \sin \frac{\pi}{4}(1 - t):
e^{s} \cos \frac{\pi}{4}(1 - t)]).$$

\nsmallskip
Thus for fixed $\,0<t<1\,$,
the function $\, \R \to \R$, defined by $\, s \mapsto f(\exp (isC' )\cdot \ell_1(t))$
is given by 
$$ (\,e^{2s} \sin^2 \frac{\pi}{4}(1 - t)- 
e^{-2s} \cos^2 \frac{\pi}{4}(1 - t)\, ) (\,e^{-2s} \sin^2 \frac{\pi}{4}(1 - t) -
e^{2s} \cos^2 \frac{\pi}{4}(1 - t)\, )\,$$
 
\nsmallskip
and vanishes 
exactly twice, namely on $\,G\cdot w_1\,$ and on $\,G\cdot w_4$.
It follows that $\,\exp (i\c ' )\cdot \ell_1(t)\,$ never crosses the boundary
of  $\,\Omega $ and consequently the complex orbit
 $\,T^\C \cdot \ell_1(t)\,$ 
is entirely contained in $\,\Omega$, as claimed. Moreover, for  every fixed $\,s
>0\,$,  one has
$$\lim_{n \to \infty} \exp (isC') \cdot  \ell_1(1/n) = \ell_{2}(s) \in \Omega,$$
$$\lim_{n \to \infty} \exp (-isC') \cdot \ell_1(1/n) =  \ell_{4}(s) \in \Omega.$$

\nmedskip
Then the orbit-convexity of $\,\Omega\,$ in $\Omega^*$
 implies that the sequence
 $\,\{ \ell_1(1/n) \}_n\,$ has a limit point in  $\,\Omega$.
 On the other hand, in $\,G^\C/K^\C\,$
 one has $\,\lim_n \ell_1(1/n) = z_2$, which is not in  $\,\Omega$.
 This yields a contradiction and proves that $\,\Omega\,$ is not Stein.
 The classification of all Stein $\,G$-invariant domains in $\,M^\C\,$ is
 now complete.
 \qed
 }\end{exa}


\bn

We conclude this section with a remark which is 
a consequence of Theorem \ref{STEIN}
and  is often used in the sequel.  

\bn

\begin{remark}
\label{STEINREMARK} 
{\rm 
Let $\,D\,$ be a domain in $\,G^\C/K^\C\,$ with smooth boundary 
$\,\partial D$. It is well known that if $\,D\,$ is not pseudoconvex 
at $\,z \in \partial D$, then no holomorphic function on $\,D\,$  diverges in the
vicinity of $\,z$. Let $\ell\colon I \to G^\C/K^\C$ be a  slice for 
principal $G$-orbits
in $G^\C/K^\C$.  By the classification of Stein, $G$-invariant domains in $G^\C/K^\C$
given  in Theorem \ref{STEIN}, the following facts hold.

\sn
\item{$(i)$} Assume that the Levi form of the orbits parametrized by $\ell$ is definite. Let $(c,d)\subset I$ be an interval
with $0\le c< d$ and $d\in I$. Then no holomorphic function on the 
invariant domain $G\cdot \ell((c,d))$
diverges in the vicinity of the boundary orbit $G\cdot \ell(d)$
($\,$for instance, if $\,I=(0,1)$, the domain $\,D_1(d)\,$ is strictly pseudoconvex
at every point of the boundary orbit $\,G \cdot \ell_1(d)$. Thus the domain 
$\,G\cdot \ell_1((c,d))\,$ is not pseudoconvex at every point of $\,G \cdot \ell_1(d)\,$).

\sn
\item{$(ii)$} Assume that the Levi form of the orbits  parametrized by $\ell$ is indefinite. Let $(c,d)\subset I$ be an interval, with $c\in I$. Then no holomorphic function on the invariant domain $G\cdot
\ell((c,d))$ diverges in the vicinity of the  boundary orbit  $G\cdot \ell(c)$. Similarly, if  $d\in I$, then no holomorphic function diverges in the vicinity of  $G\cdot \ell(d)$.
 \qed
 }
\end{remark}


\bn
\section{Univalence  over $G^\C/K^\C$}

\bn

Let $G $ be a connected, non-compact,
real simple Lie group, $K \subset G$ a
maximal compact subgroup  and $G^\C$ the universal complexification of $G$. 
Assume that the center $\Gamma$ of $G$ is finite and  $G$ is not a covering of
$SL(2,\R)$.  
 In this section we show that
a holomorphically separable, $\,G$-equivariant Riemann  domain 
$\,q:\Sigma \to G^\C/K^C\,$   is necessarily univalent, if the rank of $G/K$ is equal to one
(cf.~Thm.~\ref{UNIVALENCE2} and Rem.~\ref{FINITECENTER}). 

In most cases the map $q$ is  injective on every $G$-orbit (cf. Sect.$\,$5).  So 
 we are  reduced to prove the injectivity of $q$ over the global slices for the $G$-action 
defined by  diagrams (\ref{DIAGRAM1}),   (\ref{DIAGRAM2}), (\ref{DIAGRAM3}), (\ref{DIAGRAM4}) in Section 4.
 Recall that the slices parametrizing principal $G$-orbits are diffeomorphic to open intervals of $\, \R \,$ and that  a local diffeomorphism  of a one-dimensional 
smooth manifold  into 
 the real line $\,\R\,$ is necessarily injective. As a consequence,  $q$ is injective on every connected component  of $\Sigma$ over a domain  in $G^\C/K^\C$ consisting of principal orbits.
 
  However, 
in order to ensure monodromy around
  the singular   orbit  $G\cdot z_2$ (cf.~diagrams in Sect.$\,$4),  it is necessary  to 
combine the uniqueness property of path liftings for Riemann domains with the complex 
 geometry of the $\,G$-invariant domains in $\,G^\C/K^\C$. 
Before proving  the main result of this section
some preliminary  lemmas are needed.

\bigskip
Let $\ell\colon I \to G^\C/K^\C$ be one of the slices 
for principal $G$-orbits  defined in (\ref{SLICE13}), (\ref{SLICE24}), (\ref{SLICENR13}), (\ref{SLICENR24})
and
(\ref{SLICENR5}) of Section 4.
Define 
\begin{equation}
\label{EXTENDED}
\hbox{ $\widehat I:= (0,1],\ $ when $\ I=(0,1)\,,$ \quad and  \quad $\widehat I:=I,\ $ when $\ I=\R^{>0}$.}
\end{equation}

\noindent
 Recall that $\,I=(0,1)\,$ only when  $\ell =\ell_1$ or $\ell =\ell_3$. In those cases extend $\ell$ to   $\widehat I=(0,1]$ by defining
  \begin{equation}
\label{TOTREALORB}
\ell_1(1):=e K^\C,\qquad \ell_3(1):=\exp (iA_3) K^\C \,.
\end{equation}
 
 \noindent
 We refer to $\,\ell\colon \widehat I \to G^\C/K^\C\,$ as an {\it extended slice}. Note that 
the images of the extended slices  $\,\ell_1\,$ and  $\,\ell_3\,$ include
the points  $\,z_1\,$ and 
$\, z_3\,$, respectively. 

Let $\,q:\Sigma \to G^\C/K^\C\,$ be a  $\,G$-equivariant Riemann domain and let 
$\ell \colon \widehat I \to G^\C/K^\C$   be an extended slice.
A {\it local lifting} of
$\ell$ is a smooth path  
$$  \widetilde \ell \colon J\longrightarrow \Sigma,$$
defined on  a   non-empty  interval   $\,J\,$ open in $\,\widehat  I$,
and satisfying the condition 
  $q\circ   \widetilde \ell =\ell $ on $J$.
A local lifting $  \widetilde \ell \colon J\longrightarrow \Sigma$ is {\it maximal} if
it cannot be extended to a larger   interval  $J'$, with~$J\varsubsetneq J'\subset
\widehat I$.

\nbigskip
\begin{lem}
\label{EXTENSION} Assume that $G$ is embedded in its  universal complexification~$G^\C$
and is different from $SL(2,\R)$ and  $Spin(3,1)$. Let $\,q:\Sigma \to G^\C/K^\C\,$ be a Stein, $\,G$-equivariant Riemann domain and  
 let $  \widetilde \ell\colon J\to \Sigma$ be a maximal local lifting of an extended slice $\ell\colon \widehat I\to G^\C/K^\C$.

\sn
\item{(i)} if the Levi form of the principal orbits parametrized by $\ell$ is definite, then  the   invariant  domain
$G\cdot \ell(J)$   in  $G^\C/K^\C$ is Stein (see Theorem \ref{STEIN}).
 \sn
\item{(ii)} If the Levi form of the principal orbits parametrized by  $\ell$ is indefinite, then $\,J\,$ coincides with $\,\widehat I$. 

\end{lem}

\smallskip
\begin{proof}  
$(i)$ Consider first  the case $\widehat I=\R^{>0}$ ( see diagram \ref{DIAGRAM1},
Ex.~\ref{NONREDUCED2} and
Rem.~\ref{N=1} ). By Theorem \ref{STEIN}, we need to show that $J=(b,+\infty)$, for some $b\ge 0$. Assume by 
contradiction that $J=(b,d)$, with $0\le b<d<\infty$.  
Since the lifting  $  \widetilde\ell(J)$ is a one-dimensional real-analytic submanifold of $\Sigma$,  
 the local diffeomorphism $q|_{  \widetilde\ell(J)}$ is injective. By Proposition \ref{ORBITINJ}, the map $q$ is injective on
every $G$-orbit. Therefore the restriction $q|_{G\cdot   \widetilde\ell(J)}\colon G\cdot
  \widetilde\ell(J) \to G\cdot  \ell(J)$ is a biholomorphism. 

By the maximality of $  \widetilde \ell$,  when $n$ grows,  the sequence  
$\,\{   \widetilde  \ell(d-1/n)\}_n\,$  leaves every given 
compact subset in $\,\Sigma$. Since  $\Sigma$ is Stein,   there
exists a holomorphic function $\,f \in \O(\Sigma)\,$ such that  $\,\lim_{n \to \infty}
|f(  \widetilde  \ell(d-1/n))|= \infty.$ 
 
On the other hand,  the push-forward of  $f$ by $q|_{G\cdot   \widetilde\ell(J)}$ defines a 
holomorphic function  in $\,  \O(G\cdot  \ell(J))\,$  which diverges  in the vicinity of the boundary orbit
$G\cdot \ell(d)$. This contradicts $(ii)$ of Remark \ref{STEINREMARK}, implying that $J$ is of the form $(b,\infty)$, as claimed.

Consider now the case $\widehat I=(0,1]$. This only occurs for $\ell=\ell_1$ or, when the restricted root system of $\g$ is reduced, for  $\ell=\ell_3$ (cf.
diagrams in Sect.~4 and  \cite{Ge1}, Prop. 5.6).
By Theorem \ref{STEIN}, we need to show that $J=(a,1]$,  for  some $a\ge 0$. Assume by contradiction that $J=(a,d)$ with $0\le a<d \le 1$.  
The  argument used in the previous case shows that $J=(a,1)$   and  there exists a holomorphic  function $\,f \in \O(\Sigma)\,$ such that  $\,\lim_{n \to \infty}
|f(  \widetilde  \ell(1-1/n))|= \infty.$ 
 Moreover,  the push-forward of  $f$ by $q|_{G\cdot   \widetilde\ell(J)}$ defines a holomorphic function $\,\bar f 
\in \O(G\cdot  \ell(J))\,$ which diverges in the vicinity of the boundary orbit $G\cdot \ell(1)$. On the other hand,  such 
orbit is a totally real submanifold of   $G\cdot  \ell((a,1])$. Thus $\bar f$ extends to a
holomorphic function on  $G\cdot  \ell((a,1])$. This yields a contradiction, implying that $J=(a,1]$, as desired. 

$(ii)$ Assume first that $\widehat I=\R^{>0}$. Then $(ii)$ of Remark \ref{STEINREMARK} and an analogous argument as  in the proof of $(i)$  show  that $J=\widehat I$.  
Consider then the case  $\widehat I=(0,1]$. This only occurs when the restricted root system of $\g$ is non-reduced and $\ell=\ell_3$ (cf.
diagrams in Sect.~4 and  \cite{Ge1}, Prop. 5.6).
A similar argument as in the proof of $(i)$   shows that if a lifting $  \widetilde \ell_3\colon J\to \Sigma$ is maximal, then  either $J=(0,1]$ or $J=(0,1)$. 

In order to prove that $J=(0,1]$, suppose by contradiction that  
$J=(0,1)$.  Consider a sequence
 $\{z_n\}$  in $G\cdot \ell_3(J)$, converging to a point on the boundary orbit
$G \cdot w_5$, say $\,w_5$.  Since the Levi form of $G\cdot w_5$ is indefinite (see
Remark \ref{REMW23}), 
 no holomorphic function on $G\cdot \ell_3(J)$ diverges on $\{z_n\}$. 
Note that the restriction $q|_{G\cdot   \widetilde \ell_3(J)}\colon G\cdot   \widetilde \ell_3(J)\to G\cdot \ell_3(J)$ is a biholomorphism. Hence no holomorphic function of $G\cdot   \widetilde \ell_3(J)$ diverges on the sequence $\{\zeta_n\}$ in $\Sigma$, defined by $\zeta_n:=(q|_{G\cdot   \widetilde \ell_3(J)})^{-1} (z_n)$. By the Steinness of $\Sigma$,  there exists a subsequence of $\{\zeta_n\}$ converging to a point $\eta_5$ in $\Sigma$. Since $q$ is continuous, one has $q(\eta_5)=w_5$. 

By the $\,G$-equivariance of $\,q$, the description of the slice representation at
$\,z_3\,$ given in Remark \ref {ADJOINT} 
and Proposition \ref{ORBITINJ}, there exists  a $\,G$-invariant neighbourhood
$\,V\,$ of
$\,\eta_5\,$ in
$\,\Sigma\,$  on which
$\,q\,$ is injective. Its  image
$q(V)$ intersects the slice $\ell_5$ in $\ell_5((0,\epsilon))$, for some $\epsilon>0$.  By statement $(i)$ of this lemma, the local lifting $s\mapsto (q|_V)^{-1} (\ell_5(s))$, with $s\in (0,\epsilon)$, extends to a lifting
$   \widetilde \ell_5\colon \widehat
I_5\to
\Sigma$ of $\ell_5$. Note that $q$ maps the $G$-invariant domain $W:=G\cdot (   \widetilde \ell_3(J)\cup \eta_5 \cup
  \widetilde
\ell_5(\widehat I_5))$ in
$\Sigma$ biholomorphically onto the domain $q(W)=G\cdot (  \ell_3(J)\cup w_5 \cup  \ell_5(\widehat I_5))$ in $G^\C/K^\C$. Since
$G\cdot \ell_3(1)$ is a  totally real submanifold of   $q(W) \cup G\cdot  \ell_3(1)  $
(see Lemma 2.11 and Rem.~2.13 in \cite{Ge1}),   every holomorphic function on $q(W)$  extends to a
holomorphic function on  $q(W)\cup  G\cdot  \ell_3(1) $.  As a consequence, no
holomorphic function on $ W $  can diverge  on  the sequence 
$\{  \widetilde\ell_3(1-{1\over n})\}_n$ in $\Sigma$.

On the other hand, 
by the maximality of $  \widetilde \ell_3$, when $n$ grows  the sequence  
$\,\{   \widetilde  \ell(1-1/n)\}_n$   leaves every given 
compact subset in $\,\Sigma$. Since  $\Sigma$ is Stein,  there
exists a holomorphic function $\,f \in \O(\Sigma)\,$ such that  $\,\lim_{n \to \infty}
|f(  \widetilde  \ell_3(b-1/n))|= \infty.$ 
This yields a contradiction, implying that $J$ necessarily coincides with $(0,1]$, as claimed.

\end{proof}


\sn

Let  $\ell_1$ and $\ell_3$ be  the slices  parametrizing the principal orbits through the fundamental Cartan subset $\mathcal A$.  Denote by  $\,\mathcal C=\exp i \c \cdot z_2\,$ the standard Cartan subspace with base point  $\,z_2$ and define $\,\mathcal C^*:=\mathcal C \setminus \{z_2\}\,$.  Recall that in the reduced case, $\c=\R(X+\theta(X))$, for some  non-zero vector  $X\in\g^\alpha$, and $z_2=\exp (iA_2)K^\C$ with $\alpha(A_2)=\pi/2$. In the non-reduced case, $\c=\R(X+\theta(X))$, for some  non-zero vector  $\,X\in\g^{2\alpha} \,$, and $\,z_2=\exp (iA_2)K^\C\,$  with $\,\alpha(A_2)=\pi/4$.
In both cases, $\,\exp \c\,$ is a compact, one-dimensional,
real torus in $G$ which we denote by $\,T$. Both  $\,T\,$ and its universal complexification $\,T^\C \cong \C^*\,$ act on  $\,G^\C/K^\C\,$ by left translations.

In the next proposition, we single out two distinguished
$\,G$-invariant domains $\Omega$ and
$\Omega'$ in
$G^\C/K^\C$ containing all $T^\C$-orbits through the slices $\ell_1(I_1)$ and $\ell_3(I_3)$, respectively. 


\nbigskip
\begin{lem}
\label{GLOBALORB}
Let  $\,G/K\,$ be a non-compact, rank-one,  Riemannian symmetric space.
Consider the domain in  $G^\C/K^\C$ defined by 
$$\, \Omega := G \cdot (\,\ell_{1}(I_{1}) \cup w_{1} \cup w_4\cup \mathcal C^*\,).$$ Then for every point $z\in
\ell_{1}(I_{1}) $, the complex orbit $T^\C\cdot z$ is contained in
$\Omega$.
  \pn 
 Similarly, define 
 $$\, \Omega':= G \cdot (\,\ell_{3}(I_{3}) \cup w_{2} 
\cup w_3\cup \mathcal C^*\,).$$

\sn
Then for every point $z\in \ell_{3}(I_{3}) $, the complex orbit $T^\C\cdot z$ is contained in $\Omega'$.
\end{lem}

\smallskip
\begin{proof}
We first assume that $G=SO_0(2,1)$ and prove the statement by using the model
$\,M^\C\,$ of  $G^\C/K^\C$ constructed in Example \ref{HYPER}.
Let  $C$ be the generator of $\c$ chosen there.
Then, for  $s\in\R$ and $t\in (0,1)$,
one has 
$$\,\exp (isC )\cdot \ell_1(t) = (\,\sinh (2s) \sin \frac{\pi}{2}(1 - t),
\, i \cosh (2s) \sin \frac{\pi}{2}(1- t),  \, \cos \frac{\pi}{2}(1 - t)\,)\,.$$

Since $\,z_2 = (\,0, i,0 \,)\,$ and the entries of the matrix group $\,G\,$
are real, from the above expression one   easily verifies that 
$\,\exp i\c \cdot \ell_1(I_1) \cap G \cdot z_2 = \emptyset$.
Consider then  the $\,G$-invariant function $\,f(z)= |z_1|^2 + |z_2|^2 - |z_3|^2 -1\,$
defined on   $\,M^\C$. The function  $f$ 
 vanishes on the real hypersurface $\,G \cdot (z_2 \cup_{j=1}^4 w_j)\,$, is negative
on the sets
$\,G\cdot \ell_j(I_j)\,$, for $\,j=1,3\,$ and positive on the sets $\,G\cdot \ell_j(I_j)\,$, for $\,j=2,4$.
Moreover,  for every fixed $\, t _0\in(0,1)\,$, one sees that
$$ f( \exp (isC) \cdot \ell_1(t_0)) =
(\sinh^2 2s + \cosh^2 2s)\sin^2 \frac{\pi}{2}(1 - t_0) -\cos^2 \frac{\pi}{2}(1 - t_0)-1\, $$ 

\nmedskip
is strictly increasing for $\,|s| \to \infty$.  Thus it vanishes exactly 
twice.  
As a consequence, 
 the path $ \exp (isC) \cdot \ell_1(t_0)$ crosses the 
hypersurface $f^{-1}(0)\setminus\{G\cdot z_2\}$ exactly twice, namely on the orbits 
$\,G\cdot  w_1$ and $G\cdot   w_2 $. It follows that
$\,\exp (isC) \cdot \ell_1(t_0) \in
\Omega\,$, for every $s\in\R$. Thus the $T^\C$-orbit through  $\ell_1(t_0)$ is entirely
contained in $\Omega$, as stated.
An analogous argument proves the statement for the higher  
dimensional hyperquadrics. By $(ii)$ of Remark \ref{QUOTOP}, this settles the case when
$\g$ has a reduced restricted root system.

Consider now the case when the restricted root system of $\g$ is non-reduced.
We prove the statement by reducing to the two-dimensional case. 
Set $\, \widehat \g :=\mathfrak s \mathfrak o(2,1)\,$ and 
fix a basis in $\, \widehat \g\,$   of the form $\{\widehat X ,~ \theta (\widehat X) ,~
\widehat A=[\theta (\widehat X),\widehat X ]\} $, where $\,\widehat X\,$
is a root vector in 
$\,\widehat\g^\alpha\,$ and $\,\alpha(\widehat A)=\frac{\pi}{2}$.  Define
$\widehat C= 
\widehat X + \theta (\widehat X)$.
Choose a root vector $X\in\g^{2\alpha}$   
and   
normalize  the  triple  $\{X, \theta(X),  A=[\theta (X),X]\}$ in $\g$ so that
$\,\alpha(A)=\frac{\pi}{4}\,$.   Such a  triple  generates a three-dimensional  $\theta$-stable
subalgebra of $\g$ isomorphic to $\,\widehat \g$. In particular, there
exists  an injective  Lie algebra   homomorphism 
$$\varphi_*\colon \widehat\g \to \g $$
mapping  $\,\widehat X,~  \widehat A $ and $\theta (\widehat X)\,$ into $\, X ,~  A 
$ and $ ~  \theta (X) $ respectively. Clearly
$\,\varphi_*\,$ maps  $\,\widehat C = \widehat X+\theta (\widehat X)\,$ into $\,C=X+\theta (X)\,$as
well.
Let $\,\widehat K =SO(2)\,$ be the maximal compact subgroup of $\,
\widehat G := SO_0(2,1)\,$ and $\,\widehat \k\,$ its Lie algebra. Note that
$\,\widehat \k\,$ and $\,\k\,$ are generated by $\,\widehat C\,$ 
and $\, C$, respectively. One can check that the $\,\C$-linear extension
$\,\widehat \g^\C \to \g^\C\,$ of $\,\varphi_*\,$ induces a Lie group morphism
$\, \varphi\colon  \widehat G^\C
\to G^\C   \,$
mapping $\,\widehat K^\C\,$ to
$\,K^\C$. As a consequence one obtains
a holomorphic map  
(denoted by the same symbol)
$$\,  \varphi\colon \widehat G^\C
/\widehat K^\C \to G^\C/K^\C .$$

\REM{
Recall that $\widehat
\g$ is also isomorphic $\s\l(2,\R)$. Denote by
$\varphi$  the  Lie group morphism $\,\varphi: SL(2,\C) \to G^\C\,$, induced by
the
complexified homomorphism
$ \mathfrak s\mathfrak l (2,\C) \to \g^\C $. 
Observe that by construction 
 $\,\varphi(SL(2,\R))\subset  G$ and   the complexification of the maximal compact subgroup of  $SL(2,\R)$ is mapped into $K^\C$. In particular the center  $Z(SL(2,\C))=\{\pm I_2\} $ is mapped   into 
$ K^\C\,$ as well. Since 
$SO(3,\C) \cong SL(2,\C)/\{\pm I_2\}$, the morphism 
$\varphi$ induces a holomorphic map (denoted by the same symbol)
$$\,  \varphi\colon SO(3,\C)/SO(2,\C)\to G^\C/K^\C    .$$ 
Denote  $SO_0(2,1)$ by $\widehat G$ and  $SO(3,\C)/SO(2,\C)$ by $\widehat G^\C/\widehat K^\C$. }

\noindent
Let $\, \widehat \Omega$ be  the domain  
$$ \widehat \Omega= \widehat G \cdot (\,\widehat \ell_{1}(I_{1}) \cup \widehat w_{1} \cup 
 \widehat \ell_{2}(I_2) \cup \widehat w_4\cup \widehat \ell_{4}(I_{4})\,)   $$ 
in  $\widehat G^\C/\widehat K^\C$ considered in (i). 
 We  claim  that  $ \varphi(\widehat \Omega)\subset \Omega .$
The  map $ \varphi $  is ``equivariant" with respect to the action  of $\widehat G$, that is  $\varphi(g\cdot x)=\varphi(g)\cdot \varphi(x)$, for every $g\in \widehat G$ and $x\in  \widehat G^\C/ \widehat K^\C$.
By the definition of $\varphi_*$  one has 
$\varphi(\exp (it\widehat A) )=\exp (itA) $ and $\varphi(\exp (it\widehat C))=\exp (itC).$ It follows that 
$$\varphi (\widehat  \ell_1(I_1))= \ell_1(I_1),\quad \varphi(\widehat z_2)=z_2,\quad \varphi(\widehat {\mathcal
C})=\mathcal C.$$
 We conclude the proof of the claim by showing
that $\, \varphi(\widehat w_1) \in G
\cdot w_1\,$ and
$\, \varphi(\widehat w_4) \in G \cdot w_4$ (possibly the orbit 
$\,G \cdot w_4\,$ and $\,G\cdot w_1\,$ coincide). 
 By (\ref{LUNA})There is a commutative diagram
$$\begin{matrix} 
  \widehat G \times_{\widehat G_{\widehat z_2}} \widehat V_2  \ & \buildover{\widetilde \varphi}{\longrightarrow}  &  
G\times_{G_{ z_2}}V_2 \cr
                                 &                          &                        \cr
      \downarrow \  \ \        &            &    \downarrow \  \ \cr
                                 &                         &                         \cr
\widehat G^\C/\widehat K^\C \ \ & \buildover{ \varphi}{\longrightarrow}   &  G^\C/K^\C \,
.
\cr 
\end{matrix}\,  $$

\nmedskip
The vertical arrows correspond to the equivariant embeddings given in 
(\ref{LUNA})
and the map
$\,\widetilde \varphi\,$ is defined by $\,[\widehat g, \widehat X] \to
[\varphi(\widehat g), \varphi_*(\widehat X)]$.
 Since  $\,\varphi_*\,$ is an injective homomorphism,   
 $\, \varphi(\widehat w_1) \,$ does not lie on the singular
orbit $G\cdot z_2$. Indeed in the twisted bundle $\,G\times_{G_{ z_2}}V_2\,$
such orbit corresponds to the set $\,\{\,[g,0] \ :\ g \in G\,\}$.
On the other hand 
 $\, \varphi(\widehat w_1) \in  \overline {G \cdot \ell_1(I_1)} \cap
 \overline {G \cdot \ell_2(I_2)}$. Therefore the image
$\, \varphi(\widehat w_1)\,$ necessarily lies on the orbit
$\, G
\cdot w_1$.
 Similarly one proves that $\, \varphi(\widehat w_4) \in G \cdot w_4$.
 In conclusion   $ \widehat  \Omega$ is mapped by
$\varphi$ into $\Omega$, as claimed.
 
Observe that  $\,\,\exp   \c^\C \cdot   \ell_1(I_1) =\varphi
(\exp \widehat
\c^\C \cdot \widehat
\ell_1(I_1))\,$ and recall that in the 2-dimensional case
we already
showed that
$\exp
\widehat
\c^\C
\cdot
\widehat
\ell_1(I_1)\subset
\widehat \Omega$.
Then, by  the above claim, for every $\,z \in \ell_1(I_1)\,$ one has
    $$\,\,T^\C \cdot 
\ell_1(z) \subset \Omega\,,$$ 

\nsmallskip
as required.
 The statement regarding the domain $\,\Omega'\,$
follows from similar arguments.
\end{proof}


\nbigskip
\begin{lem}
\label{ADDITION} Assume that $G$ is embedded in its 
universal
complexification~$G^\C$ and is different from the groups $SL(2,\R)$
and  $Spin(3,1)$.  Let $\,q:\Sigma \to G^\C/K^\C\,$ be a Stein,
$\,G$-equivariant Riemann domain.

\sn
\item{(i)} Let $  \widetilde \ell_1\colon I_1\to \Sigma$ be a lifting
of the slice
$\ell_1$. Assume that the closure of $G\cdot   \widetilde \ell_1(I_1)$
in $\Sigma$ contains points $\eta_1$ and $\eta_4$ mapped by $q$ into the
non-closed orbits $G\cdot w_1$ and $G\cdot w_4$, respectively (possibly 
the orbits $\,G \cdot w_1\,$ and $\,G \cdot w_4\,$ coincide). Then the
singular orbit
$G\cdot z_2$ is contained in
$q(\Sigma)$.

\sn
\item{(ii)} Let $  \widetilde \ell_3\colon I_3\to \Sigma$ be a lifting of the
slice $\ell_3$. Assume that the closure of $G\cdot   \widetilde
\ell_3(I_3)$ in $\Sigma$ contains points $\eta_2$ and $\eta_3$ mapped by $q$
into the non-closed orbits $G\cdot w_2$ and $G\cdot w_3$, respectively
(possibly 
the orbits $\,G \cdot w_2\,$ and $\,G \cdot w_3\,$ coincide). Then the
singular orbit $G\cdot z_2$ is contained in $q(\Sigma)$.
\end{lem}

\smallskip
\begin{proof} (i)
We begin by showing that 
$\Sigma$ contains an open $G$-invariant set which is biholomorphic   to the
domain $\Omega= G \cdot (\,\ell_{1}(I_{1}) \cup w_{1} \cup w_4\cup \mathcal
C^*)$ of Lemma \ref{GLOBALORB}.
By the $\,G$-equivariance of $\,q$, the description of the slice representation at
$\,z_2\,$ given in Remark \ref {ADJOINT}, 
and Proposition \ref{ORBITINJ}, there exists  a $\,G$-invariant neighbourhood
$V$ of
$\eta_1$ in
$\Sigma$  on which
$q$ is injective. Its image 
$q(V)$ intersects the slice
$\ell_2$ in $\ell_2((0,\epsilon))$, for some
$\epsilon>0$.  By $(i)$ of Lemma  \ref{EXTENSION},  
 the map
$s\mapsto (q|_V)^{-1} (\ell_2(s))$, with $s\in (0,\epsilon)$,
extends to a lifting $ 
\widetilde \ell_2\colon  I_2\to
\Sigma$ of $\ell_2$. 
A similar argument yields a lifting $  \widetilde \ell_4\colon  I_4\to
\Sigma$ of $\ell_4$. 
 Since $q$ is injective on the set
$  \widetilde \ell_1(I_1)\cup \eta_1 \cup   \widetilde
\ell_2( I_2) \cup\eta_4\cup    \widetilde \ell_4(I_4) $, 
as well as on every
$G$-orbit (cf.~Proposition \ref{ORBITINJ}),
it is injective on the
$G$-invariant subdomain of $\Sigma$ given by  
$$W:= G\cdot (  \widetilde \ell_1(I_1)\cup \eta_1 \cup   \widetilde
\ell_2( I_2) \cup\eta_4\cup    \widetilde \ell_4(I_4)).$$

\nsmallskip
Note that $q(W)=\Omega$. In particular
$W$  is  biholomorphic to $\Omega$, as
claimed.

Let $\,\mathcal C= \exp i\c \cdot z_2\,$ be the standard Cartan subset
in $\,G^\C/K^\C\,$ starting at $\,z_2$. Recall that $\,T:=\exp \c\,$ 
is a compact torus in $\,G$. 
By Heinzner's globalization theorem (\cite{He1}, Sect. 6.6),
the space  $\Sigma$
can be 
 embedded in its $T^\C$-globalization $\Sigma^*$,  as a $T$-invariant,  orbit-convex
 domain. 
By definition, this means that  the intersection of $\Sigma$ with an $\exp i\c $-orbit
in $\,\Sigma^*\,$ is necessarily connected.

Consider now  the induced local $\,T^\C$-orbit of a point $\,\zeta
\in   \widetilde 
\ell_1(I_1)$ in $\Sigma\,$.  Since $q|_W$ is biholomorphic and $G$-equivariant,
by Lemma \ref{GLOBALORB},  such orbit  is in fact global.
Let $\,C\,$ be a
generator of the abelian subalgebra $\,\c$.
 For  every fixed $\,s >0\,$,  one has
$$\lim_{n \to \infty} \exp (isC) \cdot   \widetilde \ell_1(1/n) =   \widetilde \ell_{2}(s)
\in W\,$$
and
$$\lim_{n \to \infty} \exp (-isC) \cdot   \widetilde \ell_1(1/n) =   \widetilde
\ell_{4}(s)\in W.$$

\nmedskip
By the orbit-convexity of $\,\Sigma\,$ in its  $T^\C$-globalization,
the sequence $\{\,  \widetilde 
\ell_1(1/n)\,\}_n$ converges   to a point $\,\zeta_2\in\Sigma$. By the
continuity of $q$, one has $\,q(\zeta_2)=z_2$. Therefore $\, z_2 \in
q(\Sigma)$, as required. 

\sn
Part (ii) is proved by showing that $\Sigma$ contains an open subset
biholomorphic  to the domain $\Omega'$ of Lemma \ref{GLOBALORB} and arguing
as in the previous case. 
\end{proof}

 \bn 
 
 Let $\,G\,$ be a connected Lie group
and  let $\,  \widetilde G \to G =  \widetilde G/\Gamma\,$ be a covering
of $\,G$.  If $\,X\,$ is a $\,G$-manifold, it can be regarded
as a $\,  \widetilde G$-manifold by letting $\,\Gamma\,$ act trivially on it.

 
\bn
\begin{lem}
\label{KILLtheCENTER} Let $\,G\,$ be a  connected, real  Lie group
and  let $\,  \widetilde G \to G =  \widetilde G/\Gamma\,$ be a finite covering of $\,G$. 
Let $\,X\,$ be a complex $G$-manifold with the property that every 
Stein, $\,G$-equivariant Riemann domain over $\,X\,$ is univalent.
Let $\,q:\Sigma \to X\,$ be a Stein, $\, \widetilde G$-equivariant Riemann domain.
Then 

\begin{itemize}

\smallskip

\item[$(i)$] the image $\,q(\Sigma)\,$ is biholomorphic to the 
quotient $\,\Sigma/\Gamma\,$ and 
$\,q:\Sigma \to q(\Sigma) \,$ can be identified with the quotient map,

\smallskip

\item[$(ii)$]  $\,q\,$ is a $\,  \widetilde G$-equivariant covering.

\end{itemize}

\smallskip
\noindent
In particular $\,q(\Sigma)\,$ is Stein.
\end{lem}

\medskip
\begin{proof}
$(i)$
Since  $\,\Gamma\,$ is a finite subgroup of $\,\widetilde G$, the quotient
$\,\Sigma/\Gamma\,$ can be regarded as the categorical quotient
of $\,\Sigma\,$ with respect to $\,\Gamma$.
Then $\,\Sigma/\Gamma\,$ is a Stein space and  the
quotient map $\,\pi:\Sigma \to \Sigma/\Gamma\,$
is holomorphic (cf.~Thm.~\ref{CATEGO}).
Moreover, since $\,q\,$ is $\,\Gamma$-invariant,
there exists a  
$\,G$-equivariant holomorphic map $\,\widehat q:\Sigma/\Gamma \to X\,$
making the diagram
$$\def\mapright#1{\smash{
\mathop{\rightarrow}\limits^{#1}}}
\def\mapdown#1{
\rlap{$\vcenter{\hbox{$\scriptstyle#1$}}$}{\ \ \big\downarrow \ \, \ }}
\begin{matrix}\Sigma &\mapright{\pi}& \Sigma/\Gamma\cr 
\cr
q \mapdown \  &\swarrow \rlap
{$\vcenter{\hbox{$\widehat q$}}$}\cr
\cr
X\cr\end{matrix}\ \,$$

\nsmallskip
commute.
Since  $\,q=\widehat q \circ \pi\,$
is locally  biholomorphic,  then  $\,\pi\,$ is also locally
 biholomorphic.
In particular $\,\Sigma/\Gamma\,$ is a manifold
 and $\,\widehat q:\Sigma/ \Gamma \to X\,$ is a Stein, $\, G$-equivariant Riemann domain.
By the assumption on $\,X$, the map
$\,\widehat q\,$ is injective, implying  $\,(i)$.

 $(ii)$ Without loss of generality one may assume that $\,\Gamma\,$ acts effectively on $\,\Sigma$. Then
the statement  follows by showing that $\,\Gamma\,$
 acts freely on $\,\Sigma$.
 Assume by contradiction that this is not the case. Then there exists
 $\,\gamma \in \Gamma\,$ whose fixed point set
 $\,Fix(\gamma):=\{\,\zeta \in \Sigma \ : \ \gamma \cdot 
 \zeta = \zeta \}\,$ is not empty.
 Since  $\,Fix(\gamma)\,$ is a proper analytic subset of 
 $\,\Sigma$, it has no interior point. In particular
 there exist $\,\zeta \in Fix(\gamma) \,$
 and a sequence $\,\{\zeta_n\}_n\,$ in the complement of
 $\,Fix(\gamma)\,$ in $\,\Sigma\,$ such that $\ \zeta_n \to \zeta$.
Note that  by the continuity of $\,\gamma$, one has
$\ \gamma \cdot \zeta_n \to \gamma \cdot \zeta =\zeta$.

 Let $\,U\, $ be an open neighborhood of $\,\zeta\,$ on which $\,\pi\,$
 is injective. Then, for $\,n\,$ large enough, both  $\ \zeta_n\,$ and
  $\ \gamma \cdot \zeta_n\,$  lie in $\,U$. 
 Since $\,\Gamma\,$ acts trivially on $\,\Sigma/\Gamma$, it follows that
$\,\pi(\zeta_n)=\gamma \cdot \pi(\zeta_n)= \pi( \gamma \cdot \zeta_n)$.
 On the other hand  since $\,\zeta_n \not \in Fix(\gamma)$, 
one has $\,\gamma \cdot \zeta_n \not= \zeta_n$. This gives
 a contradiction and concludes the proof of the lemma.
\end{proof}

\bigskip
Recall  the following 
 consequence of the uniqueness of path-liftings 
on Riemann domains, which will be often used in the proof of the 
main theorem of this section.


\bn
\begin{lem}
\label{TRIVIAL}
Let  $\,q:\Sigma \to Z\,$ be  a Riemann domain
and let $\,W\,$ be a domain of $\,\Sigma\,$ such that the restriction 
$\,q|_{W}: W \to Z \,$ is bijective. Then $\,W=\Sigma$.
\end{lem}

\bn
Next comes the main result of this section.  


\nbigskip
\begin{theorem}
\label{UNIVALENCE2} 
Let $\,G/K\,$ be a non-compact, rank-one, Riemannian symmetric space.
Assume that $\,G\,$ is a connected, simple, real Lie group which is
embedded in its universal complexification~$\,G^\C\,$ and
is different from $\,SL(2,\R)$.
Then a holomorphically separable, $\,G$-equivariant  Riemann domain $\,q\colon
\Sigma \to G^\C/K^\C\,$ is univalent.
\end{theorem}

\smallskip
\begin{proof} Recall that $\,\Sigma\,$ admits a $\,G$-equivariant
holomorphic embedding into its envelope of holomorphy.
Thus we may assume that
$\,\Sigma \,$ is Stein (cf.~Sect.~2). We prove the
theorem in the case when the $\,G$-orbit diagram of  $\,G^\C/K^\C\,$ is of type 
(\ref{DIAGRAM3}), namely for $\,\g=\s\u(n,1)$.
In all remaining cases, but $\,G=Spin(3,1)$ which is discussed separately,
the statement follows
from the same arguments with fewer steps.

So  we first assume
that $\,G\,$ is different from $\,Spin(3,1)\,$ and divide
the proof 
in three subcases, depending  on  the image
of $\,\Sigma\,$ in $\,G^\C/K^\C$. Finally we discuss the case
$\,G=Spin(3,1)$.


\medskip
(i) {\it The image $q(\Sigma)$ contains the singular orbit $G\cdot z_2$.}
We begin by proving that 
there exists a $G$-invariant domain $V\subset
\Sigma$ with the  properties that
$\,q\,$ is injective on $\,V\,$ and
$$q(V)=G\cdot \Bigl ( \ell_1(1)\bigcup_{j=1}^4\bigl (\ell_j( I_j) \cup w_j \bigr) \cup z_2
\Bigr ) \,.$$
The extended slices  $\,\ell_j:\widehat I_j \to G^\C/K^\C\,$ are defined in (\ref{EXTENDED}).
Let $\,\zeta_2\,$ be an element  in
$\,q^{-1}(z_2)\,$ and let $\,U\,$ be an open
neighborhood of 
$\,\zeta_2$ in $\Sigma$ on which
the restriction   $\,q|_U\,$ is injective. Since
the map
$q$ is open, the image $q(U)$ intersects the slices  for principal orbits
starting at $z_2$ in the sets
$\ell_j((0,\epsilon))$, for $j=1,\ldots,4$ and some $\epsilon >0$. The image
$\,q(U)\,$  also intersects  all non-closed
$G$-orbits containing $G\cdot z_2$ in their closures. 
By Lemma \ref{EXTENSION}  each extended slice $\ell_j$ admits  a lifting
$\,  \widetilde\ell_j\colon \widehat I_j\to \Sigma\,$ such that
$$  \widetilde\ell_j(t)=(q|_U)^{-1}\ell_j(t),\qquad t\in (0,\epsilon).$$
For  $j=1,\ldots,4$, choose points $\eta_j\in (q|_U)^{-1}(G\cdot w_j)$. 
Consider then the open
$G$-invariant set in $\Sigma$ 
$$V :=G\cdot \Bigl ( \widetilde \ell_1(1)\bigcup_{j=1}^4\bigl (\widetilde
\ell_j( I_j) \cup \eta_j \bigr) \cup \zeta_2
\Bigr ) \,.$$
 Since   $q$ is injective on each lifted slice
$\,  \widetilde\ell_j\,$ and on all $G$-orbits 
(cf. Proposition \ref{ORBITINJ}), it  is  injective on $V$ as well.
Hence 
$\,V\,$ is the open $\,G$-invariant domain in $\,\Sigma\,$
with the required properties.

Consider a sequence
 $\{z_n\}$  in $G\cdot \ell_3(J)$, converging to a point on the boundary orbit
$G \cdot w_5$.  Recall that the Levi form of $G\cdot w_5$ is indefinite
(see Rem.~\ref{REMW23}). Then,  
by arguing as in the proof of $(ii)$ in Lemma \ref{EXTENSION}, the domain $V$ can
be enlarged to an invariant  domain $W$ in
$\Sigma$ with the properties
that the restriction $q|_W$ is injective and $q(W)=G^\C/K^\C$.
By Lemma  \ref{TRIVIAL} one has $\,W=\Sigma\,$ and the theorem follows.

\medskip
(ii) {\it The image $q(\Sigma)$ does not contain the orbit $G\cdot z_2$,
but contains a non-closed
$G$-orbit.}
Assume for example that $w_1\in q(\Sigma)$ and let $\eta_1\in q^{-1}(w_1)$.
By the $\,G$-equivariance of $\,q$, the description of the slice representation at
$\,z_2\,$ given in Remark \ref {ADJOINT} 
and Proposition \ref{ORBITINJ}, there exists  a $\,G$-invariant neighbourhood $V$ of
$\eta_1$ in $\Sigma$  on which
$q$ is injective. Its image
$q(V)$ intersects the slices
$\ell_1$ and $\ell_2$ in the sets $\ell_1((0,\epsilon))$ and
$\ell_2((0,\epsilon))$, for some $\epsilon>0$.  Arguing as in previous case,
one can construct a $\,G$-invariant domain $V\subset \Sigma$ with the 
properties that 
$\,q\,$ is injective on $V$ and
$$q(V) =
G\cdot ( \ell_1(\widehat I_1)\cup w_1 \cup  \ell_2(\widehat I_2)).$$

\sn
If  $V=\Sigma$ (this is possible by Theorem \ref{STEIN}), then
the map $q$ is injective, as desired. 
If $V\not=\Sigma$, then there exists a point $\eta$ in the closure
of $V$ in $\Sigma$  which is mapped by $q$ into one of the non-closed orbits
$G\cdot w_2$ or $G\cdot w_4$.  Assume that 
 $\,q(\eta)\,$ lies in $\, G\cdot w_4$.
Then  by $\,(i)\,$ of Lemma \ref{ADDITION} the image $q(\Sigma)$ necessarily
contains $G\cdot z_2$,  contradicting the current assumption. 

If $q(\eta)\in G\cdot w_2$, iterating the procedure of lifting
slices and  orbits, we can enlarge $V$ to an invariant  domain  
$ W $ in $\Sigma $ on which $\,q\,$ is injective and such that 
$$q(W)=G\cdot ( \ell_1(\widehat I_1)\cup w_1 \cup  \ell_2(\widehat I_2)\cup w_2
\cup  \ell_3(\widehat I_3) ).$$
In particular $\,W\,$ is biholomorphic to $\,q(W)$,
which is not Stein by Theorem \ref{STEIN}.
Hence $\,W\,$ is a proper subset of 
 $\,\Sigma\,$ and there exists a point $\eta $ in the closure of
$\,W\,$ in
$\,\Sigma\,$  whose image 
$q(\eta)$ lies either in $G\cdot w_3$ or in
$G\cdot w_4$. In both cases  Lemma
\ref{ADDITION}  implies that $q(\Sigma)$ contains $G\cdot z_2$,
contradicting the current assumption. 
In conclusion, if $q(\Sigma)$ does not contain the singular
orbit  $G\cdot z_2$ but contains
the non-closed orbit $G\cdot w_1$, then $q$ is injective.
For the other non-closed $G$-orbits, the theorem can be proved
by  arguing  in a similar way.

\medskip
(iii)
{\it The image $q(\Sigma)$ contains no non-closed $G$-orbits.} 
This assumption implies that the image 
 $q(\Sigma)$ contains none of the singular orbits lying in
the closure of a non-closed
$G$-orbit. More precisely $q(\Sigma)$ contains neither 
$G\cdot z_2$ nor $G\cdot z_3$.
Note that the hypersurfaces $G\cdot (z_2
\bigcup_{j=1}^4w_j)$ and
$G\cdot (z_3
\cup w_5)$   disconnect
$G^\C/K^\C$. Therefore there exists a slice $\ell=\ell_j$, for some
$j=1,\ldots,5$, such that $q(\Sigma)=
G \cdot \ell(J)$, for some
interval $J\subset \widehat I$ which is open in $\widehat I$.
Define $M:=q^{-1}(\ell(J))$. One has that
$\Sigma=G\cdot M$.
Moreover, since 
 $\,q\,$ is injective on $\,G$-orbits
(see Prop. \ref{ORBITINJ})  and
 every orbit in $\,q(\Sigma)\,$ intersects
 $\,\ell(J)\,$ in a single point,   every
$\,G$-orbit in $\,\Sigma\,$ intersects $\,M$ in a single point as well.
 As a consequence, the 
 surjective map $\,\Pi: \Sigma \mapsto  M$, given by
$\,\zeta \mapsto G \cdot \zeta \cap M$,
 is well defined.  
\medskip
\item[] $Claim.\ $ The map $\,\Pi\,$ is continuous.
\smallskip
\item[] 
$Proof\ of\ the \ Claim.\ $ Let $\,N\,$  be an  open set in $\,M\,$. We prove the claim
by showing  that for every $m\in N$ and $\zeta\in \Pi^{-1}(m)$, there exists an open neighbourhood of $\zeta$ in $\Sigma$ which is contained in $\Pi^{-1}(N)$.  
By construction  $\zeta=g\cdot m,$ for some $\,g \in G$.
Let $\,V\,$ be an open neighborhood of $\,m\,$ in $\, \Sigma\,$ 
on which $\,q\,$ is injective. Choose an open 
interval $\,J' \subset J\,$ such that $\,q(m) \in \ell (J')\ \subset q(V)$.
Note that $q(m)$ either sits on a principal $G$-orbit, or on the singular orbit $G\cdot z_1\cong G/K$.
Let $\,   \widetilde \ell (J')\,$  be the lifting of $\, \ell (J')\,$ via the  restriction
$\, q|_V.\,$ By schrinking $\,J'\,$ if necessary, one can find an
open  neighborhood $\,  U\,$ of the identity  in $G  $ such that $\,  U \cdot    
\ell (J')\,$ is open and 
contained in  $\,q(V)$.  This  fact is clear  if $\,q(m)\,$ lies on a principal $G$-orbit (see diagram (\ref{DIAGRAM3})).  If $\,q(m)\,$ lies on the singular orbit
$\,G \cdot z_1$, it follows from the equivariant embedding (\ref{LUNA}) 
 at $\,z_1\,$ and the compactness of  
the isotropy subgroup $\,G_{z_1} \cong K$. 

As a result,  $\,  U \cdot   \widetilde \ell (J') = (q|_V)^{-1}(  U \cdot
\ell (J'))\,$ is an open neighbourhood of $m$ in  $\Sigma$  and 
  $\,g  U\cdot   \widetilde \ell(J') $ is
an open neighbourhood of
$\zeta$ contained  in $\, \Pi^{-1}(N)\, $. Hence  $\,\Pi^{-1}(N)\,$
is open in $\,\Sigma$, as
wished  (one can show that $\,M \cong \Sigma/G\,$ and that $\,\Pi\,$ can be
identified with  the quotient map).

\sn

By the above claim, $\,M\,$ is connected and  is 
a one-dimensional real-analytic submanifold of $\,\Sigma$. It follows that 
$q$ is injective on  $M$.
Moreover $\,M\,$ and $\,q(M)\,$ are slices for the $\,G$-action in 
$\,\Sigma\,$ and $\,q(\Sigma)$, respectively. 
Since $q$ is injective on
$G$-orbits, it is injective on $\Sigma$ implying the theorem.

\sn
(iv) {\it The group $\,G\,$ is  $\, Spin(3,1)$.}
Assume by contradiction that $\,q:\Sigma \to G^\C/K^\C\,$
is not univalent. Recall that the center of $\,G\,$ acts trivially 
on $\,G^\C/K^\C\,$ and that by $\,(i)$, $\,(ii)\,$ and $\,(iii)\,$ the statement holds true for the group $\,SO_0(3,1)$. Then Lemma \ref{KILLtheCENTER},
applies to show that the restriction of $\,q\,$ to every $\,G$-orbit is
a  double covering and the image $\,q(\Sigma)\,$ is Stein.
On the other hand, by Theorem \ref{STEIN},
all Stein $\,G$-invariant domains in $\,G^\C/K^\C\,$
contain a singular orbit diffeomorphic to $\,G/K$. 
Since $\,G/K$ is simply connected, this gives a contradiction.
This concludes the proof of the theorem.
\end{proof}

\bigskip
When  $G=SL(2 ,\R)$, non-injective, Stein $G$-equivariant Riemann domains  over $G^\C/K^\C$ do exist. Next we construct one such Riemann domain explicitly. It turns out that such an  example is essentially the only possible one.
Indeed by Lemma \ref{KILLtheCENTER}, 
if $\,q : \Sigma \to G^\C/K^\C\,$ is
a Stein, $G$-equivariant Riemann
domain which is not univalent, then  
the center $\,\Gamma=\{\pm I_2\} \,$ acts freely on $\,\Sigma\,$. Moreover, $\,q\,$ is a
$G$-equivariant covering onto its image $\,q(\Sigma)\,$ which turns out to be Stein.
It follows that  the restriction
of $\,q\,$ to every $\,G$-orbit is a double covering. Thus  the singular orbits $\,G\cdot z_1\,$ and $\,G\cdot z_3$,
which are simply connected,
cannot lie in $\,q(\Sigma)$.  
Then,  by  Theorem \ref{STEIN},
the image $q(\Sigma)$  coincides with a domain
$\,S_i(b)$, for some $\,i=1,2\,$ and
$\,b \ge 0$. For every $\,S_i(b)\,$ there is exactly one  
 $G$-equivariant
double covering.
In the example below we carry out its construction for
$\,q(\Sigma)=S_1(0)$.

\bigskip
\begin{exa}
\label{NONINJ1}
{\rm Let $\,G=SL(2,\R)$. 
Consider the Stein domain
$S_1(0) $ in  $G^\C/K^\C$ defined in (\ref{DOMAINS}). 
Let   
$$\ell_2\colon \R^{>0}\rightarrow G^\C/K^\C,\quad 
\ell_2(s):= \exp (isC) z_2\,$$
be the slice map defined in
(\ref{SLICE24}) of Section 4. The  isotropy subgroup in $G$ of
every point  $\,\ell_2(s)\,$  coincides with  $\{\pm I_2\}$
(cf.~Remarks \ref{SU21} and \ref{QUOTOP} ).
It follows that 
$\,S_1(0):= G \cdot \ell_2(\R^{>0}) \,$ is topologically equivalent to
$SO_0(2,1)\times \R^{>0}$. 
Define 
$\,\Sigma:=G\times \R^{>0}$.
Since
$\,G\,$ is a double
covering of $SO_0(2,1)$, the map  
$$q\colon \Sigma\longrightarrow S_1(0), \quad
(g,s )\mapsto g\ell_2(s)$$
defines a double covering of $\,S_1(0)$.
As a consequence, with the complex
structure pulled back from $\,S_1(0)\,$
the manifold
$\Sigma$ is Stein (\cite{St})
and the map $\,q\,$ is a holomorphic covering.
In other words,  $q\colon
\Sigma\longrightarrow S_1(0)$ defines
 a non-univalent Stein,
$G $-equivariant Riemann domain over $G^\C/K^\C$.}
\qed
\end{exa}

\bn
\begin{remark} 
\label{FINITECENTER} 
By the results of Lemma \ref{KILLtheCENTER}, one can show that  
Theorem \ref{UNIVALENCE2} also holds for $\,G\,$ not embedded in $\,G^\C$,
provided that the center $\,\Gamma\,$ of $\,G\,$ is finite
and $\,G\,$ is not a covering of $\,SL(2,\R)$ (\,cf. $\,(iv)\,$ in the proof of Theorem \ref{UNIVALENCE2}\,).
If $\,G\,$ is a covering of $\,SL(2,\R)$, a 
construction similar to the one  in Example \ref{NONINJ1}
yields a non-univalent, Stein
$\,G$-equivariant Riemann domain over $\,G^\C/K^\C$. 
\qed
\end{remark} 
 
 
\bn

\section{ Univalence over $\,G^\C$}

\bn

Let $G$ be a connected, non-compact,  real simple Lie  group,
$K \subset G$ a
maximal compact subgroup and $G^\C$ its universal complexification.  In  this
section we prove a univalence result for $G\times K$-equivariant Riemann
domains over $G^\C$, when the symmetric space $G/K$ has rank one. We also 
discuss some examples.

\nbigskip
\begin{theorem}
\label{UNIVALENCE3} 
Let $\,G/K\,$ be a non-compact, rank-one, Riemannian symmetric space.
Assume that $\,G\,$ is a connected, simple, real Lie group which has 
finite center and is not a covering of $\,SL(2,\R)$.
Then a  holomorphically separable,
$\,G \times K$-equivariant Riemann domain $\,p\colon Y\to\,G^\C\,$  is univalent.
\end{theorem}

\begin{proof}
Recall that $\,Y\,$ admits a $\,G \times K$-equivariant holomorphic
embedding into its envelope
of holomorphy. Thus we may assume that $\,Y \,$ is Stein (cf. Sect. 2). 
Consider the induced Stein, $\,G$-equivariant 
Riemann domain $\,q: Y \qq K \to G^\C/K^\C\,$ constructed in Section 3. 
By Remark \ref{FINITECENTER} the map $\,q\,$ is injective. Then, by 
Corollary  \ref{INI} the Riemann domain 
$\,p\colon Y\to\,G^\C\,$ is univalent, as wished.
\end{proof}

\bigskip
When  $\,G\,$ is either $\,SL(2 ,\R)\,$  or   a non-trivial covering of $\,SL(2 ,\R)\,$, 
a construction similar  to the one 
in  Example \ref{NONINJ1} yields
examples of non-univalent,
Stein, $G \times K$-equivariant Riemann
domains over $\,G^\C$.

 \bigskip
\begin{exa}
\label{NONINJ2} 
{\rm 
Let $\,G=SL(2,\R)\,$ and let
$\,S_1(0) \,$ be the Stein, $\,G$-invariant 
domain in  $G^\C/K^\C$ defined in  (\ref{DOMAINS}).
As we observed in Example \ref{NONINJ1},  the domain $\,S_1(0)\,$ is diffeomorphic to
$\,SO_0(2,1)\times \R^{>0}$. 
Define $\Omega:= \pi^{-1}(S_1(0))$, where $\,\pi\colon G^\C\rightarrow G^\C/K^\C\,$ 
is the canonical projection. 
Since $\,\pi\,$ is holomorphic and both $\,S_1(0)\,$ and  $\,G^\C\,$
are Stein, the domain  $\,\Omega\,$ is Stein as well.
Consider  the 
slice $\, \ell_2\colon \R^{>0}\to G^\C/K^\C \,$ (cf. (\ref{SLICE24}) of Sect.~4)
and its lifting  to $\,G^\C\,$ defined by $\,\widetilde \ell_2(s):=\exp (isC)\exp (iA_2)$.
Note that the map $\,SO_0(2,1) \times \R^{>0}  \times K^\C$, given by 
$\, (g,s,k) \mapsto g \widetilde \ell_2(s)k^{-1} $, is a diffeomorphism. 
Define $\,Y:= G \times \R^{>0}  \times K^\C.$ By construction the map 
$$\,p:Y \to \Omega\,, \quad (g,s,k) \mapsto g \widetilde \ell_2(s)k^{-1}  \,$$
is a double covering.
With the complex
structure pulled back from $\,\Omega$,
the manifold
$\,Y\,$ is Stein (\cite{St})
and the map $\,p\,$ is holomorphic.
 Let  $\, G \times  K$
act  on
$\,Y\,$ by 
$\,(l,h) \cdot (g,s,k) :=(lg,s,hk)\,$
and  on $\,\Omega \,$ by  left and right translations.
Then $\,p\,$ defines
 a non-univalent,
Stein, $\,G \times K$-equivariant Riemann domain over $G^\C$.
\qed
}
\end{exa}

\bigskip
Let $\,G=K \times N \,$ be the product of a compact 
Lie group and a simply connected nilpotent Lie group. 
Then a holomorphically separable, 
$\,G$-equivariant Riemann  domain over $\,G^\C  \,$
is necessarily univalent (see \cite{CL}, \cite{Ia}, \cite{CIT}).
The above example  shows that an
analogous statement does not hold for a semisimple Lie group $\,G$.
 Next  we exhibit  a different counter-example for $\,G=SO_0(2,1)$, a group
 which meets the assumptions of Theorem \ref{UNIVALENCE3}. 
Such  example
was pointed out to us by K. Oelijeklaus.
We are not aware of similar constructions in higher dimension.
That is, if dimension of $\,G/K\,$ is greater than two, univalence of 
holomorphically separable, 
$\,G$-equivariant Riemann  domains over $\,G^\C  \,$ 
seems to be an open question.

 \bigskip
\begin{exa}
\label{KLAUS} \rm
Let $G=SO_0(2,1)$. Then $G^\C=SO(2,1,\C)$ and  $K^\C= SO(2,\C)$. 
Let $S_1(0)$ be the $G$-invariant Stein  domain in $G^\C/K^\C$ defined in (\ref{DOMAINS}) and let  $\Omega=\pi^{-1}(S_1(0))$, where 
$\,\pi\colon G^\C\to G^\C/K^\C\,$ is  the canonical projection.
As we already observed in 
 Example \ref{NONINJ2}, the domain $\Omega$ is  a Stein,
 $G$-invariant domain in $G^\C$ which is diffeomorphic to
 $\,G\times \R^{>0}\times   K^\C\,$. 
Denote by $\widetilde K^\C$ the universal covering of $K^\C$ and by
$\,\psi\colon   \widetilde K^\C\rightarrow K^\C$ the covering homomorphism.
Define $\,Y:=G\times \R^{>0}\times   \widetilde K^\C\,$ and let 
 $\,G\,$ act on $\,Y\,$ by left translations. 
Consider  the 
slice $\, \ell_2\colon \R^{>0}\to G^\C/K^\C \,$ (cf. (\ref{SLICE24}) of Sect.~4)
and its lifting  to $\,G^\C\,$ given  by $\,\widetilde \ell_2(s):=\exp (isC)\exp (iA_2)$.
Define a  $\,G$-equivariant covering of $\,\Omega\,$ by
$$p\colon Y\longrightarrow \Omega,\quad \quad 
(g,s,k)\mapsto g\widetilde \ell_2(s) \psi(k^{-1}).$$

\sn
With the complex
structure pulled back from $\,\Omega$,
the manifold
$\,Y\,$ is Stein (\cite{St})
and the map $\,p\,$ is holomorphic. In particular  
$\,p\colon Y\longrightarrow \Omega\,$ defines a non-univalent,  Stein, $G
$-equivariant Riemann domain over $G^\C$.
\qed
\end{exa}

\bigskip
\begin{rem} 
One can show that $\,\Omega\,$ is a holomorphically trivial $\,\C^*$-bundle over $\,S_1(0)$. Thus it is biholomorphic 
to $\,S_1(0) \times \C^*\,$ and consequently $\,Y\,$ is biholomorphic to
$\, S_1(0)  \times \C$.
 After identyfing $\,S_1(0) \,$ with $\,SO_o(2,1) \times \R^{>0}$, one sees
that  the map $\,SO_o(2,1) \times \R^{>0} \to G^\C$, given by $\,(g,s) \mapsto 
g \widetilde \ell_2(s)$,
defines a global $\,C^\infty$-section  of the holomorphic $\,\C^*$-bundle
$\,\pi |_\Omega:\Omega \to
S_1(0) $. Hence  such bundle   is differentiably trivial and, by Oka principle, is also holomorphically trivial  (cf.~\cite{Gr}), as claimed.
For the sake of completeness, we  explicitly   construct  a trivialization 
on the model of $\,G^\C/K^\C\,$ discussed in Example \ref{NONREDUCED1}
and Remark \ref{N=1}. 
 
Let  $\,G=SU(1,1)\,$ and
 identify $G^\C/K^\C$ with  $\,  \P^1 \times \overline {\P}^1 \setminus \{\,\langle z, w\rangle_{1,1} = 0\,\}\,$.
Note that $\,S_1(0) \,$ corresponds to the subset $\,\{\,([1:u],[\bar v:1]) 
 \ : \ u,v \in \Delta, \ u \not = v\,\}\,$ (cf.  Example \ref{NONREDUCED2}). Denote by $\,D\,$  the diagonal in  $\Delta \times \Delta$. 
 Then  the injective holomorphic map $\,\Delta \times \Delta \setminus D \to \P^1\times\overline\P^1\setminus\{\langle z,w\rangle_{1,1}=0\}\,$ defined by
 $$ (u,v)\mapsto( [1:u],[\bar v:1])$$
 
\noindent
identifies $\,\Delta \times \Delta \setminus D\,$ with $\,S_1(0)$.
The map $\Delta \times \Delta \setminus D \to G^\C$  given by 
$$(u,v) \mapsto \left ( \begin{matrix}  1  &  \frac{1}{u-v}  \cr
                                      v  &  \frac{u}{u-v} \cr 
\end{matrix} \right )\,   $$

\nmedskip 
defines a global holomorphic section of the
$\,\C^*$-bundle $\,\pi|_\Omega:\Omega \to  S_1(0)$, since one has 
$$\left ( \begin{matrix}  1  &  \frac{1}{u-v}  \cr
                                      v  &  \frac{u}{u-v} \cr 
\end{matrix} \right ) \cdot ([0:1],[0:1]) = ( [1:u],[\bar v:1])\,.$$

\mn
As a consequence the map $  (\Delta \times \Delta) \setminus D
 \times \C^* \to 
\Omega$, given by 

$$(u,v,\lambda) \to \left ( \begin{matrix}  1  &  \frac{1}{v-u}  \cr
                                      v  &  \frac{1}{v-u} \cr 
\end{matrix} \right )
\left ( \begin{matrix}  \lambda^{-1}  & 0 \cr
                                      0  &   \lambda  \cr 
\end{matrix} \right )$$

\nmedskip
defines a biholomorphism from $\, S_1(0) \times \C^*\,$
onto $\,\Omega$.
\qed
\end{rem}

\bn

\section{ Appendix. The Levi form of non-closed hypersurface orbits.}

\bn

In this section we  outline the computation of the Levi form of non-closed hypersurface $G$-orbits in $G^\C/K^\C$. The results of this computation were used in Section 6, where we completed the classification of Stein  $G$-invariant domains in $G^\C/K^\C$.
Recall that every real hypersurface $S$  in a complex manifold  inherits a $\,CR$-structure of hypersurface type. Let $J$ denote  the complex structure of the ambient manifold. For every $x\in S$, the tangent space to $S$ at $x$ decomposes as
$$TS_x=T_\C S_x\oplus NS_x,$$
where $T_\C S_x=TS_x\cap J(TS_x)$ is a complex subspace of $TS_x$, called the {\it complex tangent space}, and $NS_x$ is a one-dimensional  {\it real} subspace.
Denote by $TS=T_\C S\oplus NS $ the tangent bundle of $S$.  The subbundle $(T_\C S)^\C\subset TS^\C$ of the complexified tangent bundle $TS^\C$ decomposes as $HS\oplus AS$, where $HS$ and $AS$ denote its $(1,0)$ and $(0,1)$ components, respectively.  Let  ${  Z} $ be a tangent vector in $T_\C S_x$ and  $\widehat { Z}$ an arbitrary extension of ${ Z} $ to a local section of $T_\C S$. Then the vector fields $\half(\widehat { Z} -i J \widehat { Z})$ and $\half(\widehat { Z} +i J \widehat { Z})$ define local sections of the bundles $HS$ and $AS$, respectively.
The Levi form of $S$ at $z$ is the hermitian form 
 $L_x\colon T_\C S_x\times T_\C S_x \to (TS_x)^\C/(T_\C S_x)^\C$ defined~by
 $$L_x( Z,W, ):= {i\over 4}[ \widehat Z-i J\widehat Z, \widehat W+i J\widehat W]_x\qquad \hbox{mod $(T_\C S)^\C$}. $$
 In the hypersurface case, $(TS_x)^\C/(T_\C S_x)^\C$ is a one-dimensional complex vector space. When $Z$ varies in $T_\C S_x$, the image  of the quadratic form $ L_x( Z,Z, )$ is contained in its real part, which can be identified with $NS_x\cong \R$. 
 We say that the Levi form of $S$ is {\it definite} if  $\{L_x(Z,Z)\}$ is a halfline in $NS_x$; that it is {\it indefinite} if  $\{L_x(Z,Z)\}$ coincides with $NS_x$; that it is {\it identically zero} if $\{L_x(Z,Z)\}=\{0\}$  (for more details, see \cite{Bo}).

\bn
\subsection{ Non-closed orbits with a totally real singular orbit in their closure.} 
We first consider non-closed $G$-orbits which contain in their closure the orbit of a point $z=\exp iA K^\C\in\mathcal A_0$,  satisfying the condition $\alpha(A)=\pi/2$, with  $\alpha$ a simple restricted root (see (\ref{BASEPOINT1}) and (\ref{BASEPOINT2}) in Sect.~4).
The singular orbit $G\cdot z$ is diffeomorphic to  a rank-one, pseudo-Riemannian symmetric space $G/H$, embedded in $G^\C/K^\C$ as a totally real submanifold  of maximal dimension. Let $(\g=\h\oplus \q,\tau_z)$ be the corresponding symmetric algebra.  Non-closed $G$-orbits in $G^\C/K^\C$  containing  $G\cdot z$ in their closure  are in one-to one correspondence with the nilpotent $Ad_H$-orbits in~$\q$ (cf.~(\ref{LUNA}) and Remark~\ref{ADJOINT}).

 Let $X  $ be an  element in $\q$ and  let $x =\exp iX  \cdot z $ be  the corresponding point in $G^\C/K^\C$. Denote by  $S $ the $G$-orbit of $x $.
Denote by $\pi\colon G^\C\to G^\C/K^\C$ the canonical projection and by $\pi_*$ its differential. Then the tangent space to $S$ at $x $ is generated by the vector fields induced on $G^\C/K^\C$ by the one-parameter subgroups in $G$, via the map
\begin{equation} 
\label{VECTORFIELDS}
{}^*\colon \g\to  T(G^\C/K^\C)_{x },\quad X\mapsto X^*:=(\pi_*)_{x }(~{d\over {dt}}_{|_{t=0}}\exp tX~).\end{equation} 
Observe that $T(G^\C/K^\C)_z\cong \q^\C$ and $T(G^\C/K^\C)_{x }\cong Ad_{x }\q^\C.$ 
Hence the vector $X^*$ 
 is  the   $  Ad_{x }\q^\C $-component of $X$ in the decomposition $\g^\C=Ad_{x }\h^\C \oplus Ad_{x }\q^\C $.
 
 \bn 
 To explicitly determine base points for such non-closed orbits  and their tangent spaces,  we decompose $\g$ by an appropriate restricted root system.
Fix a maximal abelian subalgebra $\b\subset \h\cap \p$. Since $\g$ has real rank one, $\dim \b= 1$ and $Z_\g(\b)= \b\oplus Z_\k (\b)$.
Let  $\Delta_\b$ be the  restricted root system of $\g$ with respect to $\b$  and let 
$$\g=\g^0 \oplus \g^{\pm\lambda}\oplus  \g^{\pm 2\lambda},\qquad \g^0=Z_\g(\b)$$
be  the corresponding restricted root decomposition. Every root space $\g^\lambda$ is ~$\tau_z$-stable.  For every $\mu\in\Delta_\b\cup\{0\}$, we 
 indicate by  $\g^\mu_\h$ and $ \g^\mu_\q$ the  intersections of $\g^\mu$ with $\h$ and $\q$, respectively. In particular, we have  a combined decomposition
$$\g=\h\oplus \q,  \qquad\hbox{where}\quad   \h=\g^0_{\h\cap \k }\oplus \g_\h^{\pm \lambda}\oplus \g_\h^{\pm 2\lambda}\oplus \b,\quad\hbox{and}\quad  \q=\g^0_{\q\cap \k} \oplus \g_\q^{\pm \lambda} \oplus \g_\q^{\pm 2\lambda} .$$
Here $\g^0_{\h\cap \k} $ and $\g^0_{\q\cap \k}$ denote the intersections  of $Z_\k(\b)$ with $\h$ and $\q$ respectively. Note that, by the real rank one condition,  $\g^0_{\q\cap \k}$ coincides with $\g^0_{\q}$.
If the restricted roots system $\Delta_\b$ is reduced, then    $\g^{\pm 2\lambda}=\{0\}$.
\bn
\begin{lem} 
\label{RED}  Let $\g$ be a simple real  Lie algebra of real rank one, with reduced restricted root system (i.e. $\g=\s\o(n,1)$).  Then  the following facts hold:

 \sn
{(i)} $\dim \g_\q^{\pm \lambda} =1$; 

\sn
{(ii)} $[\g_\q^{ \lambda},\g_\h^{ -\lambda}] = \g_\q^{ 0}, 
\quad   [\g_\q^{ -\lambda},\g_\h^{  \lambda}] = \g_\q^{ 0}.$ 

\end{lem}  

\mn
\begin{proof}   Observe that $\theta \g_\q^{ \lambda}=\g_\q^{-\lambda}$. Hence $\g_\q[\lambda]:=\g^\lambda_\q\oplus  \g^{-\lambda}_\q$ is a $\theta$-stable subspace  of $\q$ and   $\dim \g_\q[\lambda]\cap \p =\dim \g^\lambda_\q$. Since $\g^0_\q\subset \k$  and $\dim \p\cap\q=1$ (see the proof of Lemma \ref{LOCAL1} (ii)),   statement (i) holds.  Statement (ii) can be verified directly. 

\end{proof}
 
\bn
\begin{lem} \label{NONRED} Let $\g$ be a real simple Lie algebra of real rank one, with non-reduced restricted root system (i.e. $\g=\s\u(n,1)$,  $\s\p(n,1)$, or $\f_{4}^*$). 
Then the following facts hold:

\sn
{(i)}  The root spaces $\g^{\pm2\lambda} $ are contained in $\h$.  Therefore $\g^{\pm2\lambda}_\q=\{0\}$;

\sn
{(ii)} $\dim \g^{ \pm \lambda}_\q>1$; 

\sn
(iii) Fix $X^0_\lambda\in \g^{ \lambda}_\q$ and denote by $(\g^{ \lambda}_\q)_0$ a complement of  $\R X^0_\lambda$ in $\g^{ \lambda}_\q$ (resp. by $(\g^{-\lambda}_\q)_0$ a complement of  $\R \theta X^0_\lambda$ in $\g^{- \lambda}_\q$). Then $$[X^0_\lambda,\g^{0}_{\h\cap \k}]=(\g^{  \lambda}_\q)_{0},\quad [X^0_\lambda,\g^{-\lambda}_\h]=\g^0_\q,\quad   [X^0_\lambda,\g^{-2\lambda}_\h]=(\g^{ -\lambda}_\q)_{0} .$$
\end{lem}

\mn
\begin{proof}  Real rank one Lie algebras with a non-reduced restricted root system are  equal-rank. Hence the root system $\Delta$ of $\g^\C$, with respect to a maximally split Cartan subalgebra of $\g$ extending $\b$, has a real root. Since $\dim \g^{2\lambda}$ is  odd,  the restriction of such a root  to $\b$ coincides with the  restricted root $2\lambda$ (cf. \cite{Hl}, p. 584).
Further, by  Remark 2.13 in  \cite{Ge1},  the subalgebra 
$\h$ is a non-compact real form of $Ad_z\k^\C\cong \k^\C$, with respect to the conjugation $\sigma\tau_z|_{Ad_z\k^\C}$. Precisely, if   $\g=\s\u(n,1)$,  $\s\p(n,1)$, or  $\f_{4}^*$, then $\h$ is given by $\u(n-1,1)\oplus \u(1)$, $\s\p(n-1,1)\oplus \s\p(1)$ and $\s\o(8,1)$, respectively.
Since  $\h$ is equal-rank, the root spaces $\g^{\pm 2\lambda}$ have non-trivial intersection with $\h$.
Statements (i) and (ii) follow then by looking at the dimensions of the restricted root spaces  of $\h$ and $\g$ (cf.~\cite{Hl}, p.~532).
Statement (iii) can be verified directly.

\end{proof}

\bn 
{\bf Reference points for non-closed $G$-orbits.}  
Let $\mathcal C=\exp i\c\cdot z$ be the standard Cartan subset in $G^\C/K^\C$ with base point $z$.  
Recall that $\c=\R(X +\theta (X) )$, where $X $ is a non-zero vector in $\g^\alpha$ (here $\g^\alpha$ is a restricted root space with respect to the adjoint action of $\a\subset \p$, as in Sect.~4). Normalize the triple $\{X ,\theta (X) , A:=[\theta (X) ,X ]\}$ so that $\alpha(A)=2$.  Define $B:=X -\theta (X)$ and  $\b:=\R(X -\theta (X) )$. One easily verifies that $\b$ is a maximal abelian subalgebra in $\h\cap \p$. 
If the restricted root system $\Delta_\b$ is reduced, then
\begin{equation}
\label{BASEPOINTS}
X_\lambda^0={1\over 2}(A-(X +\theta X ) ), \qquad X_{-\lambda}^0 ={1\over 2}(A+(X +\theta X ))
\end{equation}
are generators of the one-dimensional spaces $\g_\q^{ \lambda} $ and $\g_\q^{- \lambda} $,  respectively.  They  satisfy the relations 
$$[B,X^0_\lambda]=2X^0_\lambda,\quad [B,X^0_{-\lambda}]=-2X^0_{-\lambda},\quad [X^0_\lambda,X^0_{-\lambda}]=B,\qquad \theta X^0_\lambda=- X^0_{-\lambda}.$$
The vectors $$  X_\lambda^0,\quad   X_{-\lambda}^0, \quad   -X_{ \lambda}^0,\quad -X_{-\lambda^0} $$ are a complete set of representatives of the nilpotent  $Ad_H$-orbits in $\q$. The corresponding points in $G^\C/K^\C$ 
$$x_0=\exp iX^0_\lambda \cdot z,\quad x_1=\exp iX^0_{-\lambda}\cdot z,\qquad   y_0=\exp( -iX^0_\lambda) \cdot z,\quad y_1=\exp (-iX^0_{-\lambda})\cdot z   $$ lie on  non-closed $G$-orbits containing the singular orbit $G\cdot z$ in their closures. In the orbit  diagram  (\ref{DIAGRAM1}), the $G$-orbits of 
$x_0$,  $ x_1$,  $ y_0$, $y_1$  are represented by  $w_3$, $w_2$, $w_1$,  $w_4$,  respectively. If $\dim G/K>2$ the points  $x_0$ and $ x_1$ lie on  the same $G$-orbit and likewise the points $y_0$,   and $ y_1$ (cf. diagram  (\ref{DIAGRAM2})). 
When the restricted root system $\Delta_\b$ is non-reduced, all points $x =\exp iX_\lambda\cdot z$, with $X_\lambda\in \g^\lambda_\q$,  and  $y=\exp iX_{-\lambda}\cdot z$, with $X_{-\lambda}\in \g^{-\lambda}_\q$, 
 lie on the same $G$-orbit. They are represented by $w_5$ in the orbit diagrams (\ref{DIAGRAM3}) and~(\ref{DIAGRAM4}).

\mn
\begin{remark}
\label{REM1}
When the restricted root system $\Delta_\b$ is reduced,  the points
 $x_0$ and $ x_1 $ lie on the boundary of the Stein domain $D_2(0)$ (cf.~Theorem \ref{STEIN}). 
The points
 $y_0$ and $ y_1 $ 
  lie on the boundary of the Stein domain $D_1(0)$ 
(cf.~Theorem \ref{STEIN}). \qed
 \end{remark}


\bn 
{\bf The tangent space  to the $G$-orbit of $x_0$.}  
Denote by  $S$ the $G$-orbit of  the  point $x_0=\exp iX^0_{ \lambda}\cdot z$, with $X^0_{ \lambda}\in \g^\lambda_\q$. 
In the next Lemma we determine the generators of the tangent space to $S$ at $x_0$, namely the vectors $X^*\in TS_{x_0}$, for~$X\in\g$.

\sn

\begin{lem} 
\label{FIELDS}  
\item{(i)}  Let $Y_{2\lambda}\in\g^{2\lambda}_\h$ and 
$Y_{\lambda}\in\g^{\lambda}_\h$. Then $Y^*_{2\lambda}=Y^*_\lambda= 0$.
 \mn
\item{(ii)}   Let $X_{\lambda}\in\g^{\lambda}_\q$. Then $X^*_\lambda= Ad_{x_0}X_\lambda  $.
\mn
 \item{(iii)}  Let   $B\in\b$. Then $B^*=i\lambda(B) Ad_{x_0}X_{\lambda}^0 $.  
 \mn
\item{(iv)} Let $W_0\in\g^{0}_{\h \cap \k}$. Then $W_0^*=-i Ad_{x_0}[X_{\lambda}^0,W_0]$.  

\mn
\item{(v)} Let $Z_0\in\g^{0}_{\q}$. Then $Z_0^*=Ad_{x_0}Z_0 $. 
\mn
\item{(vi)}  Let $Y_{-\lambda}\in\g^{-\lambda}_\h$. Then $ Y_{-\lambda}^*=-i Ad_{x_0}[X_{\lambda}^0,Y_{-\lambda}]$.
\mn
\mn
\item{(vii)} Let  $X_{-\lambda} \in  \g^{- \lambda}_\q  $. Then
 $ X_{-\lambda}^*=Ad_{x_0}X_{-\lambda}-{1\over 2}Ad_{x_0}
[X_\lambda^0,[ X_\lambda^0,X_{-\lambda}]] $. 
\mn
\item{(viii)} Let $ Y_{-2\lambda}\in \g^{- 2\lambda}_\h$. Then $ Y_{-2\lambda}^*=-i Ad_{x_0}[X_{\lambda}^0,Y_{-2\lambda}]
+{i\over 6} Ad_{x_0}[X_{\lambda}^0,[X_{\lambda}^0,[X_{\lambda}^0,Y_{-2\lambda}]]] 
$. 

\end{lem}

\mn
\begin{proof}
All statements are obtained by combining  the formula
$Ad_{\exp iX}Y=\exp ad_{iX}Y$ with the bracket relations among root vectors. We omit the computations, which are long but straightforward.
\end{proof}

\bn
Fix $ \theta X^0_{ \lambda} \in \g^{-\lambda}_\q $ and denote by $( \g^{-\lambda}_\q)_0$ a complementary subspace to  $ \R  \theta X^0_{ \lambda}$ in $ \g^{-\lambda}_\q $.  
By (iii) of Lemma \ref{NONRED}  and Lemma \ref{FIELDS},  the  tangent space  to $S$ at $x_0$ is given~by 
 $TS_{x_0}=T_\C S_{x_0}\oplus NS_{x_0},$ 
where
\begin{equation}
\label{TANGENT}
 T_\C S_{x_0}=  Ad_{x_0}(\g^0_\q)^\C\oplus Ad_{x_0}( \g^\lambda_\q)^\C\oplus Ad_{x_0}( \g^{-\lambda}_\q)^\C_0, \qquad NS_{x_0}=\R Ad_{x_0}\theta X^0_{ \lambda}.\end{equation} 
Note that if $\Delta_\b$ is reduced, by $(i)$ of Lemma  \ref{RED}, one has $( \g^{-\lambda}_\q)_0=\{0\}$. 

\bn
\begin{remark}\label{ORTHOGONAL} There exists a  basis of $\g$  so that the above  decomposition of  $TS_{x_0}$ is orthogonal with respect to the Killing form $B$ of $\g^\C$.  If the restricted root system $\Delta_\b$ is non-reduced, one can construct  it starting from a  basis of $\g^\C/\s^\C$  consisting of root vectors  with respect to a maximally split Cartan subalgebra $\s$ of $\g$ extending $\b$. In the reduced case, this is immediate by $(i)$ of Lemma \ref{RED}.  \qed
\end{remark}

\bn 
{\bf The Levi form of the $G$-orbit of $x_0$.}  
The same arguments used in Section 4 of \cite{Ge1}  yield  the following formulas for the Levi form of  $S $ at $x_0$.  Let ${  Z},{  W} $  be  vectors in $T_\C S_{x_0}$. Then   
\begin{equation}\label{LEVI}
L_{x_0}({ Z},{  W})={1\over 2} [(*)^{-1}J{  W},{  Z}] -{i\over 2} [(*)^{-1}{  W},{  Z}]  \quad \hbox{mod~$(T_\C S_{x_0})^\C$ },
\end{equation}
where $(*)^{-1}J{  W}$ and $(*)^{-1}{  W}$ are arbitrary elements in the preimages  of $JW$ and $W$ by the map defined in (\ref{VECTORFIELDS}). 
 In the next lemma we compute the Levi form of $S$ at $x_0$. Fix  $F^0_{-\lambda}:=Ad_{x_0}\theta X^0_{ \lambda}$ as a  generator of $NS_{x_0}$.

\bn

\begin{lem} 
\label{LEVI2} 
\item{(i)}  Let $X_{-\lambda} \in (\g_\q^{-\lambda})_0$.  
Set   $F_{-\lambda}:= Ad_{x_0}X_{-\lambda} $. Then  
$$L_{x_0}(F_{-\lambda} ,F_{-\lambda} )= -{1\over 6} 
Ad_{x_0}[[X^0_\lambda,X_{-\lambda}],X_{-\lambda}]=p F^0_{-\lambda}   \hbox{\quad {\rm mod}~~$(T_\C S_{x_0})^\C$}
, \quad p\ge 0.$$

\item{(ii)}  Let   $Z_0\in \g^0_\q$.  Write  $Z_0=[X^0_\lambda,Y_{-\lambda}]$, for some $Y_{-\lambda}\in \g^{-\lambda}_\h$ (cf. Lemma \ref{NONRED}) and set  $F_0:= Ad_{x_0}Z_0$. Then 
$$L_{x_0}(F_0 ,F_0 )=-{1\over 2}Ad_{x_0}[Y_{-\lambda},Z_0]=nF^0_{-\lambda}  \hbox{\quad {\rm mod}~~$(T_\C S_{x_0})^\C$},\quad n\le 0.$$

\item{(iii)}  Let $X_\lambda\in\g^\lambda_\q$ and set $F_\lambda:= Ad_{x_0}X_\lambda$. Then 
$$ L_{x_0}(F_\lambda ,F_\lambda )= 0 .$$

\item{(iv)} Let $X_{-\lambda} \in (\g_\q^{-\lambda})_0$ and  $Z_0\in  \g^0_\q$. 
Set $F_{-\lambda}:= Ad_{x_0}X_{-\lambda}$ and $F_0:=  Ad_{x_0}Z_0$. Then 
$$L_{x_0}(F_{-\lambda},F_0 )=  0 .$$

\item{(v)} Let $X_{-\lambda} \in (\g_\q^{-\lambda})_0$ and  $X_\lambda\in\g^\lambda_\q$. 
Set  $F_{-\lambda}:= Ad_{x_0}X_{-\lambda} $ and $F_\lambda:=Ad_{x_0}X_{ \lambda} $. Then 
$$L_{x_0}(F_{-\lambda} ,F_{ \lambda} )=aF^0_{-\lambda},~~a\in\C .$$ 

\item{(vi)}  Let $X_\lambda\in\g^\lambda_\q$ and  $Z_0\in \g^0_\q$.  
Set   $F_\lambda:= Ad_{x_0}X_\lambda$ and $F_0:= Ad_{x_0}Z_0$. Then 
$$ L_{x_0}(F_0 ,F_\lambda )= 0 .$$

\end{lem}

\mn
\begin{proof}
By way of example, we prove the first two statements of the  Lemma. The remaining ones follow in a similar way, and the details are omitted.

\bn
(i) Let  $F_{-\lambda}= Ad_{x_0}X_{-\lambda} $. 
In order to compute the brackets  (\ref{LEVI}), we invert  the relations in Lemma \ref{FIELDS} and decompose the results in $\g^\C=Ad_{x_0}\h^\C\oplus Ad_{x_0}\q^\C$. Write $X_{-\lambda}=[X^0_\lambda,Y_{-2\lambda}]$, for some $Y_{-2\lambda}\in \g^{2\lambda}_\h$ (cf. Lemma \ref{NONRED}). Then
$$(*)^{-1}JF_{-\lambda}=-Y_{-2\lambda}+{1\over 6}ad_{X_\lambda^0}^2(Y_{-2\lambda})=$$
$$=-Ad_{x_0}Y_{-2\lambda}+iAd_{x_0}ad_{X_\lambda^0} (Y_{-2\lambda})+{1\over 2}Ad_{x_0}ad_{X_\lambda^0}^2(Y_{-2\lambda})-$$
$$-{i\over 6}Ad_{x_0}ad_{X_\lambda^0}^3(Y_{-2\lambda})+{3\over 8}Ad_{x_0}ad_{X_\lambda^0}^4(Y_{-2\lambda}).$$
and
$$(*)^{-1} F_{-\lambda}=ad_{X_\lambda^0} (Y_{-2\lambda})+{1\over 2} ad_{X_\lambda^0}^3(Y_{-2\lambda})= $$
$$=Ad_{x_0}ad_{X_\lambda^0} (Y_{-2\lambda}) - {i }Ad_{x_0}ad_{X_\lambda^0}^2(Y_{-2\lambda})+ {i\over 6}Ad_{x_0}ad_{X_\lambda^0}^4(Y_{-2\lambda}).$$
By formulas (\ref{LEVI}), we obtain
$$ L_{x_0}(F_{-\lambda} ,F_{-\lambda} )=-{1\over 6} Ad_{x_0}[ ad_{X_\lambda^0}^2(Y_{-2\lambda}), ad_{X_\lambda^0} (Y_{-2\lambda})]=$$
$$=-{1\over 6}Ad_{x_0}[[X^0_\lambda,X_{-\lambda}],X_{-\lambda}]    \hbox{\quad {\rm mod}~~$(T_\C S_{x_0})^\C$}.$$
To complete the proof of the statement, set $F_\lambda^0:=Ad_{x_0}X^0_\lambda$ and note that by Remark \ref{ORTHOGONAL}, the component $pF^0_{-\lambda}$ of the above brackets  in $NS_{x_0}$ is given by 
$$ B(L_{x_0}(F_{-\lambda} ,F_{-\lambda} ),F^0_\lambda)=p  B(F^0_{-\lambda},F^0_\lambda).$$
Since $B(F^0_{-\lambda},F^0_\lambda) =B( X^0_{-\lambda}, \theta X^0_{-\lambda})$ is negative,   the real number  $p$ has  the same sign as 

\begin{equation}
 \label{NEG} B( [[X^0_\lambda,X_{-\lambda}],X_{-\lambda}], X^0_\lambda)= - B( [X^0_\lambda,X_{-\lambda}], [X^0_\lambda,X_{-\lambda}]).
\end{equation}
By Lemmas \ref{RED} and  \ref{NONRED}, the brackets $[X^0_\lambda,X_{-\lambda}]$  lie in $\k$, so   
$B( [X^0_\lambda,X_{-\lambda}], [X^0_\lambda,X_{-\lambda}])$ is non-positive. It follows that $p\ge 0$, as claimed.

\bn
(ii) Write  $Z_0=[X^0_\lambda,Y_{-\lambda}]$, for some $Y_{-\lambda}\in \g^{-\lambda}_\h$ (cf. Lemma \ref{RED} and Lemma \ref{NONRED}). By computations similar to the above ones, we have
$$(*)^{-1}JF_0=-Y_{-\lambda}=-Ad_{x_0}Y_{-\lambda} +iAd_{x_0}ad_{X_\lambda^0}(Y_{-\lambda})+{1\over 2}Ad_{x_0}ad_{X_\lambda^0}^2(Y_{-\lambda}),$$
$$(*)^{-1} F_0= Z_0=Ad_{x_0}Z_0-iAd_{x_0} ad_{X_\lambda^0}( Z_0)  ,$$
and
 $$L_{x_0}(F_0 ,F_0 )=-\half Ad_{x_0}[Y_{-\lambda},[X^0_\lambda, Y_{-\lambda}]]=-\half Ad_{x_0}[Y_{-\lambda},Z_0] \quad {\rm mod}~(T_\C S_{x_0})^\C.  $$
To complete the proof of the statement, observe that  $n={{B(L(F_0,F_0), F^0_{\lambda})}\over {B(F^0_{-\lambda}, F^0_\lambda)}}$ has the same sign as 
 $$B( [Y_{-\lambda},[X^0_\lambda,Y_{-\lambda}]],X^0_{\lambda})= 
B( [X^0_\lambda,Y_{-\lambda}],[X^0_\lambda,Y_{-\lambda}]). $$
Since $[X^0_\lambda,Y_{-\lambda}]$ lies in $ \k$, the above expression  is non-positive  and  $n\le 0$, as claimed.    
  
 \end{proof}

\begin{prop}
\label{W5}
Let $S$ be the $G$-orbit of the point $x_0=\exp i X^0_\lambda\cdot z$. \pn
If the restricted root system $\Delta_\b$ is reduced,    
 then the Levi form of the  orbit $S$  is  definite provided that $\dim G/K>2$. It is identically zero if  $\dim G/K=2$.
   \pn
If the restricted root system $\Delta_\b$ is non-reduced,
 then the Levi form of the  orbit $S $  is indefinite.
\end{prop} 

\mn
\begin{proof} If the restricted root system $\Delta_\b$ is reduced, then only (ii),(iii) and (vi) of Lemma \ref{LEVI2} apply. 
By (ii) of Lemma \ref{LEVI2},  for every $F_0 \in Ad_{x_0} (\g^0_\q)^\C$, the real numbers $B(L(F_0,F_0), F^0_{\lambda})$ all have the same sign. In other words, the restriction of the Levi form to $Ad_{x_0} (\g^0_\q)^\C\subset T_\C S_{x_0}$ is either  definite or identically zero.  It is  identically zero when $ad_{X^0_{\lambda}}\colon \g^{-\lambda}_\h\to \g^0_\q$ is  the zero-map. This happens if and only if  $\g=\s\l(2,\R)$ and $\dim G/K=2$.

If the restricted root system $\Delta_\b$ is non-reduced, then $\dim G/K>2$. In this case,  the restriction of the Levi form to $Ad_{x_0} (\g^0_\q)^\C\subset T_\C S_{x_0}$ is definite. Moreover, by $(i)$ of Lemma \ref{LEVI2} and   $(iii)$ of Lemma \ref{NONRED},
  the restriction of the Levi form to $Ad_{x_0}(\g^0_\q)^\C\subset T_\C S_{x_0}$ is definite  with opposite sign.
As a result,   
the Levi form of $S$ is indefinite, as stated.  

\end{proof}

\bn 
{\bf The Levi form of the $G$-orbit $y_0$.}  
By the same methods, one can compute the tangent space and the Levi form of the $G$-orbit  $S$ of  the point $y_0=\exp  (-i X^0_\lambda)\cdot z$.  As we already remarked, the orbits $G\cdot x_0$ and $G\cdot y_0$ are distinct only when  the restricted root system of $\g$ is  reduced. So we  focus on  this case.
 For the tangent space to $S$ at $y_0$ one has $TS_{y_0}=T_\C S_{y_0}\oplus NS_{y_0}$, where 
$$T_\C S_{y_0}= Ad_{y_0}(\g^0_\q)^\C\oplus Ad_{y_0}( \g^\lambda_\q)^\C,\qquad  NS_{y_0}=\R Ad_{y_0}\theta X^0_{ \lambda}.$$
Fix $F^0_{-\lambda}:=Ad_{y_0}\theta X^0_{ \lambda}$ as a generator  of $NS_{y_0}$.
For the Levi form, one has the  following results.

\begin{lem}
\label{REDY}
\item{(i)}  Let   $Z_0\in \g^0_\q$.  Write  $Z_0=[X^0_\lambda,Y_{-\lambda}]$, for some $Y_{-\lambda}\in \g^{-\lambda}_\h$ (cf. Lemma \ref{NONRED}) and set  $F_0:= Ad_{y_0}Z_0$. Then 
$$L_{y_0}(F_0 ,F_0 )={1\over 2}Ad_{y_0}[Y_{-\lambda},Z_0]=sF^0_{-\lambda}  \hbox{\quad {\rm mod}~~$(T_\C S_{y_0})^\C$},\quad s\ge 0.$$

\item{(ii)}  Let $X_\lambda\in\g^\lambda_\q$ and set $F_\lambda:= Ad_{y_0}X_\lambda$. Then 
$$ L_{y_0}(F_\lambda ,F_\lambda )= 0 .$$

\item{(iii)}   Let   $Z_0\in \g^0_\q$ and $X_\lambda\in\g^\lambda_\q$.  Set $F_0:= Ad_{y_0}Z_0$ and $F_\lambda:= Ad_{y_0}X_\lambda$. Then 
$$ L_{y_0}(F_0 ,F_\lambda )= 0 .$$

\end{lem}

\begin{prop} 
\label{W23} Let $S$ be the $G$-orbit of the point $y_0$.\pn
If the restricted root system $\Delta_\b$ is reduced,    
 then the Levi form of the  orbit $S$  is  definite provided that $\dim G/K>2$. It is identically zero  if  $\dim G/K=2$.
\end{prop} 

\bn
\begin{remark}
\label{REMW23}
By Proposition \ref{W5} and Proposition \ref{W23}, if the restricted root system
$\Delta_\b$ is reduced, then   the Levi form of the orbits represented by $w_1$ and
$w_2$ in  diagram (\ref{DIAGRAM2}) is definite. This is consistent with the fact that
these orbits lie in the boundary of the Stein domains $D_1(0)$ and $D_2(0)$,
respectively (cf. Theorem \ref{STEIN}).  If $\dim G/K=2$, all orbits represented by
$w_1,\ldots,w_4$  in  diagram (\ref{DIAGRAM1}) are Levi flat.  We refer to Example
\ref{NONREDUCED2} for a classification of $G$-invariant Stein domains bounded by such
orbits. If the restricted root system $\Delta_\b$ is reduced, then   the Levi form of the orbit 
represented by $w_5$ in diagrams (\ref{DIAGRAM3})  and  (\ref{DIAGRAM4}) is indefinite.
As a consequence this orbit cannot lie in the boundary of a Stein $G$-invariant domain in $G^\C/K^\C$.
\end{remark}

\mn

\subsection { Non-Closed orbits with a CR singular orbit in their closure.}
We  consider now non-closed $G$-orbits containing in their closure the orbit of a point $z=g K^\C=\exp iAK^\C \in\mathcal A_0$, satisfying the condition $\alpha(A)=\pi/4$, with $\alpha$ a simple restricted root  (see (\ref{BASEPOINT2}) in Sect.~4.2).  
In this case the singular orbit $G\cdot z$  has dimension greater than $\dim G/K$.   
Recall  from Section 4.2 that the isotropy subgroup $H'$ of $z$ in $G$ is contained in 
 $G':=Z_G(g^4)$  and that  $G'/H'$ is diffeomorphic to a rank-one, pseudo-Riemannian symmetric space, totally real in $G^\C/K^\C$. 
 Let  $(\g'= \h'\oplus \q',\tau_z)$  be the associated symmetric algebra.   Non-closed $G$-orbits in $G^\C/K^\C$ containing $G\cdot z$ in their closure are in one-to-one correspondence with the nilpotent $Ad_{H'}$-orbits in $\q'$ (cf.~(\ref{LUNA}) and Remark \ref{ADJOINT}).  


\sn 

 To explicitly determine reference  points for such non-closed orbits  and their tangent spaces,  we decompose $\g$ and $\g'$ by an appropriate restricted root system.
Let  $\mathcal C'=\exp i\c'\cdot z$ be the standard Cartan subset with base point $z$. Recall that  $\c'=\R(X+\theta (X))$, where $X$ is a non-zero vector in $\g^{2\alpha}$. In particular,  $\c' $  is contained in $\g'$ (see (\ref{GIPRIMO})).
Define $\b'=\R(X-\theta(X))$. Then $\b'$ is a maximal  abelian subalgebra in $\h'\cap \p$  
and  the restricted root decompositions of  $\g$  with respect to  $\b'$ is given by
$$\g=Z_\k(\b')\oplus \b'\oplus \g^{\pm 2\lambda}\oplus  \g^{\pm \lambda}.$$
In order to determine how the above  root decomposition restricts to the subalgebra $\g'$, observe that  in general $\g'$ is not simple, but is the direct sum of a copy of   
$\s \o(m,1)$, with $m=\dim \g^{2\alpha}+1$ (even), and a compact subalgebra $\l$ entirely contained in $\h'$
\begin{equation}
\g'=\l\oplus \s \o(m,1),\qquad \h'=\l\oplus \s \o(m-1,1).
\end{equation}
Observe also that all real rank one Lie algebras with a non-reduced restricted root system are equal-rank. Hence the root system $\Delta$ of $\g^\C$  with respect to a maximally split Cartan subalgebra of $\g$ extending $\b'$ contains a real root. Since $\g^{2\lambda}$ is odd-dimensional (cf.~Table 4.0), the restriction of this real root to $\b'$ coincides the restricted root  $2\lambda$  (see \cite{Hl}, p.$\,$584).  
Since $\g'$ has a reduced restricted root system (cf.~(\ref{GIPRIMO})) and $\s\o(m,1) $, with $m$ even,   is equal-rank,  then $\g'\cap \g^{2\lambda}\not=\{0\}$.  It follows that the restricted root decomposition  of  $\g'$  with respect to  $\b'$ is given by
\begin{equation}\label{GIPRIMO2}
 \g'=Z_\k(\b')\oplus \b'\oplus \g^{\pm 2\lambda}.
 \end{equation}
Let
$$ \g'=\h'\oplus \q',  \qquad \hbox{with}\quad 
  \h'=\g^0_{\h'\cap \k }\oplus \g^{\pm 2\lambda}_{\h'}\oplus \b' \quad \hbox{and}\quad  \q'=\g^0_{\q'}\oplus \g^{\pm 2\lambda}_{\q'},$$ be the combined decomposition of $\g'$. Note that $\g'$ has real rank one as well.  Therefore $\g^0_{\q'}\subset \k$ and an analogue of Lemma \ref{RED} holds for $\g'$. Set $\g[\lambda]:=\g^\lambda\oplus \g^{-\lambda}$ and $\g[\alpha]:=\g^\alpha\oplus \g^{-\alpha}$ (here $\alpha$ is a restricted root  in $\Delta_\a$, as in Sect.~4).

\begin{lem} \label{RSPACES} The following facts hold:
\item{(i)}  $ \dim \g^{\pm 2\lambda}_{\q'}=1 $; 
\item{(ii)} $[\g_{\q'}^{2\lambda}, \g_{\h'}^{-2\lambda}]=\g^0_{\q'},\quad [\g_{\q'}^{-2\lambda}, \g_{\h'}^{ 2\lambda}]=\g^0_{\q'}$;
\item{(iii)} the decomposition $\g=\g'\oplus \g[\alpha]$
 is $ad_{\b'}$-stable. In particular $ \g[\alpha]=\g[\lambda]$. 
\end{lem}

\mn
\begin{proof}
\pn
Statement (i) follows from the fact that  $\dim \q'\cap \p=1$ (see the proof of $(ii)$ in Lemma  \ref{LOCAL2}), while $(ii)$ can be checked directly.
\pn
To prove statement $(iii)$, note that  $ad_{\b'}\g'\subset \g'$.  Moreover,  $ad_{\b'}(\g^{   \alpha}\oplus \g^{-  \alpha})\subset (\g^{   \alpha}\oplus \g^{-  \alpha})$.  By (\ref{GIPRIMO}) and (\ref{GIPRIMO2}) it follows that  the decomposition $\g=\g'\oplus \g[\alpha]$
 is $ad_{\b'}$-stable and $ \g[\alpha] =\g[\lambda]$. 
\end{proof}

\bn
{\bf Reference points for non-closed $G$-orbits.} 
Reference points for non-closed orbits containing $G\cdot z$ in their closures can be obtained  by applying the methods of the previous section  to the symmetric space $G'/H'$ (cf. (\ref{BASEPOINTS})). In this case take $X\in\g^{2\alpha}$, $\theta X $ and $A:=[\theta X,X]$, normalized so that $2\alpha(A)=2$. Then 
\begin{equation}\label{NRBASEPOINTS}
X_{2\lambda}^0={1\over 2}(A-(X +\theta X ) ), \qquad X_{-2\lambda}^0 ={1\over 2}(A+(X +\theta X ))
\end{equation}
are generators of 
$\g^{2\lambda}_{\q'}$ and  $ \g^{-2\lambda}_{\q'}$  respectively, and 
the points 
$$x_0=\exp iX^0_{2\lambda} \cdot z,\quad x_1=\exp iX^0_{-2\lambda}\cdot z,\quad   y_0=\exp( -iX^0_{2\lambda}) \cdot z,\quad y_1=\exp (-iX^0_{-2\lambda})\cdot z $$  lie on non-closed $G$-orbits in $G^\C/K^\C$  containing the singular orbit $G\cdot z$ in their closures. 
If the orbit diagram is of type (\ref{DIAGRAM3}) there are four such orbits,  represented by 
$w_3$, $w_2$, $w_1$ and $w_4$,  respectively. 
If the orbit diagram is of type (\ref{DIAGRAM4}) the points $x_0$ and $x_1$ lie on the same $G$-orbit, represented by $w_2$. Similarly, the  points $y_0$ and $y_1$   lie on the same $G$-orbit  represented by 
$w_1$. 
The $G$-orbits of $y_0$ and $y_1$ lie on the boundary of the Stein domain $D_1(0)$ (cf.Theorem~\ref{STEIN}).

\bn
{\bf The  tangent space to the $G$-orbit  of $x_0$.}  
Denote by $S$ the $G$-orbit of the point $x_0=\exp iX^0_{2\lambda}\cdot z$. In order to compute the 
  tangent space $TS_{x_0}$, observe that at  the point $z$ 
 \begin{equation}\label{INCLUSION}
T(G\cdot z)_z=\q' \oplus V_z, \quad \hbox{and}\quad T(G^\C/K^\C)_z=Ad_z\p^\C=(\q')^\C\oplus V_z, 
\end{equation}
where $\q'= T(G'\cdot z)_z$ and $ V_z=Ad_z \g[\alpha]_\p^\C$ is a complex subspace of $\g[\alpha]^\C$ (see  \cite{Ge1}, Prop. 3.2).
 It follows that 
\begin{equation}\label{INCLUSION2}
TS_{x_0}\subset Ad_{x_0}(\q')^\C\oplus Ad_{x_0}V_z .
\end{equation}
In order to determine generators for $TS_{x_0} $,
fix  a maximally split Cartan subalgebra $\mathfrak s$ of $\g$  extending $\b'$ and entirely contained in $\h' $ (one can check that in all cases under consideration $\h'$ has the same rank as $\g$ and such a Cartan subalgebra indeed exists).
Let 
 $$\g^\C=\s^\C \bigoplus_{\beta\in\Delta} \g^\beta$$ be the corresponding root decomposition of $\g^\C$ and let $\{ Z_\beta \}_{\beta\in\Delta}  $    be a complex basis of $\g^\C/\s^\C$ consisting of root vectors $Z_\beta\in  \g^\beta$.
Choose compatible orderings of $\Delta_{\b'}$ and $\Delta$ (i.e. a root $\beta\in \Delta$ is positive if its restriction to $\b'$ is). 
Fix $\lambda\in \Delta_{\b'}$(either a positive or a negative short restricted root) and denote by $\Delta_\lambda$ the set of roots in $\Delta$ which restricted to $\b'$ are equal to $\lambda$. The set  $\Delta_\lambda$   
consists of pairs of complex roots $$\beta_1,\bar\beta_1,\ldots,\beta_m,\bar\beta_m ,\qquad m=\half \dim \g^{\lambda},$$ all with the same real part, equal to $\lambda$.
For $\beta_i,\bar\beta_i\in\Delta_\lambda$,  choose root vectors $Z_{\beta_i}\in\g^{\beta_i}$ and $\sigma Z_{\beta_i}\in\g^{\bar\beta_i}$. Then the   vectors   
defined as
\begin{equation}\label{XY}
X_{\lambda}^i=Z_{\beta_i}+\sigma Z_{\beta_i}, \quad \hbox{and}\quad  Y_{\lambda}^i= -i(Z_{\beta_i}-\sigma Z_{\beta_i}),\quad i=1,\ldots,m,
\end{equation}
belong to $\g$ and  form a basis of the restricted root space $\g^\lambda.$
\begin{lem} \label{PRELIM}
The following facts hold:
\item{(i)} for all $i=1,\ldots,m$, one has $\tau_z Z_{\beta_i}=-Z_{\beta_i}$ and $ i\tau_zX_\lambda^i= Y_\lambda^i $;
\pn
\item{(ii)} for every $i=1,\ldots,m$,  the brackets $[X_{  \lambda}^i,i\tau_z X_{  \lambda}^i] $  lie in $\g^{2\lambda}_{\q'}$. For at least one index $i$, such brackets are non-zero;

\item{(iii)} for all $\, i,j=1,\ldots,m$, with $i\not= j$, the brackets $[X_{  \lambda}^i,i\tau_z X_{  \lambda}^j]$ have no components in $\g^{2\lambda}_{\q'}$. 

\end{lem}

\begin{proof} $(i)$ Since the Cartan subalgebra $\s$ lies in $\h'$, it is pointwise  fixed by $\tau_z$. As a consequence,  all root spaces $\g^\beta$, with $\beta\in\Delta $,  are $\tau_z$-stable. 
The inclusion $V_z\subset Ad_z\p^\C$ (see  (\ref{INCLUSION})) implies that $\tau_z Z_{\beta_i}=-Z_{\beta_i}$, for $i=1,\ldots,m$. Since $\sigma\tau_z=-\tau_z\sigma$ on  $V_z\subset \g[\alpha]^\C$,  one  has   $i\tau_z X^i_\lambda=Y^i_\lambda$, as desired.
\pn
$(ii)$ By  the definitions of $X^i_\lambda$ and $Y^i_\lambda$, one has 
$$[X_{  \lambda}^i,i\tau_z X_{  \lambda}^i] =[X_{  \lambda}^i,Y_{  \lambda}^i] =2i[Z_{\beta_i},\sigma Z_{\beta_i}]\in \g^{2\lambda}.$$
By $(i)$ and the fact that $\tau_z \sigma=-\sigma \tau_z$ on $\g[\lambda]^\C=\g[\alpha]^\C$, one also has 
$$\tau_z(2i [Z_{\beta_i},\sigma Z_{\beta_i}])=-2i [Z_{\beta_i},\sigma Z_{\beta_i}].$$
This implies that $[X^i_\lambda, i\tau_z X^i_\lambda]$ lies in $\g^{2\lambda}_{\q'}$, as claimed.
To prove the second part of the statement, consider the set  $\Delta_{2\lambda}$ consisting of the roots in $\Delta$ which restricted to $\b'$ coincide with $2\lambda$. Since $\Delta_{2\lambda}$ contains a real root in  $\Delta$  and such root is not simple (cf. Satake diagrams in [Hl], p.532), then  there exist $\beta, \bar\beta \in  \Delta_{ \lambda}$ such that  $\beta+\bar \beta=2\lambda$. 
This shows  that at least one of the  brackets $[X^i_\lambda, i\tau_z X^i_\lambda]$ has a non-zero component  in~ $\g^{2\lambda}$.

\pn
$(iii)$ Let $\beta_i, \beta_j$ be roots in $\Delta_\lambda$, with $\beta_j\not= \beta_i, \bar\beta_i$.  If either $  \beta_i+\beta_j$  or $ \beta_i+\bar\beta_j $ is a root in $\Delta$,  then it is a root in $\Delta_{2\lambda}$, with non-zero imaginary part. Since the root spaces relative to the real root in $\Delta_{2\lambda}$ are contained in $(\g^{2\lambda}_{\q'})^\C$ and $\dim (\g^{2\lambda}_{\q'})^\C=1$ (cf. Lemma \ref{NONRED}), then the root spaces relative to the remaining roots in  $\Delta_{2\lambda}$ are necessarily contained in $(\g^{2\lambda}_{\h'})^\C$. Hence the statement follows.
\end{proof}

\bn
For $\lambda\in \Delta^+_{\b'}$, fix  bases of $\g^{\lambda}$ and  $\g^{-\lambda}$ of the form
\begin{equation} \label{BASES} 
X_\lambda^1,i\tau_z X_\lambda^1,\ldots , X_\lambda^m,i\tau_z X_\lambda^m,\qquad X_{-\lambda}^1,i\tau_z X_{-\lambda}^1,\ldots , X_{-\lambda}^m, i\tau_z X_{-\lambda}^m,
\end{equation}
respectively. 
For $i,j=1,\ldots, m$, define
$$w_i:={1\over 2}Ad_{x_0}(X_{ \lambda}^i -\tau_z X_\lambda^i)\quad\hbox{and}\quad v_j:={1\over 2}Ad_{x_0}(X_{-\lambda}^j-\tau_zX_{-\lambda}^j).$$  
In the next lemma we compute  the images of the  vectors  
in (\ref{BASES}) 
under the map  ${}^*\colon \g\to TS_{x_0}$, defined in (\ref{VECTORFIELDS}). We omit the proof which consists of straightforward computations.

\begin{lem} \label{TANGENT2} The images of the vectors in (\ref{BASES}) under the map (\ref{VECTORFIELDS}) are given as follows.
\item{(i)} $( X_{ \lambda}^i)^*=w_i$,
\item{(ii)} $( i\tau_z X_{ \lambda}^i)^*=-iw_i$,
\item{(iii)} $( X_{-\lambda}^j)^*=v_j-iw',\quad \hbox{where }w'=Ad_{x_0} [X^0_{2\lambda},X_{-\lambda}^j]$,
\item{(iv)} $( i\tau_z X_{-\lambda}^j)^*=-iv_j-iw'',\quad \hbox{where }w''=Ad_{x_0} [X^0_{2\lambda},i\tau_z X_{-\lambda}^j]$.

\end{lem}

\bn
Denote  
by $W_{x_0}^+$  the complex subspace of  $W_{x_0}$ spanned by  the vectors $\{w_1,\ldots,w_m\}$ and by $W_{x_0}^-$  the one spanned by $\{v_1,\ldots,v_m\}$.
By (\ref{INCLUSION2}), the results of Section 9.1 applied to the symmetric space $G'/H'$ and  Lemma \ref{TANGENT2}, 
  the 
 tangent space to $S$ at $x_0$ is given by $TS_{x_0}=T_\C S_{x_0}\oplus NS_{x_0}$, where
\begin{equation}\label{TANGENTDECO}
T_\C S_{x_0}=T_\C(G'\cdot x_0)_{x_0}\oplus W_{x_0}^+\oplus W_{x_0}^- \qquad NS_{x_0}=\R Ad_{x_0}\theta X^0_{ 2\lambda}.
\end{equation}
 Fix $F^0_{-2\lambda}:=Ad_{x_0}\theta X^0_{2\lambda}$ as a generator of $NS_{x_0}$.

\begin{lem}\label{LEVIDECO} The following facts hold.
\item{(i)} The decomposition of $T_\C S_{x_0}$ given in (\ref{TANGENTDECO}) is orthogonal with respect to the Levi form.
\item{(ii)} Let $W\in W_{x_0}^+$. Then $L_{x_0}(W,W)= 0$.
\item{(iii)} Let $W\in W_{x_0}^-$. Then $L_{x_0}(W,W)=b F^0_{-2\lambda}$, with $b\ge 0$.
\item{(iv)} Let $Z\in T_\C(G'\cdot x_0)_{x_0}$. Then $L_{x_0}(Z,Z)=n F^0_{-2\lambda}$, with $n\le  0$.
\end{lem}

\mn
\begin{proof}  
(i) Let $Z\in T_\C(G'\cdot x_0)_{x_0}$ and $W\in W_{x_0}$.
In order to show that $L(Z,W)\equiv L(W,Z) \equiv 0 $, 
observe  that both $(*)^{-1}JZ,$ and $ (*)^{-1}Z$ belong to $\g'=\h'\oplus \q'$, and can be written as 
 $$(*)^{-1}JZ=Ad_{x_0}X_0 +Ad_{x_0}X_{ 2\lambda}+Ad_{x_0}X_{-2\lambda},\ \  (*)^{-1} Z= Ad_{x_0}Y_0+Ad_{x_0}Y_{  2\lambda} +Ad_{x_0}Y_{- 2\lambda}, $$
according to  the $ad_{\b'}$-root decomposition of  $\g'$ given in (\ref{GIPRIMO2}).
Similarly  by (\ref{INCLUSION}),  the vector $W\in W_{x_0}^+\oplus  W_{x_0}^+=Ad_{x_0}Ad_z\g[\lambda]_\p^\C$ can be written as
$$W=Ad_{x_0}Ad_zP_\lambda+i Ad_{x_0}Ad_zQ_\lambda,$$
where
$$Ad_zP_\lambda=U_\lambda+iV_{-\lambda}-\theta U_\lambda+i\theta V_{-\lambda},\quad \hbox{and}\quad Ad_zQ_\lambda=U_\lambda'+iV_{-\lambda}'-\theta U_\lambda'+i\theta V_{-\lambda}',$$
with $U_\lambda, U'_\lambda\in \g^\lambda$ and $V_{-\lambda}, V'_{-\lambda}\in \g^{-\lambda}$.
One can verify that none of the  brackets in $(\ref{LEVI})$ has a  component  in $Ad_{x_0}\g^{-2\lambda}_{\q'}$ and $L_{x_0}(Z,W)\equiv 0$, as required. 

\sn
Let $w_i\in W_{x_0}^+$ and $v_j\in W_{x_0}^-$. Then, modulo $(T_\C S_{x_0})^\C$, the Levi form is  given by
 $$2L_{x_0}(w_i,v_j)
\equiv -{1\over 2}Ad_{x_0}[ i\tau_zX_\lambda^i, (X_{-\lambda}^j- \tau_z X_{-\lambda}^j)]-
 {i\over 2}Ad_{x_0}[ X_\lambda^i, (X_{-\lambda}^j- \tau_z X_{-\lambda}^j)].$$
In particular, $L_{x_0}(w_i,v_j)=0$, for all $i,j=1,\ldots,m$. This concludes the proof of $(i)$. 
\pn
In the same way  one shows that   $L(w_i,w_j)=0$, for all $w_i,w_j\in W_{x_0}^+$, 
 proving $(ii)$.

 \sn
 $(iii)$ 
Similar calculations and $(iii)$ of   Lemma \ref{PRELIM} imply that  $ L_{x_0}(v_i,v_j) =0,$ 
for all $v_i, v_j\in W_{x_0}^-$, with $ i\not=j.$ 
When $i=j$,  one has 
$$L_{x_0}(v_i,v_i)=Ad_{x_0} [ X_{-\lambda}^i,  i\tau_zX_{-\lambda}^i]= 
Ad_{x_0} i[Z_{-\beta_i}, \sigma Z_{-\beta_i}]=
b_iF^0_{-2\lambda},\quad b_i\in\R.$$
In order to prove  that $b_i \ge 0$  observe that, by $(iii)$ of Lemma \ref{RSPACES} one can write $X_{-\lambda}^i=X^i_\alpha+X^i_{-\alpha}$, 
  for appropriate $X_\alpha^i \in \g^\alpha$ and $X_{-\alpha}^i \in \g^{-\alpha}$. Since $z=\exp iAK^\C$, with $A\in\a$ and $\alpha(A)=\pi/4$, one also  has     $i\tau_z X^i_{-\lambda} =\theta X^i_\alpha-\theta X^i_{-\alpha}$ and  
$$[ X_{-\lambda}^i,  i\tau_zX_{-\lambda}^i]=([X^i_\alpha, \theta X^i_\alpha]-
[X^i_{-\alpha}, \theta X^i_{-\alpha}])- ([X^i_{ \alpha}, \theta X^i_{-\alpha}]+[X^i_{-\alpha}, \theta X^i_{\alpha}])\in \a \oplus Z_\k(\a).$$
By  Lemma 5.1(i) in [Ge1], the first two terms of the above sum can be written as
$ [X^i_\alpha, \theta X^i_\alpha]=B( X^i_\alpha, \theta X^i_\alpha)A_\alpha$ and 
$ [\theta X^i_{-\alpha}, \theta (\theta X^i_{-\alpha})]=B( X^i_{-\alpha}, \theta X^i_{-\alpha})A_\alpha$, where $A_\alpha$ is an element in $\a$ satisfying the condition $\alpha(A_\alpha)>0$. By the normalization of the reference points  chosen in (\ref{NRBASEPOINTS}), one has   $\theta X^0_{2\lambda}=-X^0_{2\lambda}$. Hence 
$L_{x_0}(v_i,v_i)=b_iAd_{x_0} \theta X^0_{2\lambda}$, for some real number  $b_i\ge 0$, as claimed.
\pn 
 This concludes the  proof of  $(iii)$.

\sn
(iv) Recall that the symmetric space $G'/H'$ has a reduced restricted root system  and that the Lie algebra $\g'$ is given by (\ref{GIPRIMO2}). Then the Levi form on $T_\C(G'\cdot x_0)_{x_0}$ can be computed by the methods of Section 9.1. By (\ref{TANGENT}), one has
$$T_\C (G'\cdot x_0)_{x_0}= Ad_{x_0}({ \g }^{0}_{\q'})^\C \oplus Ad_{x_0}({\g}^{2\lambda}_{\q'})^\C\quad\hbox{and}\quad N (G'\cdot x_0)_{x_0}= \R Ad_{x_0}\theta X^0_{-2\lambda}.$$
Let $F_0\in Ad_{x_0}({\g }^{0}_{\q'})^\C $ and $F_{2\lambda}\in Ad_{x_0}({\g }^{2\lambda}_{\q'})^\C$. Then by Lemma \ref{LEVI2} one has  
$$L_{x_0}(F_{2\lambda},F_{2\lambda})=L_{x_0}(F_0,F_{2\lambda})=0,
\quad L_{x_0}(F_0,F_0)= nF^0_{-2\lambda},~~n\le 0.$$
\end{proof}

\bn
The next proposition is a direct consequence of Lemma \ref{PRELIM} and of Lemma \ref{LEVIDECO}.

\bn
\begin{prop} \label{NONREDW14} Let $S$ be the $G$-orbit of $x_0$. The Levi form of $S$ at $x_0$ is  indefinite  if $\g=\mathfrak s\mathfrak p (n,1)$ or $\g=\mathfrak f_4^*$. It is  definite if 
$\g=\mathfrak s\mathfrak u(n,1)$.
   
\end{prop}

\sn
\begin{proof} By  Lemma \ref{PRELIM} and $(iii)$ of Lemma \ref{LEVIDECO},  the Levi form $L_{x_0}$ is definite on $W_{x_0}^-$. 
If $\g=\mathfrak s\mathfrak u(n,1)$, then $\dim G'/H'=1$ and the Levi form is identically zero on  $T_\C (G'\cdot x_0)_{x_0}$. As a result, in this case the Levi form $L_{x_0}$ is  definite.

If $\g=\mathfrak s\mathfrak p (n,1)$ or $\g=\mathfrak f_4^*$, then $\dim G'/H'>2$ and the Levi form $L_{x_0}$ on $T_\C(G'\cdot x_0)_{x_0}$ is definite  of the opposite sign as on  $W_{x_0}^-$ (cf. Proposition \ref{W5} and Lemma \ref{LEVIDECO}). 
As a result, $L_{x_0}$  is indefinite, as claimed.
\end{proof}

\bn 
{\bf The Levi form of the $G$-orbit of $y_0$}. By the same methods, one can compute the tangent space and the Levi form of the $G$-orbit $S $
 of the point $y_0=\exp i (-X^0_{2\lambda})\cdot z.$
The 
 tangent space to $S$ at $y_0$ is given by $TS_{y_0}=T_\C S_{y_0}\oplus NS_{y_0}$, where
\begin{equation}\label{TANGENTDECO2}
T_\C S_{y_0}=T_\C(G'\cdot y_0)_{y_0}\oplus W_{y_0}^+\oplus W_{y_0}^- \qquad NS_{y_0}=\R Ad_{y_0}\theta X^0_{ 2\lambda}.
\end{equation}
Fix $F^0_{-2\lambda}:=Ad_{y_0}\theta X^0_{2\lambda}$ as a generator of $NS_{y_0}$.

\begin{lem}\label{LEVIDECO2} The following facts hold.
\item{(i)} The decomposition of $T_\C S_{y_0}$ given in (\ref{TANGENTDECO2}) is orthogonal with respect to the Levi~form.
\item{(ii)} Let $W\in W_{y_0}^+$. Then $L_{y_0}(W,W)\equiv 0$.
\item{(iii)} Let $W\in W_{y_0}^-$. Then $L_{y_0}(W,W)=b F^0_{-2\lambda}$, with $b\ge 0$.
\item{(iv)} Let $Z\in T_\C(G'\cdot y_0)_{y_0}$. Then $L_{y_0}(Z,Z)=p F^0_{-2\lambda}$, with $p\ge  0$;

\end{lem}

\begin{proof} The proof follows the same pattern as the proof of Lemma \ref{LEVIDECO}. One can check that the Levi form  is not identically zero on $W^-_{y_0}$ and has  the same signature as on $W^-_{x_0}$. Part (iv) follows from Lemma \ref{REDY}. 
\end{proof}

\bn
\begin{prop}
\label{W2}
Let $S$ be the $G$-orbit of $y_0$. The Levi form of $S$ at $y_0$ is  definite. \end{prop}

\mn
\begin{proof}  The proposition follows  from Lemma \ref{LEVIDECO2} and the fact that the Levi form $L_{y_0}$ on $W_{y_0}^-$ is not identically zero. 
  
\end{proof} 

\bn
\begin{remark}
\label{REMW2}
If the restricted root system $\Delta_\b$ is non-reduced, then  Proposition \ref {W2}
says that   the Levi form of the orbits represented by $w_1$ and $w_4$ in  diagram
(\ref{DIAGRAM3}) and in diagram (\ref{DIAGRAM4}) is definite. This is consistent with
the fact that these orbits lie in the boundary of the Stein domain  $D_1(0)$ (cf.
Theorem \ref{STEIN}).  When $\g=\s\u(n,1)$, by   Proposition \ref{NONREDW14}, the same is true for the Levi
form of the orbits represented by $w_2$ and $w_3$ in  diagram (\ref{DIAGRAM3}). We refer to  Example \ref{NONREDUCED2} for a classification of the
$G$-invariant Stein domains in $G^\C/K^\C$ bounded by these orbits.   Proposition
\ref{NONREDW14} also says that the Levi form of the orbit  represented by  $w_2$ in  
diagram (\ref{DIAGRAM4}) is indefinite.
Hence  this orbit cannot lie in the boundary of a Stein $G$-invariant domain in $G^\C/K^\C$.
\end{remark}

\bn

\end{document}